\documentclass[11pt,a4paper,fleqn]{article}
\usepackage{amsfonts,amsmath}
\usepackage{latexsym}
\usepackage{amssymb}
\usepackage{euscript}
\usepackage{graphicx}
\usepackage{a4wide}

\newtheorem{prop}{Proposition}[section]
\newtheorem{coro}[prop]{Corollary}
\newtheorem{lemma}[prop]{Lemma}
\newtheorem{exam}[prop]{Example}
\newtheorem{rem}[prop]{Remark}
\newtheorem{theorem}[prop]{Theorem}
\newtheorem{defi}[prop]{Definition}

\newcommand{\dis}{\displaystyle }

\newcommand{\noi}{\noindent}

\newenvironment{proof}[1]{\begin{trivlist}\item {\it
\bf Proof.}\quad} {\qed\end{trivlist}}
\newenvironment{prooff}[1]{\begin{trivlist}\item {\it
\bf Proof}\quad} {\qed\end{trivlist}}

\newcommand{\qed}{\nopagebreak\hspace*{\fill}
{\vrule width6pt height6ptdepth0pt}\par}

\begin{document}


\title{\bf LIMITING LAWS ASSOCIATED WITH BROWNIAN MOTION PERTURBED
BY NORMALIZED EXPONENTIAL WEIGHTS, I}
\author{{\small{\text{\bf Bernard ROYNETTE}$^{(1)}$}},\
{\small {\text{\bf  Pierre VALLOIS}$^{(1)}$}} {\small and} {\small
{\text{\bf Marc YOR }$^{(2)}$}}}



\maketitle {\small

\noindent (1)\,\, Universit\'e Henri Poincar\'e, Institut de
Math\'ematiques Elie Cartan, B.P. 239, F-54506 Vand\oe uvre-l\`es-Nancy Cedex\\

\noi (2)\,\, Laboratoire de Probabilit\'{e}s et
 Mod\`{e}les Al\'{e}atoires, Universit\'{e}s Paris VI et VII -  4, Place Jussieu
 - Case 188 -
 F-75252 Paris Cedex 05\\



\vskip 40 pt \noi {\bf Abstract.} Let $ (B_t;t\geq 0)$
 be a one- dimensional Brownian motion, with local time process
$(L^x_t; t\geq 0,x\in\mathbb{R})$. We determine the rate of decay
of $\displaystyle{
    Z^V_t(x):=E_x\Big[\exp\Big\{-\frac{1}{2}\int_\mathbb{R}L^y_t
    V(dy)\Big\}\Big], \  t\geq 0, x\in\mathbb{R}}$
as $t$ goes to infinity,  where $V(dy)$ is a positive Radon
measure on $\mathbb{R}$. If
$\displaystyle{\int_\mathbb{R}(1+|y|)V(dy)<\infty}$, we prove that
$Z^V_t(x)_{\stackrel{ \sim}{t \rightarrow \infty}}
\varphi_V(x)t^{-1/2}$, where the function $\varphi_V$ solves the
Sturm-Liouville equation $(\varphi_V)"(dx)=\varphi_V(x)V(dx)$,
with some boundary conditions. If
$\displaystyle{\int_{-\infty}^0(1+|y|)V(dy)<\infty}$ and $V(dy)$
is "large" at $+\infty$, the asymptotics of $Z^V_t(x)$ is the same
as previously. When $V(dy)$ is "large" at $\pm \infty,\ Z^V_t(x)$
is equivalent to $ke^{-\gamma_0 t}, t\rightarrow \infty$. If
$V(dy)=[\lambda/(\theta+y^2)]dy, \ \lambda, \theta >0 $, the rate
of decay is polynomial : $Z^V_t(x) _{\stackrel{ \sim}{t
\rightarrow \infty}} k \varphi_V(x/\sqrt{\theta})t^{-n}$, with
$n=(1+\sqrt{1+4\lambda})/4$. Taking $V(dy)=[1/(1+|y|^\alpha)]dy,
0<\alpha <2$ we only obtain a logarithmic equivalent :
$\displaystyle{\ln (Z^V_t(x)) _{\stackrel{ \sim}{t \rightarrow
\infty}} -k t^{-\frac{\alpha -2}{\alpha +2}}}$.

\noi Let $Q _{x,t}$ be the probability measure defined on the
canonical space  $\Omega={\cal C}([0,+\infty[)$, by :
$$\displaystyle{Q _{x,t}=\frac{1}{Z^V_t(x)} \exp \Big
\{-\frac{1}{2} \int _ \mathbb{R} L^y_t V(dy) \Big \} \ W_x} ,$$
\noi where $W_x$ denotes the Wiener measure. We prove that $Q
_{x,t}$ converges as $t\rightarrow\infty$ to $Q_x$ and $Q_x$ is
the law of the diffusion process $X _t^x$, solution of the
 stochastic differential equation :
 $$ X_ t=x+ B _t+\int _0
^t \frac{\varphi_V'}{\varphi_V}  (X_ s )ds; \ t \geq 0.$$


\vskip 40pt


\noi {\bf Key words and phrases} : Normalized exponential weights,
penalization,  Sturm-Liouville equation, Ray-Knight's theorems,
rate of convergence, Bessel processes.
\smallskip



\noi {\bf AMS 2000 subject classifications} : 60 F 10, 60 F 17, 60
G 44, 60 J 25, 60 J 35, 60 J 55, 60 J 57, 60 J 60, 60 J 65.


\vskip 40 pt \section{ Foreword and perspectives}

\noi This paper is the first in a series of four related papers,
numbered I to IV, a sketchy description of which may be of
interest to the reader.

 \noi Stimulated by the results obtained in I (see the above
abstract), we consided, in II, some asymptotic problems obtained
from weighting the Wiener measure with a function of the maximum,
or minimum, or local time up to time $t$, and letting
$t\rightarrow \infty$. The limit laws exist in some generality;
they are not the distribution of a Markov process $(X_t)$, but
rather the two dimensional process $\dis \big(
X_t,S_t:=\sup_{0\leq s \leq t}X_s ; t\geq 0\big)$ are Markovian.
We then say that $(X_t)$ is max-Markovian.

\noi In III, we study the existence and characterization of the
limit laws for Brownian bridges,  on the time interval $[0,t]$, as
$t\rightarrow \infty$; weighted again by a function of the
maximum, or minimum, or local time up to time $t$.

\noi In IV, we study the variants of the Pitman and Ray-Knight
theorems for the max-Markovian process obtained in the previous
papers.


\vskip 40 pt
\section{Introduction}\label{intro}

{\bf 1.1} Consider a general nice Markov process $((X_t)_{t \geq
0},({\cal F}_t)_{t \geq 0},(P_x)_{x\in E})$ taking values in
$(E,{\cal E})$, with extended generator $L$, i.e. : $\varphi \in
{\cal D}(L)$ iff : $M^\varphi _t:=\varphi (X_t)-\varphi (X_0)-\int
_0^t L\varphi (X_s)ds,\ t\geq 0, $ is a martingale, for some
function $L\varphi$, and its operator "carr\'{e} du champ" $\Gamma
(\varphi ,\psi )$, defined via :
\begin{equation}\label{001Aintro}
\frac{d<M^\varphi , M^\psi >_t}{dt}=\Gamma (\varphi ,\psi )(X_t).
\quad (\varphi ,\psi \in {\cal D}(L)).
\end{equation}
Indeed, by "nice Markov process", we mean in particular that
${\cal D}(L)$ is an algebra; hence (cf \cite{DelMey}) it follows
that :
\begin{equation}\label{002Aintro}
    \Gamma (\varphi ,\psi )=L(\varphi \psi)-\varphi L\psi -\psi
    L\varphi,\quad  \varphi ,\psi \in {\cal
D}(L).
\end{equation}

 {\bf 1.2} Associated with the family $(P_x)_{x \in E}$, we shall
consider two other families of probabilities constructed from the
$(P_x)_{x \in E}$.

\noi a) $\underline{\mbox{The first family} (Q_{x,t}^V)}$. For a
potential function $V : E \mapsto \mathbb{R}$, such that :
\begin{equation}\label{1Aintro}
    Z^V_t (x):= E_x[\exp \{ -\frac{1}{2}\int _0 ^t V(X_s)ds
    \} ]<\infty , \quad \mbox{for every }\ x\in E,
\end{equation}
we define the family of (normalized) probabilities :
\begin{equation}\label{2Aintro}
    Q_{x,t}^V = \frac{\exp \{ -\frac{1}{2}\int _0 ^t V(X_s)ds
    \} }{Z^V_t (x)}\ P_{x | {\cal F}_t}.
\end{equation}
Note that, in general, these laws are not coherent, i.e. for
$s<t$,
\begin{equation}\label{3Aintro}
Q_{x,t| {\cal F}_s}^V \not =Q_{x,s}^V .
\end{equation}
b) $\underline{\mbox{The second family} (P_x^\varphi)}$. Let
$\varphi >0$ be an element of ${\cal D}(L)$, then it is well-known
that :
\begin{equation}\label{4Aintro}
    \frac{\varphi (X_t)}{\varphi (x)}\exp \{-\int _0^t
    \frac{L\varphi}{\varphi}(X_s)ds\} \quad \mbox{ is a }
    ({\cal F}_t)\mbox{-local martingale w.r. to any } P_x,
    \mbox{for} \ x\in E .
\end{equation}
\noi We suppose that it is a martingale, so that we can define a
second (Markov) family of probabilities :
\begin{equation}\label{5Aintro}
P_{x|{\cal F}_t}^\varphi :=\frac{\varphi (X_t)}{\varphi (x)}\exp
\{-\int _0^t
    \frac{L\varphi}{\varphi}(X_s)ds\} \ P_{x|{\cal F}_t}.
\end{equation}
As a converse to (\ref{4Aintro}), if $\varphi \in {\cal D}(L)$,
and
 there exists $g$ such that
$\displaystyle \Big(\varphi(X_t)\exp \{- \int _0^t g(X_u)du\} ;
t\geq 0\Big)$ is a $((P_x),({\cal F}_t))$ local martingale, then
$L\varphi = g\varphi$.

 c) Let $V_\varphi$ be the potential
function associated with $\varphi \in {\cal D}(L), \varphi >0$ :
\begin{equation}\label{6Aintro}
    \frac{1}{2}V_\varphi = \frac{L\varphi}{\varphi} .
\end{equation}
It is clear that  the following identities hold :

\begin{equation}\label{10Aintro}
   V_{c\varphi}=V_\varphi , \qquad Q_{x,t}^{V+c}= Q_{x,t}^{V},
   \qquad \mbox{ for any } \quad  c>0,
\end{equation}
and that the two probabilities introduced in a) and b) are related
via :
\begin{equation}\label{7Aintro}
P_{x|{\cal F}_t}^\varphi :=\frac{Z_t^{V_\varphi}(x)}{\varphi
(x)}\varphi (X_t) Q_{x,t}^{V_\varphi}.
\end{equation}
It goes back to \cite{kunita} that, under $(P_x ^\varphi )$,
$(X_t)$ is a Markov process with extended infinitesimal generator

\begin{equation}\label{8Aintro}
L^\varphi =L+\frac{1}{\varphi}\Gamma (\varphi , \cdot) .
\end{equation}
There exists a simple relation between the Markovian laws
$(P_x^\varphi)$ (or, rather the associated semi-group
$(T_t^\varphi )$) and $Z^{V_\varphi}$, namely  :
\begin{equation}\label{9Aintro}
\varphi(x)T^{\varphi}_t(\frac{1}{\varphi})(x)= Z^{V_\varphi}_t(x)
= E_x[\exp \{-\frac{1}{2}\int _0^t V_\varphi (X_s)ds \}],
\end{equation}
which follows from (\ref{5Aintro}).

{\bf 1.3} We are interested in finding some conditions on $V$
which ensure the weak convergence, as $t\rightarrow \infty$, of
$Q_{x,t}^V$. We have two possibilities :

- to a given $\varphi >0$ in ${\cal D}(L)$, we may associate the
potential function $V_\varphi$ defined by (\ref{6Aintro}).

- conversely, starting from a potential function $V$, we may look
for the  solutions $\varphi$ of the Poisson equation (which is the
Sturm-Liouville equation in the Brownian case) :
\begin{equation}\label{11Aintro}
  \frac{L\varphi}{\varphi} = \frac{1}{2}V .
\end{equation}
We shall see later that a particular function $\varphi_V$  plays a
central role in our  discussion of  the convergence.

{\bf 1.4} $\underline{\mbox{A meta-theorem and its "proof"}}$.

\noindent The following statement shall be rigorously proved under
various hypotheses all throughout our paper. It may be used as a
guideline for the reader, and shall be refered to as the generic
theorem.

\begin{theorem}\label{theogene} Assume that, for some $k \geq
0$, one has :
\begin{equation}\label{12Aintro}
    \lim _{t\rightarrow \infty} \Big( t^k Z^V_t(x)\Big)=
\lim _{t\rightarrow \infty} \Big( t^k E_x[\exp \{-\frac{1}{2}\int
_0^t V (X_s)ds \}]\Big)=\varphi _V(x) .
\end{equation}
\noindent Then, $\varphi_V$ is a solution of (\ref{11Aintro}),
$\displaystyle \Big\{\frac{\varphi _V(X_t)}{\varphi _V (x)}\exp
\{-\frac{1}{2}\int _0^t V(X_u)du\} ; t\geq 0\Big\}$ is a
$((P_x),({\cal F}_t))$ martingale, and for any $\Lambda _s \in
{\cal F}_s$,
\begin{equation}\label{13Aintro}
\lim _{t\rightarrow \infty}  Q_{x,t}^V (\Lambda _s )=
P_x^{\varphi_V}(\Lambda _s ).
\end{equation}
\end{theorem}


\begin{prooff} \ {\bf of Theorem \ref{theogene}}.  Let $s>0$ be fixed,
$ \Lambda _s \in {\cal F}_s$ and $t>s$. We start from
(\ref{2Aintro}), and we write :
\begin{equation}\label{14Aintro}
Q_{x,t}^V (\Lambda _s )=\frac{E_x \Big[ 1_{\Lambda _s } \exp
\{-\frac{1}{2}\int _0^s V(X_u)du\}E_{X_s}\Big[ \exp
\{-\frac{1}{2}\int _0^{t-s} V(X_h)dh\}\Big]\Big]} {E_x[\exp
\{-\frac{1}{2} \int _0^t V(X_u)du\}]} ,
\end{equation}
and we multiply both the numerator and denominator  by $t^k$. The
result will follow after some justification for the  passage to
the limit inside the expectation (for the numerator).

\end{prooff}

\noi The "proof" shows that the normalization function $t \mapsto
t^k$ in (\ref{12Aintro}) may be replaced by any positive and
non-decreasing function $\lambda$ such that $\displaystyle \lim
_{t\rightarrow \infty}
\Big(\frac{\lambda(t+s)}{\lambda(t)}\Big)=1$, and also admits some
simple extension for $\lambda (u)=ce^{au}$, say.

 {\bf 1.5} $\underline{ \mbox{Back to the Brownian framework}}$.

\noindent Here, $E=\mathbb{R},L=\frac{1}{2}\frac{d^2}{dx^2},
\Gamma (f,g)(x)=f'(x)g'(x), P_x=W_x$ is the Wiener measure, and we
write $B_t$ instead of $X_t$ since, in this case, $\big(
(B_t;t\geq 0) ; (P_x;x\in\mathbb{R})\big)$ is a one-dimensional
Brownian motion.

\noindent Hence, we are interested in the Sturm-Liouville equation
:
\begin{equation}\label{15Aintro}
    \varphi"=V\varphi ,
\end{equation}
and we have $\displaystyle L^\varphi =\frac{1}{2}\frac{d^2}{dx^2}+
\Big(\frac{d}{dx}(\log \varphi)\Big)\frac{d}{dx}$,
$\mu(dx)=\varphi^2(x)dx$ is invariant under $(T^\varphi_t)$,
 the semigroup associated with $(P^\varphi_x)$, and :
 $$<T^\varphi_tf ,g>_\mu =<f, T^\varphi_t g >_\mu ,$$
 $$
 <L^\varphi f,g>_\mu=<f,L^\varphi g>_\mu=-\frac{1}{2}<f' , g'>_\mu
 .$$

 \noindent
 Since  Brownian motion admits a (bi-continuous)
family of local times $(L^x_t; t\geq 0,x\in\mathbb{R})$, we may
define the normalization factor $Z_t^V(x)$ (cf (\ref{1Aintro}))
when $V$ is a non-negative Radon measure on $\mathbb{R}$, in the
following way :
\begin{equation}\label{16Aintro}
Z_t^V(x)=E_x\Big[\exp\Big\{-\frac{1}{2}\int_\mathbb{R}L^y_t
    V(dy)\Big\}\Big].
\end{equation}
Abusing  notation, we use the same letter $V$, $ \underline{\mbox{
whenever} \ V  \mbox{ stands for a function or
 a Radon measure}}$.

 1) We investigate the $\underline{\mbox{integrable case}}$, i.e.
 when $V(dx)$ satisfies :
\begin{equation}\label{4intro}
    \int_\mathbb{R}(1+|x|)V(dx)<\infty.
\end{equation}
\noi In Theorem \ref{theoC}, we prove  that $\sqrt {t} Z_t^V(x)$
converges as $t\rightarrow \infty$ to a real number denoted
$\varphi_V(x)$. Moreover $\varphi_V$ is a convex function which
takes its values in $]0,\infty[$ and is the unique solution to the
Sturm-Liouville equation (\ref{15Aintro}), with boundary
conditions :
\begin{equation}\label{6intro}
\lim  _ {x\rightarrow +\infty} \varphi_V'(x)= - \lim  _
{x\rightarrow -\infty} \varphi_V'(x) =\sqrt { \frac{2}{\pi} }.
\end {equation}
\noi This leads us to relax the assumption (\ref{4intro}). We
shall discuss whether  $V(dx)$ is "large" at infinity or not.

2) Let us  examine the  first case. Suppose for simplicity that
$V$ is a function.

\noi We begin with an intermediate case. We say that $V$ is
asymmetric if :
\begin{equation}\label{9intro}
    \int_{-\infty}^0(1+|x|)V(x)dx<\infty,
\end{equation}
and
\begin{equation}\label{10intro}
    \liminf_{x\rightarrow\infty}\big(x^{2\alpha}V(x)\big)>0,
    \quad \mbox{ for some } \alpha <1.
\end{equation}
\noi Then (cf Theorem \ref{thbase1}), the rate of decay of
$Z_t^V(x)$ is unchanged : $\sqrt {t} Z_t^V(x)$ converges, as
$t\rightarrow \infty$, to $\varphi_V (x)$ : roughly speaking,
(\ref{9intro}) "dominates" (\ref{10intro}). Again the function
$\varphi_V$ solves  the Sturm-Liouville equation (\ref{15Aintro}),
but with the new boundary conditions :
\begin{equation}\label{11intro}
    \lim  _ {x\rightarrow -\infty} \varphi_V'(x)= - \sqrt { \frac{2}{\pi} }\ ;
    \qquad \lim  _ {x\rightarrow+\infty} \varphi_V(x) =0.
\end{equation}

3) Let us now investigate the case where $V(dx)$ is
$\underline{\mbox{ small at infinity}}$ but does not satisfy
(\ref{4intro}). We restrict ourselves to two  examples:
\begin{equation}\label{16intro}
V(x)= \frac{\lambda}{\theta+x^2} , \quad \mbox{where}\quad
\lambda>0,\  \theta\geq 0.
\end{equation}
\noi   and
\begin{equation}\label{15intro}
    V(x)=\frac{\lambda}{1+|x|^\alpha}, \quad
    \mbox{where}\quad \lambda>0, \ 0<\alpha<2.
\end{equation}
\noi Suppose $V$ is given by (\ref{16intro}). If $\theta =0$ then
(cf Theorem \ref{theo1Bes}) :
\begin{equation}\label{20intro}
    \lim _{t\rightarrow\infty} \Big ( t^{n}E_x\Big [ \exp
\Big \{ -\frac{\lambda}{2} \int _0 ^t\frac{ds}{B_s^2} \Big \} \Big
] \Big )=x^{2n} \frac{1}{2^{n}}\frac{\Gamma (n +\frac{1}{2})}
{\Gamma (2n +\frac{1}{2})} ,
\end{equation}
\noi where $\displaystyle n=\frac{1+\sqrt{1+4\lambda}}{4}$.

\noi When $\theta>0$ the result looks like the previous one. Let
$\varphi_V$ be the unique smooth function defined on
$[0,+\infty[$, solution of $\displaystyle {\varphi"(x)=
\lambda\frac{1}{1+x^2}\varphi(x); \quad x>0}$, such that: $
    \varphi_V(x) \sim x^{2n} , x\rightarrow +\infty$. In
    Theorem \ref{theo2Bes} we give  the explicit form of $\varphi_V$ and
we prove  :
\begin{equation}\label{21intro}
    \lim _{t\rightarrow\infty} \Big ( t^{n}E_x\Big [ \exp
\Big \{ -\frac{\lambda}{2} \int _0 ^t\frac{ds}{\theta+B_s^2} \Big
\} \Big ] \Big )=\theta^{n}\varphi_V(x/\sqrt{\theta})
\frac{1}{2^{n}}\frac{\Gamma (\mu +n +1)} {\Gamma (\mu +2n +1)},
\end{equation}
\noindent where $\mu=-1/2$.

\noindent
 We observe that if formally we take the limit $\theta\rightarrow 0$ in
(\ref{21intro})we recover (\ref{20intro}).

\noi Note that $(B_s^2;s\geq0)$ is a squared Bessel process with
dimension $1$, which led us to generalize  the asymptotic results
(\ref{20intro}) and (\ref{21intro}) to Bessel processes with any
positive dimension (see Theorem \ref{theo1Bes} for $\theta=0$ and
Theorem \ref{theo2Bes} when $\theta >0$ and $0<\lambda < 8\mu
^2+6\mu +1$).

\noi  Let us deal with the second case : $V$ is given by
(\ref{15intro}). We only obtain in Theorem \ref{theo1GD}, a
logarithmic equivalent for $Z_t^V(x)$ :
\begin{equation}\label{17intro}
    \lim_{t\rightarrow\infty} \Big( t^{\frac{\alpha-2}{\alpha+2}}
\ln \big (Z_t^V(x)\big ) \Big)=-\frac{1}{2}\Theta_0(\lambda),
\end{equation}
\noi where
\begin{equation}\label{18intro}
\Theta_0(\lambda)=\inf _{\psi \in {\cal C}_0}\Big\{ \int
_0^1\dot{\psi}^2(s)ds +\lambda \int _0^1\frac{ds}{|\psi
(s)|^\alpha} \Big \},
\end{equation}
\noi belongs to $]0,+\infty[$, and ${\cal C}_0$ is the set of
continuous functions $f :[0,1] \rightarrow \mathbb{R}$ vanishing
at $0$.

 4) If $V(dx)$ is large at $\pm \infty$, the asymptotic
behaviour of $Z_x^V(t)$ is drastically different. This case was
actually considered by Kac \cite{Kac} and Titchmarsh \cite{Titch}.
 More precisely suppose that  $V$ is an even function,
non-decreasing on $[0,+\infty[$ and converging to a finite limit
at infinity. Then we prove in Theorem \ref{thunil} that there
exists $\gamma_0>0$ such that
\begin{equation}\label{12intro}
    \lim_{x\rightarrow\infty}\Big(
    e^{\gamma_0t/2}Z_t^V(x)\Big)= \kappa \psi_V(x),
\end{equation}
\noi where $\kappa >0$ and $\psi_V$ is the positive  solution to
$\psi"(x)=\psi(x)(V(x)-\gamma_0),$ converging to $0$ at infinity
and such that $\psi_V'(0)=0$.

 5) Theorem \ref{theogene} tells us  that as soon as we obtain an explicit
behaviour of $Z_t^V(x)$ as $t$ runs to infinity, we may proceed
further to define new probability measures.

\noi In the Brownian setting, the probability $Q^V _{x,t}$ is
defined on ${\cal F}_t$ by :
\begin{equation}\label{23intro}
    Q _{x,t}(\Lambda  _t) = \frac{ E _x \Big [1 _{\Lambda  _t} \exp
\Big \{-\frac{1}{2} \int _ 0^t V(B_ v )dv \Big \} \Big ]} { E _x
\Big [ \exp \Big \{-\frac{1}{2} \int _ 0 ^t V(B_ v )dv \Big \}
\Big ] } ,\  t >0, \Lambda  _t \in {\cal F} _t,
\end{equation}
if $V$ is a function. In the case where $V(dy)$ is a Radon
measure, we have :
\begin{equation}\label{22intro}
    Q _{x,t}(\Lambda  _t) = \frac{ E _x \Big [1 _{\Lambda  _t} \exp
\Big \{-\frac{1}{2} \int _ \mathbb{R} L^y_t V(dy) \Big \} \Big ]}
{ E _x \Big [ \exp \Big \{-\frac{1}{2} \int _ \mathbb{R} L^y_t
V(dy) \Big \} \Big ] },\  t >0, \Lambda  _t \in {\cal F} _t.
\end{equation}

\noi To describe the probability measure $P_x^{\varphi_V}$, we
introduce $(X_ t^x;t\geq 0)$, the solution of the stochastic
differential equation :

\begin{equation} \label{30intro}
X _t=x+B _ t+\int  _0 ^t {\varphi'_V \over \varphi_V}(X _s)ds, \ t
\geq 0
\end{equation}

\noi The law of $(X_ t^x;t\geq 0)$ is $P_x^{\varphi_V}$.

6) Let us briefly detail the organization of the paper. Section
\ref{prel} presents some preliminaries.

\noi In section \ref{Csmall} we start with a function $\varphi$
which is  locally the difference of two convex functions, hence
$\varphi"$ is a Radon measure.  We take : $V_\varphi :=\varphi
"/\varphi$. Notice that  the sign of $V_\varphi$ is not constant.
With some additional assumptions on $\varphi$ such as : $\varphi$
"small at infinity" we prove in Theorem \ref{theoC} that
$Z_t^{V_\varphi}(x)$ converges, as $t\rightarrow\infty$ to $C
\varphi(x)$, where $C$ is a suitable constant.

\noi Section \ref{integrable case} is devoted to the proof of the
generic Theorem in the integrable case, namely when $V(dy)$
satisfies $\displaystyle{\int_\mathbb{R} (1+|y|)V(dy)<\infty}$. We
develop an analytic approach,  and two other ones based
respectively on the Ray-Knight theorem and excursions.

\noi The asymmetric case (i.e. when $V(x)$ satisfies
(\ref{9intro}) and (\ref{10intro}))  is discussed in section
\ref{casunil}.

\noi We investigate two  critical cases (i.e. $V$ satisfying
(\ref{16intro}, \ref{15intro})) in section \ref{Bes}.

\noi In section \ref{LD}, using the technique of large deviations
we deal with $V$ fullfills (\ref{15intro}).

\noi We end this paper by considering the case where $V$ is large
at $\pm \infty$  in section \ref{bilateral}.

\noi The results of this paper were announced without proofs in
\cite{RoyValYor}.

\noi 7) In a subsequent study \cite{RVY2}, we consider a similar
problem, replacing the exponential weight \\
$\displaystyle \exp \Big \{-\frac{1}{2} \int _ 0^t V(B_ v )dv \Big
\}$ by $\varphi (A_t)$, where $(A_t)$ may be equal either to the
one-sided maximum : $\displaystyle \sup_{0\leq u \leq t}B_u$, or
to the one-sided minimum, or to the local time at $0$, or to the
number of down-crossings from level $b$ to level $a$.


\section {Preliminaries }\label{prel}

\setcounter{equation}{0}
 {\bf 2.1 } Let $ \varphi \ : \mathbb R \rightarrow ]0,+\infty[$ be a function
of class $C^2$ and $\mu$ the measure on $\mathbb R$ with density
$\varphi^2 (x)$ with respect to the Lebesgue measure:
\begin {equation} \label{def mu}
\mu(dx) = \varphi^2(x)dx.
\end {equation}

\noi We denote by $L^\varphi$  the differential operator:
\begin {equation} \label{def oper L}
L^\varphi f(x)= \frac{1}{2}f"(x)+\frac{\varphi'(x)}{\varphi(x)}
f'(x),
\end {equation}
\noi defined for every function $f$ of class $C^2$.

\noi If $f$ and $g$ are two functions of class $C^2$, with compact
support, then by  integration by parts we obtain:
\begin {equation} \label{sym mu}
<L^\varphi f , g> _\mu = <f , L^\varphi g > _\mu
=-\frac{1}{2}<f',g'> _\mu=
 -\frac{1}{2}\int_\mathbb R f'(x)g'(x)d \mu (x),
\end {equation}

\noi where $\displaystyle <h,k> _\mu = \int  _\mathbb R h(x)k(x)d
\mu (x)$.

\noi The relation (\ref {sym mu}) tells us that $L^\varphi$ is a
negative and symmetric operator, defined on $C^2_K(\mathbb{R})$.
Thus, it admits a self-adjoint extension, which is the generator
of a Markovian semigroup $(T _t^\varphi;t\geq 0)$ of bounded,
positive, symmetric operators on $L^p (\mu)), \ 1 \leq p \leq
\infty$ (cf \cite{Davi}). The norm of $T _t^\varphi$ is $1$ as an
operator on any $L^p (\mu)$.


{\bf 2.2 } Let $X _t^x$ be the solution of the following
stochastic differential equation:
\begin {equation} \label{def proc X}
X_ t=x+ B _t+\int  _0 ^t \frac{\varphi'}{\varphi}  (X_ s )ds; \ t
\geq 0 ,
\end {equation}
\noi where $(B _t; t \geq 0)$ is a one-dimensional Brownian motion
started at $0$.

\noi Since $\varphi$ is of class $C^2$ and $\varphi >0$ this
stochastic differential equation has a unique strong solution up
to an explosion time. We assume that this explosion time is
infinite. This occurs for instance if $\displaystyle
\frac{\varphi'}{\varphi}$ has at most  linear growth; for more
refined conditions see \cite{McKean}.

\noi Obviously:
\begin {equation} \label{semig}
E [f(X _t^x)]=T _t^\varphi f(x), \ {\rm for \ any } \ f \geq 0 ,
\end{equation}
and by the Girsanov formula :
\begin {equation} \label {loi1 X}
E _x\Big [ f(B_t)\exp \Big \{ -\frac{1}{2} \int
_0^t\frac{\varphi"}{\varphi} (B _s )d s \Big \} \Big ] =
\varphi(x) T_ t^\varphi(f/\varphi)(x)= \varphi(x)
E\Big[\frac{f(X_t^x)}{ \varphi(X_ t^x )}\Big].
\end{equation}
In particular choosing $f=1$ we get :
\begin {equation} \label {loi2 X}
Z_t^{V_\varphi}(x)=E _x\Big [ \exp \Big \{ -\frac{1}{2} \int
_0^t\frac{\varphi"}{\varphi} (B _s )d s \Big \} \Big ] =
\varphi(x) T_ t^\varphi(\frac{1}{\varphi})(x)= \varphi(x)
E\Big[\frac{1}{ \varphi(X_ t^x )}\Big].
\end{equation}
%


\begin {rem}  If  $\varphi$ is locally the difference of two convex
functions, then it is understood that $\displaystyle \int
_0^t\frac{\varphi"}{\varphi} (B _s )d s$ is  defined as
$\displaystyle\int _ \mathbb{R} L^x_t \frac{\varphi"(dx)}{\varphi
(x)}$.
\end{rem}


\section{ The case : $\varphi$ small at infinity} \label{Csmall}
\setcounter{equation}{0}
 Let $ \varphi \ : \mathbb R \rightarrow ]0,+\infty[$ be a
function of class $C^2$. We suppose moreover that
$\varphi'/\varphi$ is bounded.

\noi We define:
\begin{equation} \label{def u}
V_\varphi(x)=\frac{\varphi"(x)} { \varphi(x)}, x \in \mathbb R.
\end{equation}
\noi More generally, if  $\varphi$ is locally the difference of
two convex functions, we set :
\begin{equation} \label{def ubis}
V_\varphi(dx)=\frac{\varphi"(dx)} { \varphi(x)}, x \in \mathbb R.
\end{equation}

\noi
  In this section, we assume that $\varphi$ is small at infinity, in the
  sense that :

\begin{equation} \label{hyp1C}
\int  _\mathbb R\varphi^p(x) dx <\infty , \ {\rm for \ some } \
0<p<1,
\end{equation}

 \noi and
\begin{equation} \label{hyp2C}
\varphi \ {\rm is \ decreasing\ (resp. \ increasing) \ at }
+\infty {\rm (resp. } -\infty).
 \end{equation}

\noi It is clear that the sign of $V_\varphi$ is not constant.

 \noi We note that (\ref{hyp1C}) and
(\ref {hyp2C}) imply that $\int  _\mathbb R \varphi^2(x)
dx<\infty$ and the change $\varphi \rightarrow \lambda \varphi$,
with $\lambda
>0$, does not modify $V_\varphi$, nor (\ref {hyp1C}), nor (\ref {hyp2C}).


\begin {theorem} \label{theoC}
We suppose that $\varphi$ satisfies (\ref {hyp1C}), (\ref {hyp2C})
and is even, i.e. $\varphi(-x)=\varphi(x), \forall x \in
\mathbb{R}$.
\begin {enumerate}
\item The generic Theorem applies with $k=0$ since :

\begin{equation} \label{ergo}
 \lim _{t\rightarrow \infty} \Big ( T^\varphi
_ t(1/\varphi)(x) \Big ) = \frac{\int  _\mathbb R
\varphi(y)dy}{\int  _\mathbb R \varphi^2(y)dy},
\end{equation}

\noi and

 \begin{equation} \label{born1}
\int  _ \mathbb R h(x)\varphi^2(x) dx <\infty \qquad \mbox { where
} \quad
 h(x)=\sup_{t \geq 0}  | T^\varphi _ t(1/\varphi)(x)
|.
\end{equation}

\item $Q^{V_\varphi}_{x,t}$ converges weakly to $P^\varphi _x$, as
$t\rightarrow \infty$.

 \item
 Let $(X_ t^x;t\geq 0)$ be the solution of the SDE :

\begin{equation} \label{def2 proc X}
X _t=x+B _ t+\int  _0 ^t {\varphi' \over \varphi}(X _s)ds, \ t
\geq 0.
\end{equation}

\noi Then  the law of $(X_ t^x;t\geq 0)$ is $P^\varphi _ x$.

\noi Moreover $(X_ t^x;t\geq 0)$ is a recurrent process with
 finite invariant measure $\mu (dx)=\varphi^2(x) dx$.
\end {enumerate}
\end{theorem}


\vskip 5 pt \noi Before proving Theorem \ref{theoC}, we give five
examples numbered from \ref{exC1} to \ref{exC5}. For these
examples Theorem \ref{theoC} applies because $(T^\varphi _t ;
t\geq 0)$ is an ultracontractive \cite{KKR} or an hypercontractive
semigroup (\cite{Nel}, \cite{Gross}, \cite{Neveu}).
\begin{exam} \label{exC1}
Let $\varphi$ be the function :
$$\varphi(x)=  e^{-\frac{|x|^\alpha }{ 2}}, x \geq 0,$$
\noi where $\alpha >2$.

 \noi Then $\varphi$ obeys (\ref{hyp1C})
(in fact, for any $p>0$), (\ref{hyp2C}), and

$$V_\varphi(x)={1 \over 4} \alpha ^2|x|^{2\alpha -2}
-{1 \over 2} \alpha (\alpha -1) |x|^{\alpha -2} , \ x\geq 0. $$

\noi $(T^\varphi _t ; t\geq 0)$ is an ultracontractive semigroup
(i.e. for any positive $t$, $T^\varphi _t$ is a bounded  operator
from $L^1(\mu)$ to $L^\infty (\mu)$) and this implies directly
Theorem  \ref{theoC}.

\noi More generally,  we can  take $\varphi(x)=e^{-v(x)/2}$ where
$v(x)$ is a convex function for large $x$ and $\int ^{+\infty}
{1\over v'(x)}dx<\infty$. Theorem \ref{theoC} remains valid since
$(T^\varphi _t ; t\geq 0)$ is still an ultracontractive semigroup
\cite{KKR}.
\end {exam}


\begin{exam} \label{exC2}
Let $\varphi(x)=e^{-x^2/2}; x \in \mathbb R$. Then
$V_\varphi(x)=x^2 -1$ and $(X _t^x;t \geq 0)$ is the Ornstein
Uhlenbeck process which solves :
$$X _t=x+ B _ t - \int  _0 ^tX _ s ds, \ t\geq 0.$$
\noi Notice that $(T^\varphi _t ; t\geq 0)$ is  an
hypercontractive semigroup (cf \cite{Nel}, \cite{Gross},
\cite{Neveu}).
\end{exam}


\begin{exam} \label{exC3}
Let $\varphi$ satisfy :

\begin{equation}\label {id2 1}
-\Big ( \frac{\varphi'}{\varphi} \Big )'=\frac{\varphi'^2 -
\varphi \varphi"}{ \varphi^2} \geq 2 \kappa .
\end{equation}

\noi For every pair of functions $f$ and $g$ of class $C^2$ with
compact support, we recall that :

$$\Gamma ^\varphi(f,g) := L^\varphi(fg)
-fL^\varphi g-gL^\varphi f,$$

$$\Gamma _2 ^\varphi (f,g) :=  L^\varphi (\Gamma (f,g))-
\Gamma (L^\varphi f,g) - \Gamma (f,L^\varphi g ).$$

\noi Then (\cite{BakEme}) the operator $L^\varphi$ enjoys the
spectral gap property in $L^2(\mu)$ as soon as

\begin{equation}\label {Gamma2}\Gamma _2^\varphi (f,f) \geq \kappa \Gamma
^\varphi(f,f).
\end{equation}

\noi It is easy to check that (\ref{id2 1}) implies
(\ref{Gamma2}). Theorem \ref{theoC} follows immediately.
\end{exam}

\begin{exam} \label{exC4}
 Let $a\geq 0$ and $\varphi$ such that:
$$\varphi(x)= \left \{ \begin{array}{ll}
e^{-|x|} & \ {\rm if} \ |x|>a \\
  e^{-a}(1+a-|x|) & {\rm otherwise.}
\end{array} \right. $$
\noi Then
$$V_\varphi (dx)=\frac{\varphi"(dx)}{ \varphi(x)}
= 1_ {\{ |x|>a\}} dx -\frac{2}{ 1+a}\delta _0(dx),$$
\noi where $\delta _0(dx)$ denotes the Dirac measure at $0$.

\noi Consequently:

$$E  _x\Big [ \exp \Big \{ -\frac{1}{2} \int  _\mathbb R L _t^y V_\varphi(dy) \Big \} \Big ]
= E  _x\Big [ \exp \Big \{ -\frac{1}{2} \int  _0^t1 _{\{ |B
_s|>a\}} ds +\frac{1 }{ 1+a}
 L _t^0 \Big \} \Big ] $$

\noi and $(X _t^x;t \geq 0)$ solves:

$$X _t = x+ B _t-\int _0^t {\rm sgn}(X _s)1 _{\{|X _s|>a\}}ds - \int _0^t
 \frac{{\rm sgn}(X _s)}{  1+a-|X _s |}1 _{\{|X _s|\leq
a\}}ds, \quad t \geq 0.$$

\end {exam}


\begin{exam} \label{exC5}
 Let $\varphi(x)=e^{-\lambda |x|}$,
with $\lambda >0$. Then $V_\varphi(dx)= \lambda^2 dx -2 \lambda
\delta _0(dx)$,
$$\displaystyle
\lim_{t \rightarrow \infty} \bigg\{ \frac{ E _x \Big [ 1 _{\Lambda
_s} \exp \big \{ \lambda L^0 _t\big \} \Big ] }{ E _x \Big [  \exp
\big \{ \lambda L^0 _t\big \} \Big ] } \bigg\} = e^{\lambda |x|} E
_x \Big [ 1 _{\Lambda  _s} \exp \big \{ -\lambda |B _s| +\lambda
L^0 _s -\frac{\lambda ^2 s}{ 2}\big \} \Big ],$$
 \noi and
$(X _ t^x ; t \geq 0)$ solves:
$$X_ t=x+B_ t - \lambda \int _0^t {\rm sgn}(X_ s)ds.$$
\noi $(X _ t^x ; t \geq 0)$ is the so-called bang-bang process
with parameter $\lambda >0$ (cf \cite{Fitzs}, \cite{GraSveShi},
\cite{ShiChe}) which satisfies, for $x=0$:
\begin{equation}\label{bang-bang}
    (|X^0_t| ; t\geq 0) \stackrel {(d)} {=}\ (S^{(\lambda)}_t-
    B^{(\lambda)}_t ; t\geq 0),
\end{equation}
where $\displaystyle B^{(\lambda)}_t=B_t + \lambda t, \
S^{(\lambda)}_t=\sup_{0 \leq u \leq t} B^{(\lambda)}_u$.
\end{exam}

\vskip 10 pt \noi The proof of Theorem \ref{theoC} is based on two
preliminary results which we present in  Lemmas \ref{lemprel1} and
\ref{lemprel2}.

\begin {lemma}\label{lemprel1}
Let $\rho : [0,+\infty[ \times \mathbb R \rightarrow \mathbb R,\
\rho (t,x)= T^\varphi _t(1/\varphi)(x)$. Then
\begin {enumerate}
\item For any $t\geq 0$, $x \rightarrow \rho (t,x)$ is even and
non-decreasing on $[0,+\infty[$. \item  $t \rightarrow \rho (t,0)$
is  non-decreasing.
\end{enumerate}
\end {lemma}


\begin{prooff}\  {\bf of Lemma \ref{lemprel1}}.

i) It is well-known that $\rho$ solves:

\begin{equation} \label{edp rho}
\left \{ \begin{array}{l}
 \frac{ \partial \rho }{ \partial t}
-\frac{1}{2} \frac{\partial^2 \rho }{ \partial x^2}
-\frac{\varphi' }{
\varphi} \frac{\partial \rho }{ \partial x} = 0 \\
 \rho (0,x)  =\frac{1 }{\varphi(x)}
\end{array}
\right .
\end{equation}

 \noi We set : $\displaystyle \theta (t,x)=\frac{
\partial \rho}{\partial x}(t,x)$.

\noi We take in (\ref{edp rho}) the partial derivative with
respect to  $x$ :

\begin{equation} \label{edp theta}
\left \{ \begin{array}{l} \frac{ \partial \theta }{ \partial t}
-\frac{1}{2} \frac{\partial^2 \theta }{
\partial x^2} -\frac{\varphi"\varphi-\varphi'^2 }{ \varphi^2 }
\theta -\frac{\varphi' }{ \varphi}
\frac{\partial \theta }{ \partial x} = 0 \\
 \theta (0,x)
=\frac{\partial }{ \partial x} \Big (\frac{1}{ \varphi(x)} \Big )
\cr \theta (t,0)=0
\end{array}
  \right .
  \end{equation}

\noi Since the restriction of $\varphi$ to $[0,+\infty[$ is
non-increasing, $\theta (0,x)\geq 0$ if $x\geq 0$.

\noi It is clear that $\rho (t,.)$ is even, i.e. $\rho (t,-x)=\rho
(t,x); \forall x \in \mathbb{R}$, consequently $\theta (t,0)=0$.

\noi But $w(t,x)\equiv 0$ is a solution to (\ref{edp theta}) on
$[0,+\infty[ \times [0,+\infty[$; then, the maximum principle
implies that $\frac{
\partial \rho }{ \partial x}(t,x) \geq 0$ if $x\geq 0$. This
proves point i).

\noi ii) Since $(X_ t^x;t\geq 0)$ is a Markov process and $\rho
(t,.)$ is even :
$$\rho (t+s,0) =T_ {t+s}^\varphi\big ( \frac{1}{ \varphi} \big )(0)
=E\big [T^\varphi _t\big ( \frac{1}{ \varphi} \big) (X _ s^0)\big
]= E\big [\rho(t, |X _ s^0|)\big ] \geq \rho (t,0).$$ \noi This
proves ii).
\end{prooff}


\begin {lemma}\label{lemprel2}
Let  $h : \mathbb R \rightarrow \mathbb R$ be the function :
$\displaystyle h(x):=\sup  _{t\geq 0}^{} T^\varphi _t \big
(\frac{1}{\varphi} \big )(x)$, then :

\begin{equation} \label{maj T1}
\int_\mathbb{R}h(x)\varphi^2(x)dx<\infty.
\end{equation}
\end{lemma}

\begin {prooff} \ {\bf of Lemma \ref{lemprel2} }.

\noi Let $x \geq 0$. Due to point 1. of Lemma \ref{lemprel1}, we
have:

\begin{equation} \label{min T2}
<T^\varphi _t \big (\frac{1}{\varphi} \big ) , 1 _ {[x,+\infty[} >
_\mu \ \geq \Big(T^\varphi _t \big (\frac{1}{\varphi}\big
)(x)\Big) \mu([x,+\infty[).
\end{equation}

 \noi Recall
that $\mu(dy)=\varphi^2(y)dy$ and   $(T^\varphi _t;t \geq 0)$ is a
symmetric semigroup, then

$$<T^\varphi _t \big (\frac{1}{\varphi} \big ) , 1 _
{[x,+\infty[} > _\mu= <\frac{1}{\varphi}, T^\varphi _t \big ( 1 _
{[x,+\infty[} \big )> _\mu =\int _\mathbb R \frac{T^\varphi _t
\big ( 1 _ {[x,+\infty[} \big )(y) } {\varphi(y)} \mu (dy).$$

 \noi Let $p'=2-p>1$ and $q'$ be the
conjugate number ($1/p'+1/q'=1$).

\noi We use H\"{o}lder's inequality and the fact that $T _
t^\varphi$ is a bounded operator from $L^{q'}(\mu)$ to itself with
norm equal to $1$ :

$$T^\varphi _t \big (\frac{1}{\varphi} \big )(x) \  \leq  \tilde{h}(x).$$

\noi where

$$\tilde{h}(x) =C \Big ( \mu([x,+\infty[) \Big )^{1/q'-1},$$

\noi and

$$C= \Big (\int _\mathbb R \frac{1 }{ \varphi(y)^{p'}}\mu (dy) \Big )^{1/p'}=
\Big (\int _\mathbb R  \varphi(y)^p dy \Big )^{1/p'}<\infty ,$$

 \noi since $ 2-p'=p$ and $\varphi$ satisfies (\ref{hyp1C}).

\noi As for (\ref{maj T1}), we have:

$$\frac{1 }{ C} \int _ \mathbb R \tilde{h}(x)\mu (dx)=
  \int _ \mathbb R \Big ( \mu([x,+\infty[) \Big )^{1/q'-1} \mu (dx)=
-q' \Big [  \mu([x,+\infty[ )^{1/q'}\Big ]
_{x=-\infty}^{x=+\infty}=q'<\infty.$$

\end {prooff}


\begin{prooff} \ {\bf of Theorem \ref{theoC}}.

a) We start with the proof of point 1.

\noi Let us introduce the function $\underline \theta$:

$$\underline \theta (x)=  \underline { \lim} _ {t \rightarrow \infty} \Big (
T _t^\varphi(1/\varphi)(x)\Big ).$$

\noi Lemma \ref{lemprel2} implies that $\underline \theta$ is
$\mu$-integrable. Moreover:

$$T _s^\varphi(\underline \theta)=T _s^\varphi
\Big ( \liminf _ {t \rightarrow \infty} \Big \{ T
_t^\varphi(1/\varphi)\Big \} \Big ) \leq  \liminf _ {t \rightarrow
\infty} \Big (T _{t+s}^\varphi(1/\varphi)\Big )=\underline
\theta.$$

\noi Consequently:

\begin{equation}\label{ineg thb}
T _s^\varphi(\underline \theta) \leq \underline \theta .
\end{equation}

\noi The semigroup $(T_ t^\varphi;t\geq 0)$ being $\mu$-symmetric,
we have:

$$< T _t^\varphi(\underline \theta),1> _\mu =< \underline \theta , T _t^\varphi(1) > _\mu
=< \underline \theta ,1>  _\mu =\underline \theta \int
_\mathbb{R}\varphi^2(x)dx.$$

\noi This equality, together with the inequality (\ref{ineg thb})
implies that $\underline \theta =T _t^\varphi(\underline \theta)$.

\noi If we take the derivative with respect to $t$, we obtain:
$L^\varphi(\underline \theta)=0$.

\noi Using (\ref{sym mu}) we have:

$$< L^\varphi(\underline \theta),\underline \theta> _\mu = -\frac{1}{2}
 \int _\mathbb R  \underline \theta'(x)^2 d\mu (x)=0.$$

\noi Consequently $\underline \theta =\underline C$, where
$\underline C$ is a constant.

\noi
 We introduce:

$$\overline \theta (x)=  \limsup _ {t \rightarrow \infty} \Big (
T _t^\varphi(1/\varphi)(x)\Big ).$$

 \noi In the same way as before, we easily check that
$\overline \theta =\overline C$.


\vskip 5 pt  b) We now prove that $\displaystyle \overline \theta
=\underline \theta = \frac{\int _\mathbb R \varphi(x)dx}{\int
_\mathbb{R}\varphi^2(x)dx}$.

\noi Let $(x  _n;n \geq 1)$ and  $( \varepsilon  _n;n \geq 1)$ be
two  sequences such that $( \varepsilon  _n;n \geq 1)$ is
positive, decreasing and

$$ \lim  _{n \rightarrow \infty} x _n
=\lim  _{n \rightarrow \infty} \varepsilon _n =0.$$

\noi Suppose first that $x _ n$ and $\varepsilon  _n$ are given.
Using the definition of $\overline \theta$, there exists $t _n$
such that:

$$T _{t_ n}^\varphi (1/\varphi)(x _n) \geq \overline \theta - \varepsilon  _n.$$

\noi Moreover we can choose $t _n$ in  such a way that $(t  _n;n
\geq 1)$ is  an increasing  sequence converging to $+\infty$ as
$n\rightarrow \infty$.

\noi Let $x>0$ fixed. Since $x \rightarrow T^\varphi
_t(1/\varphi)(x)$ is non-decreasing on $[0, +\infty[$ (cf Lemma
\ref{lemprel1}) , if $n$ is large enough:

$$T_ {t_ n}^\varphi(1/\varphi)(x )\geq T_ {t_ n}^\varphi(1/\varphi)(x _n ).$$

\noi Taking the limsup on both sides we obtain:
$$\lim  _{n \rightarrow \infty}  \Big (
T_ {t_ n}^\varphi(1/\varphi)(x ) \Big )= \overline \theta .$$
\noi Thanks  to Lemma \ref{lemprel2}, we can apply the dominated
convergence theorem, hence

$$ \lim  _{n \rightarrow \infty}  \Big (
< T_ {t_ n}^\varphi(1/\varphi) , 1 > _\mu \Big )= <\overline
\theta , 1
> _ \mu =\overline \theta \int
_\mathbb{R}\varphi^2(x)dx.$$

 \noi But  recall that since $( T^\varphi _t;t\geq 0)$
is $\mu$-symmetric, then :

$$< T_ {t_ n}^\varphi(1/\varphi) , 1 > _\mu
= < \frac{1}{\varphi} ,  T_ {t_ n}^\varphi(1) > _\mu = \int
_\mathbb R \varphi(x) dx.$$

\noi Consequently

$$\overline \theta = \frac{\int  _\mathbb R \varphi(x) dx}{\int _\mathbb R
\varphi^2(x)dx}.$$

\noi Replacing $(x  _n;n \geq 1)$ by $(y  _n;n \geq 1)$ such that
$y _n<0$ and

$$\lim  _{n \rightarrow \infty} y _n=0,$$

\noi we prove similarly that $\underline \theta =\overline
\theta$.


\vskip 5pt  c)  Let $t>0$. Recall that $Q _{x,t}^{V_\varphi}$ is
the probability defined on ${\cal F} _t$ by :
\begin{equation}\label{phismall11}
Q _{x,t}^{V_\varphi}(\Lambda  _t) = \frac{ E _x \Big [1 _{\Lambda
_t} \exp \Big \{-\frac{1}{2} \int _ 0 ^t V_\varphi(B_ r )dr \Big
\} \Big ] }{ E _x \Big [ \exp \Big \{-\frac{1}{2} \int _ 0 ^t
V_\varphi(B_ r )dr \Big \} \Big ] },\ \Lambda _t \in {\cal F} _t.
\end{equation}

\noi Suppose that $s>0$ is fixed and pick $t>s$; then, replacing
in (\ref{phismall11}) $\Lambda _t$ by $\Lambda _s \in {\cal F}
_s$, and, using the Markov property at time $s$ together with
(\ref{loi1 X}), we obtain:
%
$$Q _{x,t}^{V_\varphi}(\Lambda  _s) = \frac{ E _x \Big [1 _{\Lambda  _s}
\exp \Big \{-\frac{1}{2} \int _ 0 ^s V_\varphi(B_r )dr \Big \} \
\varphi(B _s) T^\varphi _{t-s}(1/\varphi)(B _ s)\Big ] }{
\varphi(x)T^\varphi _t(1/\varphi)(x) }.$$

\noi The numerator can be written as $E _x[1 _{\Lambda  _s}Y_
{s,t}]$, where :
\begin{equation}\label{ajout11}
    Y_{s,t}=\exp \Big \{-\frac{1}{2} \int _ 0 ^s V_\varphi(B_ r )dr \Big \}
\ \varphi(B _s) T^\varphi _{t-s}(1/\varphi)(B _ s) .
\end{equation}

\noi On one hand, using (\ref{born1}), we get an upper bound for
$Y _{s,t}$:

$$0 \leq Y _{s,t} \leq Y_s, \ {\rm for\  any\ } 0<s<t $$

\noi where

$$  Y_s= \varphi(B _s) h(B _s) \exp \Big \{-\frac{1}{2} \int _ 0 ^s
V_\varphi(B_ r )dr \Big \}.$$

\noi Identity (\ref{loi1 X}) tells us that:

$$E _x[Y_s] = \varphi(x) T _s^\varphi h(x) <\infty. $$

\noi On the other hand, (\ref{ergo}) and (\ref{born1}) imply that

$$  \lim  _{t \rightarrow \infty}  Y_ {s,t}=
\lambda \varphi(B _s)1 _{\Lambda  _s} \exp \Big \{-\frac{1}{2}
\int _ 0 ^s V_\varphi (B_ r )dr \Big \},$$

$$  \lim  _{t \rightarrow \infty}  T^\varphi _t(1/\varphi)(x)  = \lambda,
\quad \mbox{where } \ \lambda = \frac{\int  _ \mathbb R
\varphi(x)dx}{\int  _ \mathbb R \varphi^2(x)dx}.$$

\noi Then, for any $s>0$,  $Q _{x,t}^{V_\varphi}(\Lambda  _s)$
converges to $ P^\varphi _{x}(\Lambda  _s)$ as $t \rightarrow
\infty$, where $P^\varphi_ x$ is the probability defined on ${\cal
F}_\infty$ by :

$$ P^\varphi _{x}(\Lambda  _s) = \frac{1 }{ \varphi (x)} \
 E _x \Big [1 _{\Lambda  _s} \varphi(B _s) \exp \Big \{-\frac{1}{2}
  \int _ 0 ^s V_\varphi (B_ s )ds \Big \}\Big ], $$

 \noi  for $s>0$ given and any $\Lambda  _s$ in   ${\cal F} _s$.

\noi Point iii) of Theorem \ref{theoC} is a direct consequence of
Girsanov formula.

\noi Moreover $(X_ t^x ;t\geq 0)$ is recurrent since  $\mu$ is its
invariant  measure.
\end{prooff}


 \section { The integrable case}\label{integrable case}
\setcounter {equation}{0}
 Throughout this section, $V(dx)$ shall
always denote a finite positive Radon measure  on $\mathbb R$,
different from $0$, with finite first moment;  hence :
\begin{equation}\label{propu}
\int _\mathbb R (1+|x| ) V(dx)<\infty.
\end{equation}
\noi Recall that in the previous section the initial data was the
function $\varphi$, whereas now the data is the potential $V$.

\vskip 10 pt \noi
\begin{theorem} \label{thbase1}
Let $V(dx)$ be a finite positive Radon measure on $\mathbb R$
fulfilling (\ref{propu}).
\begin{enumerate}
 \item
 The generic Theorem applies with $k=1/2$.

\item
 $\varphi_V$ is a convex function which takes its values in
$]0,\infty[$ and is the unique solution to the Sturm-Liouville
 equation
\begin{equation}\label{sturm}
\varphi"(dx)=\varphi(x)V(dx),
\end {equation}
\noi with boundary conditions :
\begin{equation}\label{sturmcond}
\lim  _ {x\rightarrow +\infty} \varphi_V'(x)= - \lim  _
{x\rightarrow -\infty} \varphi_V'(x) =\sqrt { \frac{2}{\pi} }.
\end {equation}
\noi  As a consequence
\begin{equation}\label{maj1C}
 \varphi_V(x) \sim_{|x| \rightarrow
\infty} \sqrt { \frac{2}{\pi}  } |x|
\end{equation}
\noi and
\begin{equation}\label{maj2C}
 \ \varphi_V(x)\leq C (1+|x|).
\end{equation}
\item
 Let $M ^{\varphi_V}$ be the process:
\begin{equation}\label{maj2CA}
    M ^{\varphi_V}(s)= \varphi_V(B _ s) \exp \Big \{-{\frac{1}{2}}
\int  _ \mathbb R L _s^y V(dy) \Big \} , s \geq 0 .
\end{equation}
\noi Then $(M ^{\varphi_V}(s);s\geq 0)$ is a  martingale such that
$E[(M ^{\varphi_V}(s))^2]<\infty$ for any $s\geq 0$ (recall that
$\displaystyle V_{\varphi_V}=V$).

 \item
 Let  $(X_t ^x;t\geq 0)$ be the solution to :
\begin{equation} \label{defX}
X _t=x+B _ t+\int  _0 ^t \frac{\varphi_V'}{\varphi_V}(X _s)ds, \ t
\geq 0 .
\end{equation}
 Then the law of $(X_t ^x;t\geq 0)$ is $P _ x^{\varphi_V}$.

\item
 The process $(X_t ^x;t\geq 0)$ is transient. More precisely, denoting
$$\displaystyle{ \rho= \int  _\mathbb R\frac{dy}{\varphi_V^2(y)}
<\infty},$$ \noi then :
\begin{equation} \label{lim1X}
P\Big ( \lim  _{t \rightarrow \infty} X _t^x=-\infty \Big ) =
\frac{1}{\rho} \int  _x ^{+\infty}\frac{dy}{\varphi_V^2(y)} ,
\end{equation}
\begin{equation} \label{lim2X}
P \Big ( \lim  _{t \rightarrow \infty} X _t^x=+\infty \Big ) =
\frac{1}{\rho} \int  ^x  _{-\infty}\frac{dy}{\varphi_V^2(y)}.
\end{equation}

\end{enumerate}
\end{theorem}

\begin{rem} \label{genBes01 } Theorem \ref{thbase1} may be
generalized replacing the Brownian motion $(B_t ;t \geq 0)$ by a
Bessel process $(R_t ;t \geq 0)$ of dimension $0<d<2$. In this
case, the generic Theorem applies with  a function $V$ with
compact support and $k=1-\frac{d}{2}$. On the other hand, we have
not been able to settle the case $d=2$.

\end{rem}

 \vskip 15 pt \noi We actually
develop two proofs of Theorem \ref{thbase1}. The first one is
based on the study of the function $\displaystyle t\mapsto
Z_t^{V}(x):=E _x \Big [ \exp \big \{-\frac{1}{2} \int _\mathbb R
L^y _ t V(dy) \big \} \Big]$. The second one  relies upon the
excursion theory and the Ray-Knight theorem which describes  the
distribution of $(L^y_S;y\in\mathbb{R})$, where $S$ is an
exponential r.v. independent of the Brownian motion.
\subsection {An analytical approach}\label{int1}

 \vskip 20 pt
 \noi Let us briefly describe our first proof of Theorem
\ref{thbase1}. The crucial point is an a priori inequality
concerning $ Z_t^{V}(x)$ stated in Lemma \ref{etap1} below. To
demonstrate that $\sqrt{t}\ Z_t^{V}(x)$ has a limit, when
$t\rightarrow\infty$ we prove that  the normalized Laplace
transform $A(\lambda ,x)$ of $Z_t^{V}(x)$ converges as
$\lambda\rightarrow 0$ (cf Lemma \ref{Lap1}).This can be done
(Lemma \ref{Lap2}) through properties involving $A(\lambda ,x)$
and its derivatives. \vskip 10 pt

\begin {lemma} \label{etap1}
Let $V(dy)\not = 0$ be a positive Radon measure on the whole line,
satisfying (\ref{propu}). Then there exists a constant $C$ such
that:
\begin{equation}\label {maj1}
\sqrt {1+t} E _x \Big [ \exp \big \{-\frac{1}{2} \int  _\mathbb R
L^y _ t V(dy) \big \} \Big] \leq C(1+|x|),\ t \geq 0,\ x\in
\mathbb R .
\end{equation}
\end{lemma}

 \begin{proof}


  1) We start with $V(dy)=\gamma \delta  _x(dy)$, where $\delta  _x(dy)$ denotes
the Dirac measure at $x$. We claim that:
\begin{equation}\label{majtempsloc}
E_ 0\big [\exp \{ -\gamma L^x  _t\}\big ]\leq \sqrt{ \frac{2}{\pi
t} } \big (|x|+\frac{1}{\gamma} \big ), \  x \in \mathbb R, \gamma
>0, t \geq 0.
\end{equation}
Observing that  $L^x_t$ is distributed as $(L_t-|x|)_+$ (cf
\cite{RevYor}), then
$$E_0[\exp{-\gamma L^x_t}]=E_0[\exp{-\gamma (L_t-|x|)_+)}]
=P(L_t\leq |x|)+\sqrt{\frac{2}{\pi t}}\int_{|x|}^{\infty}
e^{-\gamma (y-|x|)}e^{\frac{-y^2}{2t}} dy,$$
$$=P(|B_1|\leq \frac{|x|}{\sqrt{t}})+
\sqrt{\frac{2}{\pi }}e^{\gamma |x|+\gamma ^2t/2} \int_{
\frac{|x|}{\sqrt{t}}+\gamma \sqrt{t}} ^{\infty}e^{\frac{-z^2}{2}}
dz,$$
$$ \leq \sqrt{\frac{2}{\pi }} \frac{|x|}{\sqrt{t}}+
\sqrt{\frac{2}{\pi }} e^{\gamma |x|+\gamma ^2t/2}
e^{-\frac{1}{2}(\frac{|x|}{\sqrt{t}}+\gamma \sqrt{t})^2}
\frac{1}{\frac{|x|} {\sqrt{t}}+\gamma \sqrt{t}}$$
$$\leq \sqrt{\frac{2}{\pi }}\frac{|x|}{\sqrt{t}}
+ \sqrt{\frac{2}{\pi }} \frac{1}{\frac{|x|} {\sqrt{t}}+\gamma
\sqrt{t}} e^{-\frac{x^2}{2t }}\leq \sqrt{\frac{2}{\pi t}}(|x|+
\frac{1}{\gamma}).$$


2) Let $V(dy)\not =0$ be a positive Radon measure on $\mathbb R$.
We choose $a$ and $b$ such that $a<b$ and $\mu=V([a,b])/2>0$. We
have:
$$\exp \big \{-\frac{1}{2} \int  _\mathbb R
L^y _ t V(dy) \big \}\leq \exp \big \{-\frac{1}{2} \int  _a^b L^y
_ t V(dy)\}.$$
\noindent
 Since $x\mapsto e^{-\mu x}$ is convex:
$$\exp \big \{-\frac{1}{2} \int  _a^b L^y _ t
V(dy) \big \}\leq \frac{1}{2\mu} \int  _a^b \exp\{-\mu L^y_t\} \
 V(dy),$$
\noindent Taking the expectation and  applying
(\ref{majtempsloc}), we obtain
$$E _x \Big [ \exp \big \{-\frac{1}{2} \int  _\mathbb R
L^y _ t V(dy) \big \} \Big] \leq \frac{C_1}{\sqrt{1+t}} \int _a^b
(|x-y|+\frac{1}{\mu})V(dy).$$
 \noindent Then (\ref {maj1})
follows.
\end{proof}

\begin {lemma} \label{Lap1}
Let $V(dx)$ be a finite positive Radon measure as in Theorem
\ref{thbase1} and $A$ the Laplace transform:
\begin{equation}\label {defA}
A(\lambda,x)=\int_0^\infty e^{-\lambda t} Z_t^{V}(x) dt,
\end{equation}
\noindent where
\begin{equation}\label {defphi}
Z_t^{V}(x)=E _x\Big [ \exp \Big \{-{\frac{1}{2}} \int  _ \mathbb R
L _t^y V(dy) \Big \} \Big ].
\end{equation}
\noindent Then $\displaystyle{\lim_{t\rightarrow \infty}\sqrt{t}
Z_t^{V}(x)= \varphi_V(x)}$ if and only if
\begin{equation}\label {convA}
\lim_{\lambda\rightarrow 0}\sqrt{2\lambda}A(\lambda,x)=
\sqrt{2\pi}\varphi_V(x).
\end{equation}
\end{lemma}

 \begin{proof} \ We set
 %
\begin{equation}\label {defAti}
\tilde{A}(\lambda,x)=\sqrt{2\lambda}A(\lambda,x).
\end{equation}
\noindent The inequality (\ref{maj1}) implies that
\begin{equation}\label {ineg1}
\tilde{A}(\lambda,x)\leq \kappa (1+|x|), {\rm for\  all\ }\lambda
\geq 0.
\end{equation}
\noindent Since $t \mapsto Z_t^{V}(x)$ is a decreasing function, a
classical version of the Tauberian theorem (cf \cite{Fell}, Chap.
XIII, section 5 ) implies that $\displaystyle{\lim_{t\rightarrow
\infty}\sqrt{t}Z_t^{V}(x)= \varphi_V(x)}$ if and only if
(\ref{convA}) holds.
\end{proof}


\noi \ Lemma \ref{Lap1} leads us to investigate the asymptotic
properties of $\tilde{A}(\lambda ,x)$, as $\lambda\rightarrow 0$.

\begin {lemma} \label{Lap2}
Let $V(dx)$ be a finite positive Radon measure as in Theorem
\ref{thbase1} and $\tilde{A}$ be the function defined by
(\ref{defAti}).
\begin {enumerate}
 \item
 The measure $(\tilde{A})"(\lambda ,dx)-\tilde{A}(\lambda ,x)V(dx)$
 admits a density function $\theta (\lambda ,x)$ with respect to
 Lebesgue measure and
\begin{equation}\label {pthet}
\lim _{\lambda \rightarrow 0} \Big ( \sup _{x \in \mathbb R}
|\theta (\lambda ,x)| \Big )=0.
\end{equation}
($(\tilde{A})'(\lambda ,x)$ denotes the first $x$-derivative of
 $\tilde{A}(\lambda ,\cdot)$ and
 $(\tilde{A})"(\lambda ,dx)$ the second one, in the sense of distributions).

 \item The $x$-derivative of $\tilde{A}(\lambda ,x)$ is bounded:
\begin{equation}\label {supAprime}
 \sup_{x \in \mathbb R,
\lambda \geq 0} |\tilde{A}'(\lambda ,x)|  <\infty.
\end{equation}
\item \noindent We have:
\begin{equation}\label {limApri}
 \lim_{\lambda\rightarrow 0,
x\rightarrow \pm \infty }(\tilde{A})'(\lambda ,x)=\pm 2 .
\end{equation}
\end{enumerate}
\end{lemma}


 \begin{proof} \ a) It is well known that the function
 $(t,x)\mapsto Z_t^{V}(x)$ is a
 solution in the distribution sense to:
\begin{equation}\label {edpphi}
 \left\{
\begin{array}{ll}
\frac{\partial Z}{\partial t} - \frac{1}{2}\frac{\partial^2
Z}{\partial x^2}
+\frac{1}{2}VZ =0\\\\
Z_0(x)=1,
\end{array}
\right.
\end{equation}
\noindent and that $Z$ can be expressed through the Brownian
motion semigroup $(P_t(x,dy)=p_t(x,y)dy ; t\geq 0)$ :
$$Z_t^{V}(x)=P_t(1)-\frac{1}{2}\int_0^t ds \int_\mathbb R
p_{t-s}(x,y)Z_s^{V}(y) V(dy),$$
\begin{equation}\label {eqphi}
\qquad \qquad \qquad \qquad \qquad \quad=1-\frac{1}{2}\int_0^t ds
\int_\mathbb R \frac{1}{\sqrt{2\pi
(t-s)}}e^{-\frac{(x-y)^2}{2(t-s)}}Z_s^{V}(y)V(dy).
\end{equation}
\noindent We take the Laplace transform in time on both sides;
this yields to
\begin{equation}\label {eqA2}
A(\lambda ,x)=\frac{1}{\lambda}-\frac{1}{2}\int_\mathbb R
A(\lambda ,y)V(dy)\Big (\int_0^\infty \frac{1}{\sqrt{2\pi
v}}e^{-\frac{(x-y)^2}{2v}}e^{-\lambda v} dv \Big ).
\end{equation}
\noi Recall that
$$\int_0^\infty \frac{ds}{\sqrt{s}}\exp\{-\frac{\gamma
s}{2}-\frac{a^2}{2s}\}=
\sqrt{\frac{2\pi}{\gamma}}e^{-|a|\sqrt{\gamma}}; \ \gamma >0, a
\in \mathbb R.$$
\noi Hence
\begin{equation}\label {eqAti}
\tilde{A}(\lambda ,x)=\sqrt{2\lambda}A(\lambda ,x)=
\frac{1}{\sqrt{2\lambda}}\big [ 2-\frac{1}{2}\int_\mathbb R
\tilde{A}(\lambda ,y)e^{-|x-y|\sqrt{2\lambda}}V(dy) \big ],
\end{equation}
\begin{equation}\label {eqAti2}
\sqrt{2\lambda} \tilde{A}(\lambda ,x)=
 2-\frac{1}{2}\int_\mathbb R
\tilde{A}(\lambda ,y)e^{-|x-y|\sqrt{2\lambda}}V(dy).
\end{equation}
\noi Using (\ref{propu}), (\ref{ineg1}) and (\ref{eqAti2})  we
obtain :
\begin{equation}\label {limAti}
\lim_{\lambda\rightarrow0} \Big (\int_\mathbb R \tilde{A}(\lambda
,y)e^{-|x-y|\sqrt{2\lambda}}V(dy)\Big )= \lim_{\lambda\rightarrow
0} \Big (\int_\mathbb R \tilde{A}(\lambda ,y)V(dy)\Big )=4.
\end{equation}
%

 b)\ Let $h$ be a smooth function with compact support. As we
multiply both sides of (\ref {eqAti}) by $h"(x)$, and  integrate
with respect to $dx$, we obtain:
\begin{equation}\label {eqAti2A}
(\tilde{A})"(\lambda,h):=\int_\mathbb R \tilde{A}(\lambda,x) h"(x)
dx=- \frac{1}{2\sqrt{2\lambda}} \int_{\mathbb R ^2}
\tilde{A}(\lambda ,y)e^{-|x-y|\sqrt{2\lambda}}h"(x)dx V(dy),
\end{equation}
\noi where $(\tilde{A})"(\lambda,dx)$ denotes the second
derivative in the distribution sense of $\tilde{A}(\lambda,x)$
with respect to the $x$ variable.

 \noindent
 Let $U^{\lambda}(g)$ be the Brownian $\lambda$-potential of the
function $g$ :

$$U^{\lambda}(g)(x)=E_x\big [ \int _0 ^\infty g(B_s)e^{-\lambda
s} ds \big ]= \frac{1}{\sqrt{2\lambda}} \int _\mathbb R
e^{-|x-y|\sqrt{2\lambda}} g(y) dy.$$

\noindent
 Since $U^{\lambda}(g)$ solves (cf \cite{KaratShrev}):
$$U^{\lambda}(g")(x)=(U^{\lambda}g)"(x)=-2g(x)+2\lambda
U^{\lambda}(g)(x),$$
\noindent
 then
\begin{equation}\label {eqAti3}
(\tilde{A})"(\lambda,h) -\int _\mathbb R
\tilde{A}(\lambda,y)h(y)V(dy) = -\lambda \int _\mathbb R
\tilde{A}(\lambda,y)U^{\lambda}(h)(y)V(dy).
\end{equation}
\noi  This implies that the distribution $(\tilde{A})"(\lambda,dy)
-\tilde{A}(\lambda,y)V(dy)$ is a measure and
$$(\tilde{A})"(\lambda,dy) -\tilde{A}(\lambda,y)V(dy)=\theta
(\lambda,y)dy,$$ \noi where:
\begin{equation}\label {eqtheta1}
\theta (\lambda,y)=-\sqrt{\frac {\lambda}{2}} \int _\mathbb R
\tilde{A}(\lambda,x)e^{-|x-y|\sqrt{2\lambda}}V(dx).
\end{equation}
\noindent Using the inequalities (\ref{ineg1}) and (\ref {propu}),
we obtain :
\begin{equation}
\sup_{x \in \mathbb R} |\theta (\lambda,x)| \leq \sqrt{\frac
{\lambda}{2}} \int _\mathbb R \tilde{A}(\lambda,x) V(dx) \leq
k_1\sqrt{\lambda}.
\end{equation}
\noi This proves part 1. of Lemma  \ref{Lap2} .


\noindent c)  Obviously, (\ref{eqAti}) can be written as follows:
$$\tilde{A}(\lambda ,x)= \frac{1}{\sqrt{2\lambda}}\Big [
2-\frac{1}{2} \Big ( e^{-x\sqrt{2\lambda}}\int_{]-\infty ,x]}
\tilde{A}(\lambda ,y)e^{y\sqrt{2\lambda}}V(dy)$$
 $$+ e^{x\sqrt{2\lambda}}\int_{]x,+\infty[} \tilde{A}(\lambda
,y)e^{-y\sqrt{2\lambda}}V(dy) \Big )\Big ].$$
\noi Taking the  derivatives on both sides with respect to $x$, we
get
\begin{equation} \label{relAtilpri}
(\tilde{A})'(\lambda,x)=\frac{1}{2}\int _{]-\infty ,x]}
\tilde{A}(\lambda,y)e^{-|x-y|\sqrt{2\lambda}} V(dy)-
\frac{1}{2}\int
_{]x,+\infty[}\tilde{A}(\lambda,y)e^{-|x-y|\sqrt{2\lambda}} V(dy).
\end{equation}
\noi Consequently
$$\sup_{x \in \mathbb R} |(\tilde{A})'(\lambda,x)| \leq
\int _\mathbb R \tilde{A}(\lambda,y) V(dy) \leq k_2.$$
%
d) Due to (\ref{relAtilpri}), we have :
\begin{equation} \label{relAtilpri2}
(\tilde{A})'(\lambda,x)=\int _{]-\infty ,x]}
\tilde{A}(\lambda,y)e^{-|x-y|\sqrt{2\lambda}} V(dy)-
\frac{1}{2}\int _\mathbb R
\tilde{A}(\lambda,y)e^{-|x-y|\sqrt{2\lambda}} V(dy).
\end{equation}
\noi  Since (\ref{limAti}) and (\ref{ineg1}) hold,
\begin{equation} \label{relAtilpri3}
\lim_{\lambda \rightarrow 0} \Big (- \frac{1}{2}\int _\mathbb R
\tilde{A}(\lambda,y)e^{-|x-y|\sqrt{2\lambda}} V(dy) \Big )=-2,
\end{equation}
\noi and
$$|\int _{]-\infty ,x]}
\tilde{A}(\lambda,y)e^{-|x-y|\sqrt{2\lambda}} V(dy)|\leq \int
_{]-\infty ,x]} \kappa (1+|y|)V(dy).$$
\noi  As a result, $\displaystyle { \int _{]-\infty ,x]}
\tilde{A}(\lambda,y)e^{-|x-y|\sqrt{2\lambda}} V(dy) }$ goes to $0$
as $x\rightarrow -\infty$, uniformly with respect to $\lambda \geq
0$.

\noindent Consequently $(\tilde{A})'(\lambda,x)$ converges to
$-2$, as $x\rightarrow -\infty, \lambda \rightarrow 0$.

\noindent In the same way,
$$(\tilde{A})'(\lambda,x)=-\int _{]x,+\infty[}
\tilde{A}(\lambda,y)e^{-|x-y|\sqrt{2\lambda}} V(dy)+
\frac{1}{2}\int _\mathbb R
\tilde{A}(\lambda,y)e^{-|x-y|\sqrt{2\lambda}} V(dy),$$
$$\lim_{\lambda \rightarrow 0, x\rightarrow +\infty}(\tilde{A})'(\lambda,x)=2.$$
This ends the proof of Lemma \ref{Lap2}.
\end {proof}


\begin {rem}\label{rem1} If $V(dx)$ has  compact support,
say supp$(V(dx)))\subset [a,b]$, it is easy to check directly:
\begin{equation} \label{relA55}
\lim_{\lambda\rightarrow 0}(\tilde{A})'(\lambda,y)=
-\lim_{\lambda\rightarrow 0}(\tilde{A})'(\lambda,x)=2 ,\  \mbox
{for any } x\leq a \ \mbox { and }\  y\geq b.
\end{equation}
\end {rem}

\begin{prooff}\ {\bf of Remark \ref{rem1}}.


\noi Let $x\leq a$ and $T_a$ be the stopping time : $T_a=\inf \{t
\geq 0; B_t >a\}$.

\noindent We have:
$$E _x \Big [ \exp \big \{-\frac{1}{2} \int  _\mathbb R
L^y _ t V(dy) \big \} \Big]=E _x \Big [ \exp \big \{-\frac{1}{2}
\int  _a^b  L^y _ t V(dy)\big \} \ 1_{\{T_a >t\}}  \Big]$$
$$+E _x \Big [ \exp \big \{-\frac{1}{2}
\int  _a^b  L^y _ t V(dy)\big \}\ 1_{\{T_a \leq t\}}  \Big],$$
$$=P_x(T_a >t)+E _x \Big [ \exp \big \{-\frac{1}{2}
\int  _a^b  L^y _ t V(dy)\big \}\ 1_{\{ T_a \leq t\}}  \Big].$$
\noindent Using the strong Markov property at time $T_a$, and
$\displaystyle {P_x(T_a \in
 ds)= \frac{|x-a|}{\sqrt{2\pi s^3}}e^{-\frac{(x-a)^2}{2s}}\
 1_{\{s>0\}} ds}$ we obtain :
  \begin{equation}\label{idrem}
E _x \Big [ \exp \big \{-\frac{1}{2} \int  _\mathbb R L^y _ t
V(dy) \big \} \Big]=\int _t^\infty \frac{|x-a|}{\sqrt{2\pi
s^3}}e^{-\frac{(x-a)^2}{2s}}\
  ds+ \int _0^t \frac{|x-a|}{\sqrt{2\pi s^3}}e^{-\frac{(x-a)^2}{2s}}\
  Z_{t-s}^{V}(a)ds .
  \end{equation}
  \noi We take the Laplace transform on both sides with respect to time :
  $$A(\lambda ,x)=\frac{1}{\lambda} \big (
  1-e^{(x-a)\sqrt{2\lambda}} \big )
  +A(\lambda ,a)e^{(x-a)\sqrt{2\lambda}} , x\leq a.$$
  Then
  $$\frac{A(\lambda ,x)-A(\lambda ,a)}{x-a}=
  \frac{1-e^{(x-a)\sqrt{2\lambda}}}{x-a} \
  \frac{1-\lambda A(\lambda ,a)}{\lambda}, x<a ,$$
$$(\tilde{A})'(\lambda ,a)=\sqrt{2\lambda} A'(\lambda ,a)=-2+2\lambda
A(\lambda ,a), \ x<a.$$

\noi This proves $\displaystyle \lim_{\lambda\rightarrow
0}(\tilde{A})'(\lambda,x)=-2$ for any $ x\leq a $.
 \noindent If $x\geq b$, we prove by the same way
that $\displaystyle {\lim_{\lambda\rightarrow
0}(\tilde{A})'(\lambda,x)=2}$.

\end{prooff}


\vskip 20 pt \noindent {\bf Proof of Theorem \ref {thbase1}}.

\begin {enumerate}
\item
 The It\^o-Tanaka formula tells us that :

$$M ^\varphi(s)= \varphi(B _ s) \exp \Big \{-{\frac{1}{2}}
 \int  _ \mathbb R L _s^y V_\varphi(dy) \Big \} , s \geq 0$$
 \noi is a continuous local  martingale.

\noi As a consequence of (\ref{maj2C}):
$$M ^{\varphi_V}(s)\leq k(1+\sup_{0\leq u \leq s}|B_u|).$$
\noi Consequently, there exists $\gamma >0$, such that:
$$E\Big [ \exp \big \{ \gamma \big (\sup_{0\leq u \leq s}
M ^{\varphi_V}(u)\big )^2\big \} \Big ] <\infty .$$

 \noi A fortiori $(M ^{\varphi_V}(s);s\geq 0)$
 is a continuous martingale such that $E[M ^{\varphi_V}(s)^2]<\infty$.

\item The function $\tilde{A}(\lambda ,.)$ solves the following
ordinary differential equation, depending on the parameter
$\lambda
>0$:
$$ \left \{ \begin {array}{ll}
(\tilde{A})"(\lambda ,dx)-\tilde{A}(\lambda
,.)V(dx)=\theta(\lambda,x)dx \\
\displaystyle \lim _{x\rightarrow \pm \infty} \Big
((\tilde{A})'(\lambda,x)\Big )=\pm 2 +o(\lambda).
\end{array}
\right .$$
 \noi Property (\ref {pthet}) implies that
$\displaystyle \frac{\tilde{A}(\lambda,x)}{\sqrt{2\pi}}$
converges, as $\lambda \rightarrow 0$, to a function $\varphi_V$,
solution to (\ref{sturm}), and (\ref{sturmcond}). We draw from
this three conclusions :
\begin{enumerate}
\item  $\varphi_V$ is a non-negative function, being
  a limit of non-negative
 functions. Since $\varphi_V$ solves  (\ref {sturm}),
 $\varphi_V$ is a convex function.
\item
 Lemma \ref{Lap1} implies that $\sqrt{t}Z^V _t(x)$ converges
 to $\varphi_V(x)$, as $t\rightarrow\infty$.
 \item From (\ref{limAti}), we have : $\displaystyle \sqrt{2\pi}
 \int_\mathbb{R}\varphi_V(y)V(dy)=4$.
\end{enumerate}

 \item We claim that $\varphi_V$ is strictly positive.
Indeed, (\ref {sturm}) implies that $\varphi_V\not\equiv 0$. As a
result $\varphi_V(B_t)\geq 0$, and $\varphi_V(B_t)$ is not a.s.
equal to $0$. But $(M^{\varphi_V}(s);s\geq 0)$ is a  martingale,
then
$$\varphi_V(x)=E_x\Big [\varphi_V(B _t) \exp \Big \{-{\frac{1}{2}}
\int  _ \mathbb R L _{t}^y V(dy) \Big \} \Big ]>0.$$
\item The proof of the convergence of $Q_{x,t}(\Lambda _s)$ to
$Q_x(\Lambda _s)$ is similar to the one given in section
\ref{Csmall}.

\item Point 5 is a direct consequence of the Girsanov theorem (cf
section \ref{Csmall}).

\item The integral $\displaystyle {\int _\mathbb R
\frac{dy}{\varphi_V^2(y)}}$ is finite because $\varphi_V(y)$ is
equivalent to $k|y|$, as $|y|\rightarrow \infty$.

\noi We remark that $\displaystyle {\beta (x)=\int _0^x
\frac{dy}{\varphi_V^2(y)}}$ is a scale function for the diffusion
process defined by (\ref{defX}) ( \cite{KaratShrev}, Chap. 5,
section 5). Indeed
$$L^{\varphi_V}(\beta)=\frac{1}{2}\beta "+\frac{\varphi_V'}{\varphi_V}\beta '=\frac{1}{2}
\Big ( -\frac{2 \varphi_V'}{\varphi_V^3} \Big
)+\frac{\varphi_V'}{\varphi_V}\Big ( \frac{1}{\varphi_V^2} \Big
)=0. \qquad \qquad \blacksquare$$
\end{enumerate}



\begin{exam} \label{ex1integr}
Let $V(dx)=\gamma ^2 1_{[a,b]}(x)dx$, where $a<b$. Then
$$\lim _{t\rightarrow \infty} \Big (
\sqrt{t}\ E_x \Big [ \exp \big \{-\frac{\gamma^2}{2}\int_0^t
1_{[a,b]} (B_s)ds \big \} \Big ] \Big )=\varphi_V(x),$$
\noi with
 $$\varphi_V(x)= \left
\{
\begin {array}{lll}
\sqrt{\frac{2}{\pi}} \Big (\frac{1}{\gamma \tanh (\gamma
\frac{b-a}{2})}+x-b\Big )
 & {\rm if } & x>b \\
\sqrt{\frac{2}{\pi}} \Big (\frac{\cosh (\gamma
[x-\frac{a+b}{2}])}{\gamma \sinh (\gamma \frac{b-a}{2})}\Big )
& {\rm if } & x \in [a,b ]\\
\sqrt{\frac{2}{\pi}} \Big (\frac{1}{\gamma \tanh (\gamma
\frac{b-a}{2})}+a-x\Big )& {\rm if } & x<a .
\end {array} \right .$$

\end {exam}

\begin{exam} \label{ex2integr}
Let $V(dx)=\gamma ^2 (\delta_a (dx)+\delta_b (dx))$, where $a\leq
b$. Then
$$\lim _{t\rightarrow \infty} \Big (
\sqrt{t}\ E_x \Big [ \exp \big \{-\frac{\gamma^2}{2} (L^a_t+L^b
_t) \big \} \Big ] \Big )=\varphi_V(x),$$
\noi with
 $$\varphi_V(x)= \left
\{
\begin {array}{lll}
\sqrt{\frac{2}{\pi}} \Big (\frac{1}{\gamma^2}+x-b\Big ) & {\rm if } & x>b \\
\sqrt{\frac{2}{\pi}} \frac{1}{\gamma ^2}
& {\rm if } & x \in [a,b ]\\
\sqrt{\frac{2}{\pi}} \Big (\frac{1}{\gamma ^2 }+a-x\Big )& {\rm if
} & x<a.
\end {array} \right .$$
 \noi In particular if $a=b$,
$$\lim _{t\rightarrow \infty} \Big (
\sqrt{t}\ E_x \Big [ \exp \big \{-\gamma^2 L^a_t \big \} \Big ]
\Big )=\sqrt{\frac{2}{\pi}}(\frac{1}{\gamma^2}+|x-a|),$$
\noi and the process $(X_t^x;t\geq 0)$ defined by (\ref {defX})
solves:
$$X_t=x+B_t +\gamma^2\int_0^t\frac{{\rm sgn} (X_s -a)}{1+\gamma^2|X_s -a|}ds.$$
\end {exam}


\subsection {The Ray-Knight theorem and the excursion theory viewpoints
  }\label{int2}

1) Our first approach is based on the Ray-Knight theorem which
gives the law of $(L^y_S;y\in \mathbb{R})$, for $S$  an
exponential r.v. independent of the underlying Brownian motion
$(B_t ;t\geq0)$. We also use the explicit expressions  of Laplace
transforms of certain Bessel quadratic functionals in terms of
solutions to certain Sturm-Liouville equations. For the
convenience of the reader we present the relevant material from
\cite{BiaYor}, \cite{PitYor}, \cite{RevYor} without proofs, thus
making our exposition self-contained.


\begin{defi}\label{def10integ}
\begin{enumerate}
    \item
    Let $f : \ \mathbb{R}\mapsto \mathbb{R}$ and $v$ be a positive measure on
     $\mathbb{R}$ . We denote
    by $<f,v>$ or $v(f)$ the integral of  $f$ with respect to $v$, namely :
    $$
    <f,v> =v(f)=\int_\mathbb{R}f(t)v(dt) .
    $$
    \noi We set : $v_+(dx)=1_{\{x>0\}}v(dx)$ and
    $v_-(dx)$ the image of $1_{\{x<0\}}v(dx)$ by the map $x\mapsto
    -x$.

    \item
Let $\delta\geq 0$. We define $Q_x^{(\delta)}$, the distribution
of the square of the $\delta$-dimensional Bessel process, started
at $x$.
 \end{enumerate}
\end{defi}

\noi We now present some important properties of the family
$(Q_x^{(\delta)})$

\begin{prop}\label{prop21integ}
\begin{enumerate}
    \item The family of probability measures $(Q_x^{(\delta)}; \delta,x\geq
    0)$ obeys the additivity property :
    \begin{equation}\label{100integ}
Q_x^{(\delta)}\star Q_{x'}^{(\delta')}=Q_{x+x'}^{(\delta+\delta')}
,\qquad \delta, \delta',x,x' \geq 0.
\end{equation}
\item If $\lambda(ds)$ is a positive Radon measure  on $\mathbb
R_+$, with finite first moment, then \cite{PitYor}:
\begin{equation}\label{10integ}
Q_x^{(\delta)} \big [ \exp \big\{-\int _0^{+\infty}
Y(s)\lambda(ds) \big\}\big]= Q_x^{(\delta)} \big [ \exp
\big\{-<Y,\lambda >\big\}\big]= \exp \big\{ -xM(\lambda)-\delta
N(\lambda)\big\},
\end{equation}
\noi where $x\geq 0$, $(Y(s);s\geq 0)$ denotes the canonical
process on ${\cal C}(\mathbb{R}_+)$ ($Y(s)(\omega)=\omega (s)$),
$M$ and $N$ are the two positive $\sigma$-finite measures on
${\cal C}(\mathbb{R}_+)$ (\cite{PitYor},\cite{Pitman}) which allow
to express the L\'{e}vy-Khintchine representation of any
$Q^{(\delta)} _x$, i.e.
 one has :
$$M(\lambda )=\int (1-e^{-<\lambda , \omega >})M(d\omega ) , \quad
N(\lambda )=\int (1-e^{-<\lambda , \omega >})N(d\omega ).$$
\item
 \noi Introducing  $\phi_\lambda$  the unique solution of :
 \begin{equation}\label{101integ}
    \frac{1}{2}\phi"=\lambda \phi \quad \mbox{on} \quad
    (0,\infty); \quad \phi(0)=1, \ 0\leq \phi \leq 1,
\end{equation}
we have \cite[Theorem (2.1)]{PitYor}:
\begin{equation}\label{102integ}
    Q_x^{(\delta)} \big [ \exp \big\{-\int _0^{+\infty}
Y(s)\lambda(ds)
\big\}\big]=\big(\phi_{\lambda}(\infty)\big)^{\delta/2} \exp
\big\{\frac{x}{2}\phi'_{\lambda}(0)\big\},
\end{equation}
\noi where $\phi_{\lambda}(\infty)$ and $\phi'_{\lambda}(0)$ are
respectively the limit at $\infty$, and the right derivative at
$0$ of $\phi_\lambda$.

\noi Comparing (\ref{10integ}) and (\ref{102integ}), we have :
\begin{equation}\label{103integ}
    M(\lambda)=-\frac{1}{2}\phi'_{\lambda}(0), \qquad
    N(\lambda)=-\frac{1}{2}\ln (\phi_{\lambda}(\infty)).
\end{equation}
\end{enumerate}
\end{prop}

\noi Our approach is  based on the  knowledge of the law of
$(L^y_S;y\in \mathbb{R})$, for $S$  an exponential r.v.
independent of  $(B_t ;t\geq0)$. This distribution is given in
Proposition \ref{prop22integ} below, through the family of
measures $\big(P^{(\theta)}_{a,l}; a\in\mathbb{R},l>0\big)$
defined in the next Definition \ref{def11integ}.

\begin{defi}\label{def11integ} Let $l\geq 0$ and  $a>0$. We define $P^{(\theta)}_{a,l}$ to
be  the unique probability measure on ${\cal C}(\mathbb{R})$ such
that :
\begin{enumerate}
    \item
    $Y_0=l,$
    \item
    $(Y_{-t};t\geq 0)$ is the diffusion process with infinitesimal
    generator : $\displaystyle
    {2y\frac{d^2}{d y^2}-2\theta y\frac{d}{dy}}$,
    \item
$(Y_{t};0 \leq t\leq a)$ is the diffusion process with
infinitesimal generator : $\displaystyle
   {2y\frac{d^2}{d y^2}-2\theta y\frac{d}{dy}+2\frac{d}{dy}}$,
    \item
$(Y_{t}; t\geq a)$ is the diffusion process with infinitesimal
generator : $\displaystyle
   {2y\frac{d^2}{dy^2}-2\theta y\frac{d}{dy}}$,
\end{enumerate}
\noi When $a<0, \ (Y_{t}; t\in \mathbb{R})$ is distributed under
$P^{(\theta)}_{a,l}$ as $(Y_{-t}; t\in \mathbb{R})$  under
$P^{(\theta)}_{-a,l}$.

\end{defi}
\begin{prop}\label{prop22integ}
Suppose that $V(dy)$ is a positive measure on
     $\mathbb{R}$ satisfying (\ref{propu}), $S_\theta$ is an exponential
     r.v.with parameter $\theta^2/2$ (i.e. with expectation
     $2/\theta^2$) and independent of $(B_t ;t\geq 0)$. Then (cf \cite[theorem
     1]{BiaYor}):
     \begin{equation}\label{11integ}
    E_0\Big[ \exp -\big \{\int _\mathbb{R}L^y_{S_\theta}V(dy)\big
    \}\Big]=\theta \int_0^\infty e^{-\theta l}  \Big(\frac{1}{2} \int _\mathbb{R}\theta
    e^{-\theta |a|} P^{(\theta)}_{a,l}\big[e^{-<Y,V>} \big] da
    \Big)dl .
\end{equation}
\end{prop}

\vskip 15 pt \noi We are now able to state the main result of this
subsection.

\begin{prop}\label{prop1integ}
Suppose that the positive measure $V(dy)$ has  compact support.
\begin{enumerate}
    \item
Then :
\begin{equation}\label{12integ}
    \lim_{\theta\rightarrow 0}\Big(\frac{1}{\theta}
    E_0\Big[ \exp -\big \{\int _\mathbb{R}L^y_{S_\theta}V(dy)\big
    \}\Big] \Big)= H(V),
\end{equation}
\begin{equation}\label{12Ainteg}
    \lim_{t\rightarrow \infty}\Big(\sqrt{t}
    E_0\Big[ \exp -\big \{\int _\mathbb{R}L^y_t V(dy)\big
    \}\Big] \Big)= \sqrt{\frac{2}{\pi}}H(V),
\end{equation}
where $H(V)$ is defined as :
\begin{equation}\label{13integ}
    H(V)=\frac{1}{2}\int _0^{+\infty}
    \Big(Q_l^{(0)}\big[e^{-<Y,V_->}\big]Q_l^{(2)}\big[e^{-<Y,V_+>}\big]
    +Q_l^{(2)}\big[e^{-<Y,V_->}\big]Q_l^{(0)}\big[e^{-<Y,V_+>}\big]\Big)dl.
\end{equation}
\noi In terms of $M(V_\pm)$ and $N(V_\pm)$ (resp.
$\phi_{\lambda_\pm}(\infty), \ \phi'_{\lambda_\pm}(0)$) , we have
:
\begin{equation}\label{14integ}
    H(V)=\frac{1}{2(M(V_+)+M(V_-))}\big(e^{-2N(V_+)}+e^{-2N(V_-)}\big)
    = \frac{\phi_{\lambda_+}(\infty)+\phi_{\lambda_-}(\infty)}
{\phi'_{\lambda_+}(0)+\phi'_{\lambda_-}(0)}.
\end{equation}
\item In particular if $V(dx)$ is a symmetric measure (i.e.
$V(dx)$ coincides with its image by the map $x\mapsto -x$ ), then
:
\begin{equation}\label{15integ}
    H_{\mbox{sym}}(V):=H(V)=\frac{1}{2}\int _0^{+\infty}Q_l^{(2)}
    \big[e^{-<Y,V_+>}\big]
    \ dl=\frac{1}{2M(V_+)}e^{-2N(V_+)}=\frac{\phi_{\lambda_+}(\infty)}
    {\phi'_{\lambda_+}(0)}.
\end{equation}
\end{enumerate}
\end{prop}

 \begin{prooff} \ {\bf of Proposition
\ref{prop1integ}}. We give two proofs of Proposition
\ref{prop1integ}; the first one uses the Ray-Knight theorem for
Brownian local times up to an exponential time (cf Proposition
\ref{prop22integ}); the second one uses excursion theory.

\noi  {\bf First proof of Proposition \ref{prop1integ}}. 1) We set
:
$$
\Delta=E_0\Big[ \exp -\big \{\int
_\mathbb{R}L^y_{S_\theta}V(dy)\big
    \}\Big].
    $$
\noi Relation (\ref{11integ}) implies that $\Delta$ may be written
as follows :
\begin{equation}\label{16integ}
    \Delta= \theta \int_0^\infty e^{-\theta l}  \Big(\frac{1}{2} \int _\mathbb{R}
    e^{- |b|} \Delta _\theta(b,l) db
    \Big)dl ,
\end{equation}
\noi where :
$$
\Delta _\theta(b,l)=P^{(\theta)}_{b/\theta,l}\big[e^{-<Y,V>}\big].
$$
\noi Suppose that $b>0$. We decompose $<Y,V>$ in the following way
:
$$
<Y,V>=\int_{-\infty}^0 Y_z V(dz)+\int_{0}^{b/\theta} Y_z V(dz)+
\int^{+\infty}_{b/\theta} Y_z V(dz).
$$
\noi Using Proposition \ref{prop22integ} and taking the limit,
$\theta\rightarrow0$ we obtain :
$$
\lim_{\theta\rightarrow 0} \Delta _\theta (b,l)=
Q_l^{(0)}\big[e^{-<Y,V_->}\big]Q_l^{(2)}\big[e^{-<Y,V_+>}\big].
$$
\noi Moreover $0 \leq \Delta _\theta (b,l) \leq 1$. Applying the
same reasoning  to the case $b<0$ we obtain :
$$
\lim_{\theta\rightarrow 0} \Big(\frac{1}{2} \int _\mathbb{R}
    e^{- |b|} \Delta _\theta(b,l) db \Big)=H(V),
$$
\noi where $H(V)$ is defined by (\ref{13integ}).

\noi Identity  (\ref{14integ}) follows directly from
 (\ref{10integ}).

2) We  now turn to the symmetric case. Since $V_-=V_+$, additivity
property (\ref{100integ}) directly implies that :
$$
Q_l^{(0)}\big[e^{-<Y,V_->}\big]Q_l^{(2)}\big[e^{-<Y,V_+>}\big]=
Q_{2l}^{(2)}\big[e^{-<Y,V_+>}\big].
$$
\noi This proves (\ref{15integ}).

3) Since $S_\theta$ is independent from $B$, it follows that :
$$
E_0\Big[ \exp -\big \{\int _\mathbb{R}L^y_{S_\theta}V(dy)\big
    \}\Big] =\frac{\theta^2}{2}\int_0^\infty \Big( E_0\Big[ \exp -\big \{
    \int _\mathbb{R}L^y_t V(dy)\big
    \}\Big] \Big)e^{-\theta^2 t/2} dt.
    $$
    \noi As $\displaystyle {t\mapsto E_0\Big[ \exp -\big \{
    \int _\mathbb{R}L^y_t V(dy)\big
    \}}\Big]$ is decreasing, we conclude  from (\ref{12integ})
    and the Tauberian theorem (cf \cite{Fell}, Chap.
XIII, section 5 ) that (\ref{12Ainteg}) holds.

\noi {\bf Second proof of Proposition \ref{prop1integ}}.

\noi  We suppose for simplicity that $V$ is a positive function
with compact support. We start as in the previous approach
considering $S_\theta$  an exponential r.v.with parameter
$\theta^2/2$ (i.e. with expectation $2/\theta^2$) and independent
of $(B_t ;t\geq 0)$. We again consider :
%
$$
\Delta= E_0\Big[ \exp -\big \{\int
_\mathbb{R}L^y_{S_\theta}V(y)dy\big
    \}\Big].
    $$
    %
%
%

 \noi We express
$\Delta$ with the help of excursion theory
%
%
\begin{lemma}\label{excurs001} $\Delta$ is equal to the ratio
$\displaystyle \frac{{\cal N}^{(\theta )}}{{\cal D}^{(\theta )}}$,
where,
$${\cal N}^{(\theta )}=\frac{\theta ^2}{2}\int \Big[
\int_0^{\zeta (\varepsilon)}\exp \Big \{-\frac{\theta
^2t}{2}-\frac{1}{2}\int_0^t V(\varepsilon _s)ds \Big\} \ dt\Big]
n(d\varepsilon ),
$$
$${\cal D}^{(\theta )}= \int \Big[
1-\exp\Big\{-\frac{\theta ^2}{2}\zeta (\varepsilon)
-\frac{1}{2}\int_0^{\zeta (\varepsilon)}V(\varepsilon _s)ds\Big\}
\Big]n(d\varepsilon ) ,
$$
$n(d\varepsilon )$ denotes It\^{o}'s measure of excursions and
$\zeta (\varepsilon )=\inf \{ s>0 ; \varepsilon _s=0\}$.
\end{lemma}

\begin{prooff} \ {\bf of Lemma \ref{excurs001}.} It follows easily
from the general integral representation formula ( cf
\cite{RevYor}, Exercise 4.18, Chap. XII):
\begin{equation}\label{excur05}
    \int _0^\infty P_0^t dt=\int _0^\infty P_0^{\tau_l}dl \ \circ
    \ \int _0^\infty n^u(\cdot \cap \{u < \zeta \})du,
\end{equation}
\noi where, for any random time $T$, $P_0^T$ denotes the Wiener
measure restricted to the $\sigma$-field ${\cal F}_T$, $n^u$
denotes the It\^{o} measure restricted to the corresponding
$\sigma$-field ${\cal F}^\star _u=\sigma(\varepsilon_s; 0 \leq s
\leq u)$  for excursions $\varepsilon$, and $\circ$ indicates the
concatenation, operation acting on measures on path space (see
\cite{RevYor}, Chap. XII, section 4, for details). Finally
 $(\tau_l ; l\geq 0)$ is the inverse local time at $0$.

\noi As a consequence of (\ref{excur05}), we get
$\Delta=\Delta_-\Delta_+$ where :
$$
\Delta_+=\frac{\theta^2}{2}\int _0 ^\infty dt \Big(\int \exp
\Big\{ -\frac{\theta^2}{2}t-\int_0^tV(\varepsilon_s)ds
\Big\}1_{\{t <\zeta (\varepsilon )\}} n(d\varepsilon )\Big),
$$
$$
\Delta_-=\int_0^\infty dl E_0[\exp
\{-\int_0^{\tau_l}V(B_s)ds\}].$$
\noi Using Fubini's theorem, we find : $\Delta_+={\cal
N}^{(\theta)}$; concerning $\Delta_-$, we get from excursion
theory (cf \cite{RevYor}, Proposition (2.7), Chap XII) :
%
$$
 E_0[\exp
\{-\int_0^{\tau_l}V(B_s)ds\}]=\exp \Big\{ - l\int n(d\varepsilon)
\Big( 1-\exp\{-\int_0^{\zeta (\varepsilon )}V(\varepsilon _s)ds\}
\Big)\Big\},
$$
\noi and consequently $\Delta_-=1/{\cal D}^{(\theta )}$.
\end{prooff}

 \noi Let us provide now the second proof of Proposition
 \ref{prop1integ}.

 \noi
As $\theta \rightarrow 0$, the denominator ${\cal D}^{(\theta )}$
tends to :
$$\int \Big[
1-\exp\Big\{ -\frac{1}{2}\int_0^{\zeta (\varepsilon)}V(\varepsilon
_s)ds\Big\} \Big]n(d\varepsilon ) = -\frac{1}{2}(\phi'_{V_+}(0_+)+
\phi'_{V_-}(0_+)).
$$
Let us consider the numerator, which we may write as :
$${\cal N}^{(\theta )}=\frac{\theta ^2}{2}
\int_0 ^\infty e^{-\theta ^2 t/2} \Big[ \int 1_{\{ \zeta
(\varepsilon)>t\}}\exp \Big \{-\frac{1}{2}\int_0^t V(\varepsilon
_s)ds \Big\} n(d\varepsilon ) \Big] dt.
$$
Now recall that (Ex 4.18, Chap XII in \cite{RevYor}) :
$$
1_{\{ \zeta (\varepsilon)>t\}}  n_{\pm}(d\varepsilon
):=\frac{1}{2} \frac{1}{\sqrt{2\pi t}}M^t(d\varepsilon ),
$$
where $M^t$ denotes the law of the Brownian meander with length
$t$.

\noi Thus, we find :
$${\cal N}^{(\theta )}=\frac{\theta ^2}{2}
\int_0 ^\infty \frac{dt}{\sqrt{2\pi t}} e^{-\theta ^2 t/2}
\frac{1}{2}\Big\{M^t\Big( \exp \{-\frac{1}{2}\int_0^t
V_+(\varepsilon _s)ds \}\Big)$$
$$
\qquad\qquad\qquad\qquad\qquad\qquad + M^t\Big( \exp
\{-\frac{1}{2}\int_0^t V_-(\varepsilon _s)ds \}\Big) \Big\},
$$
(recall that  $V_+$ and $V_-$ are defined by the rules given in
Definition \ref{def10integ}).

\noi For simplicity, we now write $\alpha =\theta ^2 /2$, and we
make the change of variables : $\alpha t=u$; then :
$${\cal N}^{(\theta )}=\sqrt{\alpha}
\int_0 ^\infty \frac{du}{\sqrt{2\pi u}} e^{-u}
\frac{1}{2}\Big\{M^{u/\alpha}\Big( \exp
\{-\frac{1}{2}\int_0^{u/\alpha} V_+(\varepsilon _s)ds \}\Big)$$
\begin{equation}\label{16A1integ}
    \qquad\qquad\qquad\qquad\qquad\qquad + M^{u/\alpha}\Big( \exp
\{-\frac{1}{2}\int_0^{u/\alpha} V_-(\varepsilon _s)ds \}\Big)
\Big\}.
\end{equation}
Now we use the fact (cf again Ex 4.18, Chap XII in \cite{RevYor})
that :
$$
M^t=\sqrt{\frac{\pi}{2}}\frac{\sqrt{t}}{R_t}P^{(3)}_{0 |{\cal
F}_t},$$
where $P^{(3)}_0$ denotes the law of the three dimensional Bessel
process started at $0$.

\noi It is not difficult to show that :
\begin{equation}\label{excur01}
    \lim _{t\rightarrow \infty }M^t\Big( \exp
\{-\frac{1}{2}\int_0^t V_\pm(\varepsilon _s)ds \}\Big)
=E_0^{(3)}[\exp \{-\frac{1}{2} < V_\pm ,Y > \}],
\end{equation}

where  under $P^{(3)}_0$, $Y$ stands for the three dimensional
Bessel process starting from $0$, and
$$ < V_\pm ,Y
>=\int_0^\infty Y_sV_\pm(s)ds .$$
\noi Hence, from (\ref{16A1integ}), we deduce :
$${\cal N}^{(\theta )}\sim \frac{\theta}{2}
E_0^{(3)}\Big[\exp \{-\frac{1}{2} < V_+ ,Y > \}+ \exp
\{-\frac{1}{2} < V_- ,Y
> \}\Big] , \quad \theta\rightarrow 0 ,$$
and, from the Ray-Knight Theorem for the three dimensional Bessel
process, the right hand-side of (\ref{excur01}) is :
\begin{equation}\label{excur02}
\frac{\theta}{2}\Big\{Q^{(2)}_0[\exp \{-\frac{1}{2} < V_+ ,Y> \}]
+ Q^{(2)}_0[\exp \{-\frac{1}{2} < V_+ ,Y>
\}]\Big\}=\frac{\theta}{2}
(\phi_{V_+}(\infty)+\phi_{V_-}(\infty)),
\end{equation}
from, e.g.  \cite{PitYor}. Hence, we have proven (\ref{12integ}).

\end{prooff}

\begin{rem} 1)Donati-Martin and Hu \cite{DoMaHu} prove the
convergence in law for the  Wiener measure perturbed by the
exponential martingale density associated with $\displaystyle
\int_0^t\frac{dB_s}{B_s}1_{\{|B_s| \geq \varepsilon\}}$, as
$\varepsilon\rightarrow 0$; the limiting law is that of the
symmetrized $BES(3)$ process, i.e. : a process taking values in
$\mathbb{R}_+$ with probability $1/2$, and in $\mathbb{R}_-$ with
probability $1/2$.

\noi 2) Let  $V$ be a function with compact support. Y. Hu
(private communication) has studied the asymptotic behaviour of
$Z_t^V(x)$, as $t\rightarrow\infty$, where in (\ref{1Aintro}),
$(X_t;t\geq0)$ is a one-dimensional diffusion process . With some
additional assumptions, using excursion theory for Brownian
motion, Y. Hu has determined the rate of decay of $Z_t^V(x)$, as
$t\rightarrow\infty$. In particular, Y. Hu has recovered the
result concerning Bessel processes with dimension $0<d<2$.

 \end{rem}


\section {The unilateral case} \label{casunil}
\setcounter{equation}{0}

 In this section the given  positive Radon
measure $V(dy)\not = 0$ on $\mathbb R$, is supposed to be strongly
asymmetric : it is "small" at $-\infty$ and "big" at $+\infty$.
More precisely we suppose :
\begin{equation}\label{prop2-u}
\int ^0_{-\infty} (1+|y| ) V(dy)<\infty.
\end{equation}
\begin{equation}\label{prop2+u}
\liminf_{x\rightarrow +\infty} \Big ( x^{2\alpha} V_a(x) \Big )>0,
{\rm for \ some }\ \alpha <1,
\end{equation}
\noi where $V(dy)=V_a(y)dy+V_s(dy)$ is the Lebesgue decomposition
of $V(dy)$, and in (\ref{prop2+u}) the liminf may be equal to
$+\infty$

\noi
 We remark that if $\displaystyle \lim_{x\rightarrow +\infty} \Big ( x^{2\alpha} V_a(x) \Big
 )$ exists for \ some $\ \alpha >1$,   then $V$ fulfills
(\ref{propu}). This case has been studied in the previous section.

\noi Let us state the main result of this section.

\vskip 10 pt \noi
\begin{theorem} \label{thunil}
Let $V(dy)$ be a positive Radon measure on $\mathbb R$ fulfilling
(\ref{prop2-u}) and (\ref{prop2+u}).
\begin{enumerate}
 \item
 The generic Theorem applies with $k=1/2$, i.e. :
$$\lim _{t\rightarrow\infty} \Big(\sqrt {t}\ E _x\Big [ \exp \Big \{-{\frac{1}{2}}
\int  _ \mathbb R L _t^y V(dy) \Big \} \Big ]\Big):=\varphi_V(x) \
\mbox{exists in } \mathbb{R}.$$

\item
 $\varphi_V$ is a convex function which takes its values in
$]0,\infty[$ and is the unique solution to the Sturm-Liouville
equation :
\begin{equation}\label{sturmunil4}
\varphi"(dx)=\varphi(x)V(dx),
\end {equation}
\noi with boundary conditions:
\begin{equation}\label{sturmunilboun}
\lim  _ {x\rightarrow -\infty}\varphi_V'(x) =-\sqrt {
\frac{2}{\pi} } \qquad ; \quad \lim  _ {x\rightarrow +\infty}
\varphi_V(x)=0.
\end {equation}
 \noi
 Moreover there exist two positive constants    $C,C'$ such that
 \begin{equation}\label{majC+unil}
  \varphi_V(x)\leq C(1+|x|); \quad x \leq 0,
\end{equation}
\begin{equation}\label{majC-unil}
  \varphi_V(x)\leq C e^{-C'x^{1-\alpha}};\quad  x \geq 0 .
\end{equation}
\item
 Let $M^\varphi$ be the process defined by (\ref{maj2CA}), then
 $(M^\varphi_t; t\geq 0)$ is a continuous martingale.

 \item
 Let  $(X_t ^x;t\geq 0)$ be the solution to (\ref{defX}), then
  the law of $(X_t ^x;t\geq 0)$ is $P^{\varphi_V} _ x$,
$(X_t ^x;t\geq 0)$ is transient, i.e.
\begin{equation} \label{limXunil}
P\Big ( \lim  _{t \rightarrow \infty} X _t^x=-\infty \Big ) = 1.
\end{equation}
\end{enumerate}
\end{theorem}

\noi Our proof of Theorem \ref{thunil} consists of two main steps.
We begin by establishing  an a priori upper-bound for $t\mapsto
\sqrt{t}Z_t^V(x)$ (cf Lemma \ref{etape1unil}) . In a second step
we show that we may reduce the discussion to the case where $V$
has a compact support.
\begin {lemma}\label{etape1unil}
Let $\tilde{\varphi}$ be the function defined as follows:
\begin{equation} \label{defCtiunil}
\displaystyle \tilde{\varphi}(x):=\limsup_{ t\rightarrow\infty}
\Big (\sqrt {t} \ E _x\Big [ \exp \big \{-{\frac{1}{2}} \int  _
\mathbb R L _t^y V(dy) \big \} \Big ]\Big ).
\end{equation}
\noi Then  there exist two positive numbers $C,C'$ such that
\begin {equation} \label{etudeCtil}
\tilde{\varphi}(x) \leq Ce^{-C'x^{1-\alpha}}, \mbox{ for any }
x\geq 0.
\end{equation}
\noi In particular
\begin{equation} \label{limCtiunil}
 \displaystyle \lim  _{x \rightarrow +\infty} \tilde{\varphi}(x)=0.
\end{equation}
\end {lemma}

\begin{prooff} \ {\bf of Lemma \ref{etape1unil}}.
Assumption (\ref{prop2+u}) implies that there exist $\kappa, a>0$
such that
$$\int  _ \mathbb R L _t^y V(dy)\geq \frac{1}{2}\frac{\kappa^2}{b^{2\alpha}}\int_a ^b L _t^y dy,
\quad b>a .$$
\noi Then for any $y\in [a,b]$, Example \ref{ex1integr} implies
that
$$\displaystyle \limsup_{ t\rightarrow\infty} \Big
(\sqrt {t} \ E _y\Big [ \exp \big \{-{\frac{1}{2}} \int  _ \mathbb
R L _t^y V(dy) \big \} \Big ]\Big ) \leq \sqrt{\frac{2}{\pi}}
\frac{\cosh \Big[\frac{\kappa}{b^\alpha}(y-\frac{a+b}{2})\Big]}
{\frac{\kappa}{b^\alpha}\sinh \Big
[\frac{\kappa}{b^\alpha}\frac{b-a}{2}\Big]}.$$
\noi Let $y>a$, we choose $b=2y-a$. This brings
\begin{equation}
\tilde{\varphi}(y)\leq \frac{1}{\kappa}
\sqrt{\frac{2}{\pi}}\frac{(2y-a)^\alpha}{\sinh
\Big[\kappa\frac{y-a}{(2y-a)^\alpha}\Big]},
\end{equation}
\noi which proves (\ref{etudeCtil}).
\end{prooff}

\begin{lemma}\label{etape2unil}
Let  $Z^V$ be the function :
\begin{equation}\label{defphiunil}
Z^V _t(x)=E _x\Big [ \exp \big \{-{\frac{1}{2}} \int  _ \mathbb R
L _t^z V(dz) \big \} \Big ], t\geq 0 , x \in \mathbb R .
\end{equation}
\noi Then, if $y>\max\{0,x\}$, $Z^V _t(x) =Z _1(t,y;x)+Z_2
(t,y;x)$ where $Z _1(t,y;x), \ Z_2 (t,y;x)$ are two non-negative
functions and
\begin{equation}\label{inegphi1unil}
$$\displaystyle \limsup_{ t\rightarrow\infty} \Big
(\sqrt {t} \  Z _1(t,y;x)\big \} \leq 2 \tilde{\varphi}(y),
\end{equation}
\begin{equation}\label{inegphi2unil}
\sqrt {t}\  Z _2(t,y;x) \mbox{ converges as } t\rightarrow\infty ,
\end{equation}
\noi where $\tilde{\varphi}$ is the function defined by
(\ref{defCtiunil}).

\end{lemma}

\begin{prooff} \ {\bf of Lemma \ref{etape2unil}}. We decompose :
$$ Z^V _t(x)=Z _1(t,y;x)+Z_2 (t,y;x),$$
$$\displaystyle
Z _1(t,y;x)=E _x\Big [ \exp \big \{-\frac{1}{2} \int  _ \mathbb R
L _t^z V(dz) \big \} \ 1_{\{T_y<t\}} \Big ],$$
$$\displaystyle
Z _2(t,y;x)=E _x\Big [ \exp \big \{-\frac{1}{2} \int  _ \mathbb R
L _t^z V(dz) \big \} \ 1_{\{T_y\geq t\}} \Big ],$$
\noi where $T_y=\inf \{t\geq 0 ;B_t =y\}$.

a) We start with the study of $Z _1(t,y;x)$.

\noi Using the strong Markov property at time $T_y$, we get:
$$Z _1(t,y;x)= E _x\Big [1_{\{T_y<t\}}\exp \big \{-\frac{1}{2} \int  _
\mathbb R L _{T_y}^z V(dz) \big\} \ Z^V _{t-T_y}(y)\Big].$$
\noi Since $V(dx)$ is a positive measure, then
\begin{equation}\label{ineg1phi1}
\sqrt{1+t}\ Z _1(t,y;x)\leq E _x\Big [1_{\{T_y<t\}}
\frac{\sqrt{1+t}}{\sqrt{1+t-T_y}}\ \sqrt{1+t-T_y}\ Z_{t-T_y}^V
(y)\Big ].
\end{equation}

\vskip 5 pt  i) We claim that for fixed $y$
\begin{equation}\label{familleuni}
\Big \{  \frac{\sqrt{1+t}}{\sqrt{1+t-T_y}}1_{\{T_y<t\}} ; \  t\geq
1\Big \} \mbox{ is uniformly integrable}.
\end{equation}
\noi It suffices to prove that $\displaystyle
\frac{\sqrt{1+t}}{\sqrt{1+t-T_y}}1_{\{T_y<t\}}$
 are bounded r.v.'s  in $L^2$, uniformly with respect to
$t\geq 1$.

\noi We have:
$$E \Big [  \frac{1+t}{1+t-T_y}1_{\{T_y<t\}}\Big ]=\int_0^t
a_t(s)ds,$$
 \noi where:
$$a_t(s)=\frac{1+t}{1+t-s} \frac{y}{\sqrt{2\pi s^3}}
\exp \{-\frac{y^2}{2s}\}.$$
\noi We distinguish two cases:

$\alpha$ ) $s\in  [0,t/2]$, then $\displaystyle \frac{1+t}{1+t-s}
\leq 2$, hence $\displaystyle a_t(s) \leq \frac{2y}{\sqrt{2\pi
s^3}} \exp \{-\frac{y^2}{2s}\}$ and
$$\int_{0}^{t/2} a_t(s)ds \leq 2,$$
since $\displaystyle s \mapsto \frac{y}{\sqrt{2\pi s^3}} \exp
\{-\frac{y^2}{2s}\}$ is a  density function.

  $\beta$) $s\in  [t/2,t]$, then $\displaystyle a_t(s) \leq
\frac{1+t}{\sqrt{2\pi (t/2)^3}}\frac{y}{1+t-s}$ and
 $$\int_{t/2}^{t} a_t(s) ds\leq \frac{2 y}{\sqrt{2\pi}}
 \frac{1+t}{\sqrt{ t^3}}\ln (1+t/2).$$
\noi Finally
$$\displaystyle
\sup_{t \geq 1}\Big (\int_0^t a_t(s)ds\Big )<\infty.$$
 \noi This proves (\ref{familleuni}).

ii) The definition of $\tilde{\varphi}$(cf (\ref {defCtiunil}))
implies the existence of a positive number $a$ (depending on $y$)
such that:
\begin{equation} \label{ineg3phiunil}
\sqrt{1+t}Z^V_t (x)\leq 2\tilde{\varphi}(y), \quad \mbox{for any }
t\geq a.
\end{equation}
 \noi  On the right hand-side of
(\ref{ineg1phi1}) the decomposition of  $\{T_y <t\}$ as the
disjoint union of $\{t-T_y>a\}$ and $\{t-a \leq T_y<t\}$, leads to
$$\sqrt{1+t}\ Z _1(t,y;x)\leq
Z_{1,1} (t,y;x)+Z_{1,2} (t,y;x),$$
$$Z_{1,1}(t,y;x)=E _x\Big [1_{\{t-T_y>a\}} \frac{\sqrt{1+t}}{\sqrt{1+t-T_y}}\
\sqrt{1+t-T_y}\ Z_{t-T_y}^V (y)\Big ],$$
$$Z_{1,2}(t,y;x)=E _x\Big [1_{\{t-a\leq T_y <t\}} \frac{\sqrt{1+t}}{\sqrt{1+t-T_y}}\
\sqrt{1+t-T_y}\ Z_{t-T_y}^V (y)\Big ].$$
 \noi The inequality (\ref{ineg3phiunil}) and the property (\ref{familleuni}) imply
$$\displaystyle \limsup_{ t\rightarrow\infty} \Big (
Z_{1,1}(t,y;x) \Big ) \leq 2\tilde{\varphi}(y).$$
\noi As for $Z_{1,2}(t,y;x)$, the function $Z^V$ being less than
$1$,
$$Z_{1,2}(t,y;x) \leq \sqrt{1+t}\ P(t-a\leq T_y <t)=
\frac{\sqrt{1+t}}{\sqrt{2\pi}}\int_{t-a}^te^{-y^2/2s}\frac{ds}{s^{3/2}},$$
$$Z_{1,2}(t,y;x) \leq
\frac{\sqrt{1+t}}{\sqrt{2\pi}}\frac{a}{(t-a)^{3/2}}.$$
\noi Consequently:
$$\displaystyle \limsup_{ t\rightarrow\infty} \Big (
Z_{1,2}(t,y;x) \Big )=0, $$
\noi hence, finally :
\begin{equation}\label {ineg4phi1unil}
$$\displaystyle \limsup_{ t\rightarrow\infty} \Big (
Z_1(t,y;x) \Big ) \leq 2\tilde{\varphi}(y).
\end{equation}

\vskip 10 pt  b) We now prove  (\ref{inegphi2unil}).

\noi Recall that $y>x$. The key observation is the following : on
$\{T_y\geq t\}$, $V(dz)$ can be replaced by $V^{(y)}(dz)$ where:
$$V^{(y)}(dz)=1_{]-\infty ,y]}(z)V(dz),$$
\noi which allows us to reduce the discussion to the integrable
case since $V^{(y)}(dz)$ satisfies (\ref{propu}).

\noi More precisely, we have:
$$Z _2(t,y;x)=E _x\Big [ \exp \big \{-\frac{1}{2} \int  _
\mathbb R L _t^z V^{(y)}(dz) \big \}  \ 1_{\{T_y \geq t\}} \Big
].$$
Hence $ Z _2(t,y;x)=Z _{2,1}(t,y;x) -Z_{2,2}(t,y;x),$ with
$$Z _{2,1}(t,y;x)=E _x\Big [ \exp \big \{-\frac{1}{2} \int  _
\mathbb R L _t^z V^{(y)}(dz) \big \}  \Big ],$$
and
$$Z _{2,2}(t,y;x)=E _x\Big [ \exp \big \{-\frac{1}{2} \int  _
\mathbb R L _t^z Z^{(y)}(dz) \big \} \ 1_{\{T_y <t\}} \Big ].$$
\noi Theorem \ref{thbase1} tells us that $\sqrt{1+t}\ Z
_{2,1}(t,y;x)$ converges as $t\rightarrow\infty$.

\noi As for $Z _{2,2}(t,y;x)$, we use the strong Markov property
at time $T_y$:
$$\sqrt{1+t}\ Z _{2,2}(t,y;x)=E _x\Big [1_{\{T_y <t\}} \  \exp \big \{-\frac{1}{2}
\int  _ \mathbb R L _{T_y}^z V^{(y)}(dz) \big \}  \
\sqrt{\frac{1+t}{t+1-T_y}}$$
$$\qquad \qquad \times \sqrt{t+1-T_y} Z^{V^{(y)}}_{t-T_y}
(y) \Big ].$$
 \noi The conjonction of Theorem \ref{thbase1}, inequality (\ref{maj1}) (i) and
(\ref{familleuni}) implies that
$$ \Big ( 1_{\{T_y <t\}}
\sqrt{\frac{1+t}{t+1-T_y}}\sqrt{t+1-T_y} Z^{V^{(y)}}_{t-T_y} (y)
;t\geq 1\Big )$$
\noi  is a family of uniformly integrable r.v.'s converging a.s.
to  the constant $\varphi_{V^{(y)}}(y)$, as $t\rightarrow\infty$.
Hence, it converges in $L^1$.

\noi As a result, $\sqrt{1+t}\ Z _{2,2}(t,y;x)$ converges as
$t\rightarrow\infty$.

\noi This ends the proof of Lemma \ref{etape2unil}.
\end{prooff}


\vskip 15 pt \noi
\begin{prooff}\ {\bf of Theorem \ref{thunil}}.  a) Let $y>\max\{0,x\}$. Using
Lemma \ref{etape2unil} we get:
$$\limsup_{t\rightarrow\infty}\Big (\sqrt{1+t}\ Z^V_t(x)\Big )\leq
2\tilde{\varphi}(y) + \lim_{t\rightarrow\infty}\Big (\sqrt{1+t}\
Z_2 (t,y;x)\Big ),$$
$$\liminf_{t\rightarrow\infty}\Big (\sqrt{1+t}\ Z^V_t(x)\Big )\geq
\lim_{t\rightarrow\infty}\Big (\sqrt{1+t}\ Z_2 (t,y;x)\Big ).$$
\noi Hence
$$0\leq \limsup_{t\rightarrow\infty}\Big (\sqrt{1+t}\ Z^V_t(x)\Big )
-\liminf_{t\rightarrow\infty}\Big (\sqrt{1+t}\ Z^V_t(x)\Big ) \leq
2\tilde{\varphi}(y).$$

 \noi The parameter $y$ being arbitrary,
property (\ref{limCtiunil}) implies point 1. of Theorem
\ref{thunil} :
$$\varphi _V (x):=\lim_{t\rightarrow \infty} \big(\sqrt{t}Z^V _t(x)\big)
=\lim_{t\rightarrow \infty} \Big(\sqrt{t}E _x\Big [ \exp \big
\{-{\frac{1}{2}} \int  _ \mathbb R L _t^z V(dz) \big \} \Big
]\Big).$$
%
 b)(\ref{majC-unil}) (resp. (\ref{majC+unil})) is a direct
consequence of (\ref{etudeCtil}) (resp. the inequality :
$V(dz)\geq 1_{]-\infty ,0]}(z)V(dz)$ and (\ref{maj2C})).

\noi c) Obviously, (\ref{majC-unil}) implies $\displaystyle
\lim_{x\rightarrow +\infty} \varphi_V(x) =0.$

\noi In order to end the proof of Theorem \ref{thunil} we have to
check:
\begin{equation}
\lim_{x\rightarrow - \infty} \varphi_V'(x) =-\sqrt{\frac{2}{\pi}}.
\end{equation}
%

\noi Let $x<0$. We have successively:
$$\varphi_V'(0)-\varphi_V'(x)=\int_x^0 \varphi_V''(dy)=\int_x^0 \varphi_V(y)V(dy),$$
$$|\varphi_V'(0)-\varphi_V'(x)|\leq C \int_x^0 (1+|y|)V(dy).$$
\noi Assumption (\ref{prop2-u}) implies that $|\varphi_V'(x)|$ is
bounded. But $\varphi_V$ is a convex function, hence
$\displaystyle \lim_{x\rightarrow - \infty} \varphi_V'(x)$ exists.
Moreover
\begin{equation}\label{limconvunil}
\lim_{x\rightarrow - \infty} \varphi_V'(x)= \lim_{x\rightarrow -
\infty}\frac{\varphi_V(x)}{x}.
\end{equation}
Let $V^{[a]}(dy)$ denote the positive measure :
$V^{[a]}(dy)=1_{]-\infty,a[}V(dy)$, with $a<0$.

\noi It is clear that $V^{[a]}(dy)$ fulfills (\ref{prop2+u}) and
(\ref{prop2-u}). Let $\varphi^{[a]}(x)$ be the limit of $\sqrt{t}\
Z^{V ^{[a]}} _t(x)$, as $t\rightarrow \infty$, where :
$$Z^{V ^{[a]}} _t(x)=
\ E_x \Big [ \exp \big \{ -\frac{1}{2} \int_\mathbb R L^y_t
V^{[a]}(dy)\big\}\Big ],$$
\noi Let $x<a$. We have :
$$Z^V _t (x)=E _x\Big [ \exp \big \{-\frac{1}{2} \int  _ \mathbb R L _t^y V(dy)
\big \}  \ 1_{\{T_a < t\}} \Big] +E _x\Big [ \exp \big
\{-\frac{1}{2} \int  _ \mathbb R L _t^y V^{[a]}(dy) \big \}  \
1_{\{T_a \geq t\}} \Big] ,$$
\noi We use the strong Markov property at time $T_a$ :
$$Z^V _t (x)=E _x\Big [ \exp \big \{-\frac{1}{2} \int  _ \mathbb R L _{T_a}^y V(dy)
\big \}  \ 1_{\{T_a < t\}} Z^V_{t-T_a} (a)\Big] +E _x\Big [ \exp
\big \{-\frac{1}{2} \int  _ \mathbb R L _t^y V^{[a]}(dy) \big
\}\Big]
$$
$$-E _x\Big [ \exp \big \{-\frac{1}{2} \int  _ \mathbb R L _{T_a}^y V^{[a]}(dy)
\big \}  \ 1_{\{T_a < t\}} Z^{V^{[a]}} _{t-T_a} (a)\Big].$$
\noi We multiply both sides by $\sqrt{t}$ and we take the limit as
$t \rightarrow\infty$ to obtain :
$$ \varphi_V(x)=\varphi_V(a)h(x)+\varphi^{[a]}(x)+\varphi^{[a]}(a)h^{[a]}(x),
$$
\noi where
$$h(x)=E _x\Big [ \exp \big \{-\frac{1}{2} \int  _ \mathbb R L _{T_a}^y V(dy)
\big \}\Big], \quad h^{[a]}(x)= E _x\Big [ \exp \big
\{-\frac{1}{2} \int _ \mathbb R L _{T_a}^y V^{[a]}(dy) \big
\}\Big]. $$
 \noi The functions $h$ and $h^{[a]}$ are bounded, and $V^{[a]}(dy)$
 satisfies (\ref{propu}), hence
$$\lim_{x\rightarrow -\infty}\frac{\varphi_V(x)}{x}=
\lim_{x\rightarrow -\infty}\frac{\varphi^{[a]}(x)}{x}=
 \lim_{x\rightarrow - \infty} (\varphi^{[a]})'(x)=-\sqrt{\frac{2}{\pi}}.
 $$
\end{prooff}

\begin{exam}\label{ex1unil}
Let $V(dy)=\lambda^21_{[0,+\infty[}(y)dy$.
$$\displaystyle \lim _{t\rightarrow \infty}
\Big (\sqrt {t} \ E _x\Big [ \exp \big \{-{\frac{\lambda^2}{2}}
\int _ 0^{t }1_{\{B_s >0\}}ds \big \} \Big ]\Big )=\varphi_V(x),$$
 \noi where
$$
 \varphi_V(x)= \left\{
\begin{array}{ll}
\displaystyle \frac{1}{\lambda}\sqrt{\frac{2}{\pi}}e^{-\lambda x}
& {\rm \ if \ } x \geq 0 \\
 & \\
 \displaystyle \sqrt{\frac{2}{\pi}}(\frac{1}{\lambda}-x) &
{\rm \ if \ } x < 0 .
\end{array}
\right.
$$
\noi Moreover an explicit formula for $\displaystyle E _x\Big [
\exp \big \{-{\frac{\lambda^2}{2}} \int _ 0^{t }1_{\{B_s
>0\}}ds \big \} \Big ]$ is given in (\cite{BorSal} p 136]):
$$E _x\Big [ \exp \big \{-{\frac{\lambda^2}{2}}
\int _ 0^{t }1_{\{B_s >0\}}ds \big \} \Big ]$$
$$= \left\{
\begin{array}{ll}
\displaystyle e^{-\frac{\lambda ^2 t}{2}}(1-{\rm
Erfc}(-\frac{x}{\sqrt{2t}})+\frac{1}{\pi} \int_0^1
\frac{du}{\sqrt{u(1-u)}}\exp\{-\frac{\lambda^2 tu}{2}
-\frac{x^2}{2tu}\}
& {\rm  if \ } x \geq 0  \\
 & \\
\displaystyle (1-{\rm Erfc}(-\frac{x}{\sqrt{2t}})+\frac{1}{\pi}
\int_0^1 \frac{du}{\sqrt{u(1-u)}}\exp\{-\frac{\lambda^2 tu}{2}
-\frac{x^2}{2t(1-u)}\} & {\rm   if \ } x < 0 .
\end{array}
\right.
$$
\end{exam}

\begin{exam} \label{ex2unil}
Let $V(dy)=(y^2-1)1_{\{y\geq 1\}}dy$. Then
$$\displaystyle \lim _{t\rightarrow \infty}
\Big (\sqrt {t} \ E _x\Big [ \exp \big \{-{\frac{\lambda^2}{2}}
\int _ 0^t(B_s^2-1)1_{\{B_s >1\}}ds \big \} \Big ]\Big
)=\varphi_V(x),$$
\noi where
$$
 \varphi_V(x)= \left\{
\begin{array}{ll}
\displaystyle \sqrt{\frac{2}{\pi}}e^{1/2}e^{-x^2/2}

& {\rm \ if \ } x \geq 1 \\
 & \\
 \displaystyle \sqrt{\frac{2}{\pi}}(2-x) &
{\rm \ if \ } x < 1 .
\end{array}
\right.
$$

\end{exam}


\section {Some critical cases} \label {Bes}
\setcounter {equation}{0}

In this section we consider :
\begin{equation}\label{defuBes}
    V(x)=\frac{\lambda}{\theta +x^2} ,
\end{equation}
\noi where $\lambda>0$ and $\theta \geq 0$.

\noi Denoting by $(R_t;t\geq 0)$ the reflecting Brownian motion :
$R_t=|B_t|;\ t\geq 0$, then
 \begin{equation}\label{identBes}
    \int_0^t V(B_s) ds=\lambda \int_0^t \frac{1}{\theta
 +R_s^2}ds.
\end{equation}
\noi This led us to investigate more generally the asymptotic
behaviour of $\displaystyle E^\mu_x \Big [ \exp \Big \{
-\frac{\lambda}{2} \int _0 ^t\frac{ds}{\theta+R_s^2} \Big \} \Big
]$ when $t\rightarrow\infty$, where $(R_t;t\geq 0)$ is under
$P^\mu_x$, a Bessel process started at $x$, with  index $\mu >-1$
(the dimension is $d_\mu =2( \mu+1)$).

\noi Throughout this section, $n_\mu$ stands for :
\begin{equation}\label{nmuBes}
    n_\mu =\frac{-\mu+\sqrt{\mu^2+\lambda}}{2}
\end{equation}
This parameter will play a central role in the formulation of our
results.

 \noi We begin with the case $\theta =0$.


\begin {theorem} \label{theo1Bes}
Suppose $\mu > -1$. Then :
\begin{equation}\label{ident1Bes}
\lim _{t\rightarrow\infty} \Big ( t^{n_\mu}E^{\mu}_x\Big [ \exp
\Big \{ -\frac{\lambda}{2} \int _0 ^t\frac{ds}{R_s^2} \Big \} \Big
] \Big )=x^{2n_\mu} \frac{1}{2^{n_\mu}}\frac{\Gamma (\mu +n_\mu
+1)} {\Gamma (\mu +2n_\mu +1)}.
\end{equation}
\end{theorem}

\begin{rem} \label {rem1Bes}
\begin{enumerate} \item
 In particular if $\mu =-1/2$, that is if $d_\mu
=1$, then :
\begin{equation}\label{ident2Bes}
    \lim _{t\rightarrow\infty} \Big ( t^{n}E_x\Big [ \exp
\Big \{ -\frac{\lambda}{2} \int _0 ^t\frac{ds}{B_s^2} \Big \} \Big
] \Big )=x^{2n} \frac{1}{2^{n}}\frac{\Gamma (n +\frac{1}{2})}
{\Gamma (2n +\frac{1}{2})} ,
\end{equation}
\noi where $\displaystyle
n=n_{-1/2}=\frac{1+\sqrt{1+4\lambda}}{4}$.

\item Taking  $\mu =0$ (i.e. $d_\mu=2$), we obtain :
\begin{equation}\label{001Bes}
    \lim _{t\rightarrow\infty}\Big( t^\lambda E_x\Big [ \exp \Big \{
-\frac{\lambda^2}{2} \Big(4 \int _0 ^t\frac{ds}{R_s^2}\Big) \Big
\}\Big] \Big) =\frac{x^{2\lambda}}{2^\lambda}\frac{\Gamma (\lambda
+1)}{\Gamma (1+2\lambda )}=
\frac{x^{2\lambda}}{8^\lambda}\frac{\sqrt{\pi}}{\Gamma (1/2
+\lambda )} , \quad \lambda \geq 0.
\end{equation}
To obtain the last equality  in (\ref{001Bes}) we have used the
Legendre duplication formula (\cite{WhittWat}, p 240) :
$$
\Gamma(2z)=\frac{1}{\sqrt{\pi}}2^{2z-1}\Gamma (z)\Gamma
(z+\frac{1}{2}).
$$
The  formula (\ref{001Bes}) led Roynette and Yor \cite{RoyYor}  to
define and study a family of positive r.v's $(H_{c,\alpha} )$ such
that
$$
E\Big[\exp \Big \{ -\frac{\lambda^2}{2} H_{c,\alpha} \Big\}\Big]=
\frac{\Gamma (\alpha)}{\Gamma(\alpha+\lambda)}\exp\{c\lambda\} ,
\quad \lambda \geq 0 ,$$
where $\alpha >0$ and $c \leq \Gamma'(\alpha)/\Gamma(\alpha)$.

\noi We observe that these Laplace transforms appear in
(\ref{001Bes}) for $\alpha=1/2$, and $c=2\log x -\log 8$.

\end{enumerate}
\end {rem}


\begin {prooff} \ {\bf of Theorem \ref{theo1Bes}}
Our approach is based on the well-known identity \cite[chapter XI,
ex 1.22, page 430]{RevYor}, or  \cite{Yor}.
\begin{equation}\label{equiBes}
E^{\mu}_x\Big [Y \Big (\frac{x}{R_t} \Big )^\mu \exp \Big \{
-\frac{\nu^2}{2} \int _0 ^t\frac{ds}{R_s^2} \Big \} \Big ]=
E^{\nu}_x\Big [Y \Big (\frac{x}{R_t} \Big )^\nu \exp \Big \{
-\frac{\mu^2}{2} \int _0 ^t\frac{ds}{R_s^2} \Big \} \Big ],
\end{equation}
for any ${\cal F}_t$-measurable r.v. $Y \geq 0$.

\noi Choosing : $\nu =\mu +2n_\mu=\sqrt{\mu^2+\lambda}$ and
$$Y=\Big (\frac{x}{R_t} \Big )^{-\mu}
\exp \Big \{ \frac{\mu^2}{2} \int _0 ^t\frac{ds}{R_s^2} \Big \}
,$$
\noi we get :
$$E^{\mu}_x\Big [ \exp \Big \{
-\frac{\lambda}{2} \int _0 ^t\frac{ds}{R_s^2} \Big \} \Big ]=
\frac{x^{2n_\mu}}{t^{n_\mu}} E_{x/\sqrt{t}}^{\nu}\Big [ \Big (
\frac{1}{R_1^2} \Big )^{n_\mu}\Big ].
$$
\noi Then, we obtain :
$$\lim _{t\rightarrow\infty} \Big ( t^{n_\mu}E^{\mu}_x\Big [ \exp \Big \{
-\frac{\lambda}{2} \int _0 ^t\frac{ds}{R_s^2} \Big \} \Big ] \Big
)=x^{2n_\mu} E_{0}^{\nu}\Big [ \Big ( \frac{1}{R_1^2}\Big
)^{n_\mu} \Big ].$$
But, under $P^{\nu}_0$, the distribution of $R_1^2/2$ is
gamma$(\nu +1)$. A straightforward calculation yields to :
\begin{equation}\label{0Bes}
E_{0}^{\nu}\Big [ \Big ( \frac{1}{R_1^2}\Big )^{n_\mu} \Big ]=
\frac{1}{2^{n_\mu}} \frac{\Gamma (\nu -n_\mu +1)}{\Gamma (\nu
+1)}=\frac{1}{2^{n_\mu}} \frac{\Gamma (\mu +n_\mu +1)} {\Gamma
(\mu +2n_\mu +1)}.
\end{equation}
\end{prooff}

\ We now investigate the case $\theta>0$. In the sequel,
$L^{(\mu)}$ denotes the infinitesimal generator of the Bessel
process with index $\mu$ :
\begin{equation}\label{def1Bes}
L^{(\mu)}(f)(x)=\frac{1}{2}f"(x)+\frac{2\mu+1}{2x}f'(x).
\end{equation}
%

\begin{theorem} \label{theo2Bes}
Suppose $\mu \geq -1/2, \lambda>0$ such that :
\begin{equation}\label{002Bes}
    \lambda<8\mu^2+6\mu+1 .
\end{equation}
 Let $\varphi^{(\mu)}_\lambda$ be the
unique smooth function defined on $[0,+\infty[$, solution of
\begin{equation}\label{ident3Bes}
L^{(\mu)}(\varphi)(x)=\frac{1}{2}\varphi"(x)+\frac{2\mu+1}{2x}\varphi'(x)=
\frac{\lambda}{2}\frac{1}{1+x^2}\varphi(x); \quad  x>0,
\end{equation}
such that:
\begin{equation}\label{1Bes}
    \varphi^{(\mu)}_\lambda(x) \sim x^{2n_\mu} , x\rightarrow +\infty.
\end{equation}
\noi Then :
\begin{equation}\label{ident4Bes}
    \lim _{t\rightarrow\infty} \Big ( t^{n_\mu}E^{\mu}_x\Big [ \exp
\Big \{ -\frac{\lambda}{2} \int _0 ^t\frac{ds}{\theta+R_s^2} \Big
\} \Big ] \Big )=\theta^{n_\mu}\varphi^{(\mu)}_\lambda
(x/\sqrt{\theta}) \frac{1}{2^{n_\mu}}\frac{\Gamma (\mu +n_\mu +1)}
{\Gamma (\mu +2n_\mu +1)}.
\end{equation}
\end{theorem}

\begin{rem} We observe that if we take the limit $\theta\rightarrow 0$ in
(\ref{ident4Bes})we recover (\ref{ident1Bes}).
\end{rem}
The function $\varphi^{(\mu)}_\lambda$ is defined in terms of
hypergeometric functions. Let $F(\alpha ,\beta, \gamma ;x)$ be the
hypergeometric function with parameters $\alpha, \beta, \gamma$
(cf \cite{lebed}):
\begin{equation}\label{defhyperBes}
\displaystyle  F(\alpha ,\beta, \gamma ;x)=\sum _{k=0}^\infty
\frac{(\alpha)_k (\beta)_k}{(\gamma)_k k!}x^k,
\end{equation}
where $(\rho)_k=\rho (\rho +1)\times \cdots \times (\rho +k-1)$.

\noi The series in (\ref{defhyperBes}) converges for any $x$ such
that $|x|<1$.

\begin {lemma} \label {lem1Bes}

 The function $\varphi^{(\mu)}_\lambda$ in Theorem \ref{theo2Bes} which
 solves (\ref{ident3Bes}) and satisfies (\ref{1Bes}) is given by :
\begin{equation}\label{def2Bes}
    \varphi^{(\mu)}_\lambda (x)=
    \left \{
    \begin{array}{ll}
k_\mu F(n_\mu +\mu, -n_\mu,\mu+1;-x^2) & \ \mbox { if } \ n_\mu
\mbox {
is an integer }\\
k_\mu (1+x^2)^{-n_\mu-\mu} F(n_\mu
+\mu,\mu+1+n_\mu,\mu+1;\frac{x^2}{1+x^2}) & \ \mbox { otherwise, }
\end{array}
\right.
\end{equation}
\noi where
$$k_\mu=\frac{\Gamma (\mu +n_\mu )}{\Gamma (\mu +2n_\mu)}
\frac{\Gamma (\mu +n_\mu+1)}{\Gamma (\mu +1)}.
$$
\end {lemma}

\begin{rem}\label{rem2Bes}
We observe that $ n_\mu$ is a positive real number. If $n_\mu$ is
an integer, then $\varphi^{(\mu)}_\lambda (x)$ is a polynomial
function with degree $2n_\mu$.
\end{rem}


\begin {prooff} \ {\bf of Lemma \ref{lem1Bes}}.

  1) Recall that  $F(\alpha ,\beta, \gamma ;\cdot)$ fulfills :
\begin{equation}\label{equadifBes}
    x(1-x)u"(x)+(\gamma -(\alpha+\beta
    +1)x)u'(x)=\alpha \beta u(x).
\end{equation}
\noi Let $v$ be the function : $v(t)=F(\alpha ,\beta, \gamma
;t^2)$. Then \cite[p 164]{lebed}, $v$ is a solution to :
$$t(1-t^2)v"(t)+2(\gamma -\frac{1}{2} -(\alpha +\beta
    +\frac{1}{2})t^2)v'(t)=4 \alpha \beta v(t).
    $$
\noi Finally setting $w(t)=v(it)$, we obtain :
\begin{equation}\label{eq1Bes}
    \frac{1}{2}w"(t) +\frac{1}{t(1+t^2)}\Big [
\gamma -\frac{1}{2} -(\alpha +\beta
    +\frac{1}{2})t^2) \Big ] w'(t)=-2\alpha \beta
    \frac{w(t)}{1+t^2}.
\end{equation}
\noi If we choose $\alpha = n_\mu +\mu, \beta=-n_\mu$ and
$\gamma=\mu+1$ , then it is easy to check that $F(n_\mu +\mu,
-n_\mu,\mu+1;-x^2)$ solves (\ref{ident3Bes}), with $x \in ]0,1[$.

2) If $n_\mu$ is an integer, it is obvious that $F(n_\mu +\mu,
-n_\mu,\mu+1;-x^2)$ is a polynomial function with degree $2n_\mu$,
and then solves (\ref{ident3Bes}) for every $x>0$. Writing :
\begin{equation}\label{4Bes}
F(n_\mu +\mu,-n_\mu,\mu+1;-x^2)=\sum_{k=0}^{n_\mu}a_kx^{2k},
\end{equation}
\noi we obtain :
\begin{equation}\label{2Bes}
   a_k=\frac{(n_\mu +\mu)_k\times n_\mu\times \cdots(n_\mu-k+1)}{(\mu
+1)_k k!}.
\end{equation}
\noi Hence $a_k>0, \ F(n_\mu +\mu,-n_\mu,\mu+1;-x^2)>0$ and
\begin{equation}\label{3A0Bes}
    a_{n_\mu}=\frac{\Gamma (\mu +2n_\mu )}{\Gamma (n_\mu +\mu)}
\frac{\Gamma (\mu +1)}{\Gamma (\mu +n_\mu+1)}.
\end{equation}
 \noi Consequently $\varphi^{(\mu)}_\lambda (x)$ satisfies (\ref {1Bes}).

 3) Suppose that  $n_\mu$ is not an integer. To
obtain a function defined on the half-line $[0,+\infty[$ we use a
fractional linear transformation of hypergeometric functions.
Recall \cite[(9.5.1)]{lebed}  :
$$
F(\alpha ,\beta ,\gamma ;z)=(1-z)^{-\alpha}F(\alpha ,\gamma -\beta
,\gamma ;\frac{z}{z-1}), \quad |\mbox {arg}(1-z)|< \pi.
$$
\noi In our context this identity becomes :
\begin{equation}\label{hypBes}
    F(n_\mu +\mu,-n_\mu,\mu+1;-x^2)=(1+x^2)^{-n_\mu-\mu}
F(n_\mu +\mu,\mu+1+n_\mu,\mu+1;\frac{x^2}{1+x^2}), \quad x\in
\mathbb{R}.
\end{equation}
\noi An analytic continuation argument shows that
$\varphi^{(\mu)}_\lambda $ is a solution of (\ref{ident3Bes}) for
every $x>0$.

\noi It is easy to check that the coefficients in the series are
positive, thus, $\varphi^{(\mu)}_\lambda >0$.

\noi We conclude from \cite[page 297 , ex 8]{WhittWat} :
$$F(\alpha ,\beta ,\gamma ;\frac{x^2}{1+x^2})\sim
\frac{\Gamma (\alpha +\beta-\gamma)\Gamma (\gamma)}{\Gamma
(\alpha)\Gamma (\beta)}\ \frac{1}{(1+x^2)^{\gamma -\alpha
-\beta}}, \quad x\rightarrow +\infty ,
$$%
that (\ref {1Bes}) holds.

\end {prooff}

\begin {lemma} \label {lem2Bes}
\begin{enumerate}
\item
 The function $\varphi^{(\mu)}_\lambda$ defined in Lemma \ref{lem1Bes}
fulfills :
\begin{equation}\label{3Bes}
    (\varphi^{(\mu)}_\lambda )' (x)\geq 0, \quad x \geq 0,
\end{equation}
\begin{equation}\label{5Bes}
    2n_\mu -\frac{\rho_\mu}{1+x^2}\leq \frac{x (\varphi^{(\mu)}_\lambda)' (x)}
    {\varphi^{(\mu)}_\lambda (x)}
    \leq 2n_\mu, \quad x \geq 0,
\end{equation}
\noi where $\rho_\mu>0$.
\item Let $D^{(\mu)}_\lambda $ be the function defined on
$[1,+\infty[\times [0,+\infty[$ by :
\begin{equation}\label{5ABes}
D^{(\mu)}_\lambda (t,x)=\frac{1}{t^{n_\mu}}\varphi^{(\mu)}_\lambda
(x\sqrt{t}).
\end{equation}
\noi Then :
\begin{equation}\label{5BBes}
    \lim_{t\rightarrow +\infty}D^{(\mu)}_\lambda (t,x)=x^{2n_\mu} ;
    \quad  D^{(\mu)}_\lambda (t,x) \geq \hat{\rho}_\mu x^{2n_\mu} ,
    \forall t,x \geq0, \
    \mbox{ for some }\ \hat{\rho}_\mu >0.
\end{equation}
\begin{equation}\label{5CBes}
0\leq \frac{1}{D^{(\mu)}_\lambda(t,z)}-
\frac{1}{D^{(\mu)}_\lambda(t,y)}\leq k \frac{y-z}{z^{2n_\mu+1}} ;
0 <z<y.
\end{equation}
where $\displaystyle k=\frac{2n_\mu}{\hat{\rho}_\mu}$.
\end{enumerate}
\end {lemma}

\begin{prooff}\ {\bf of lemma \ref{lem2Bes}}
\noi We give the proof only for the case $n_\mu \in\mathbb{N}$,
the other cases are left to the reader.

\noi Property (\ref{2Bes}) gives (\ref{3Bes}).

\noi Using (\ref{4Bes}) we get :
$$2n_\mu \varphi^{(\mu)}_\lambda (x) -x(\varphi^{(\mu)}_\lambda )' (x)=2k_\mu \Big (
\sum _{k=0}^{n_\mu}(n_\mu -k)x^{2k}\Big )\geq 0.$$
\noi Hence
$$2n_\mu  -x\frac{(\varphi^{(\mu)}_\lambda )' (x)}{\varphi^{(\mu)}_\lambda (x)}
=\frac{Q(x)}{\varphi^{(\mu)}_\lambda (x)},
$$
\noi where $Q$ is a polynomial function with degree less than
$2(n_\mu -1)$. This implies (\ref{5Bes}).

\noi Property (\ref{5BBes}) is due to the fact that
$D^{(\mu)}_\lambda (t,\dot)$ is a polynomial function with degree
$2n_\mu$ and positive coefficients.

\noi As for (\ref{5CBes}), we take the $x$-derivative of
$1/D^{(\mu)}_\lambda(t,x)$, we obtain :
$$
\Big|\frac{\partial}{\partial x}\frac{1}{D^{(\mu)}
_\lambda(t,x)}\Big|= t^{n_\mu+1/2}
\frac{(\varphi^{(\mu)}_\lambda)'(x\sqrt{t})}{(\varphi^{(\mu)}_\lambda)^2(x\sqrt{t})}.$$
Hence, the inequalities (\ref{5Bes}) and (\ref{5BBes}) directly
imply :
$$
\Big|\frac{\partial}{\partial x}\frac{1}{D^{(\mu)}
_\lambda(t,x)}\Big| \leq
\frac{2n_\mu}{\hat{\rho}_\mu}\frac{1}{x^{1+2n_\mu}}.
$$
Then (\ref{5CBes}) follows immediately.
\end{prooff}

\begin {lemma} \label {lem3Bes}
For any positive functional $F$ and $x\geq 0,\ t>0$, we have :
\begin{equation}\label{6Bes}
E_x^\mu \Big[ F(R_s;0 \leq s \leq t)\varphi^{(\mu )}_\lambda (R_t)
\exp\Big\{ -\frac{\lambda}{2}\int _0^t\frac{ds}{1+R_s^2}\Big\}
\Big ] =\varphi^{(\mu)}_\lambda (x)E_x\Big[ F(X_s;0\leq s \leq
t)\Big ] ,
\end{equation}
\noi where the function $\varphi^{(\mu)}_\lambda$ is defined in
Lemma \ref{lem1Bes} and $(X(t);t\geq 0)$ is the process solution
of :
\begin{equation}\label{solXBes}
    X_t=x+B_t+\frac{2\mu+1}{2}\int _0^t\frac{ds}{X_s}+\int _0^t
    \frac{(\varphi^{(\mu )}_\lambda)'}{\varphi^{(\mu )}_\lambda}(X_s)ds, \quad t \geq 0.
\end{equation}
 \noi In particular :
\begin{equation}\label{7Bes}
E_x^\mu \Big[\exp\Big\{ -\frac{\lambda}{2}\int
_0^t\frac{ds}{1+R_s^2}\Big\} \Big ]=\varphi^{(\mu )}_\lambda
(x)E_x\Big[\frac{1}{\varphi^{(\mu )}_\lambda (X_t)}\Big].
\end{equation}
\end {lemma}

\begin{rem} \label{RBesHY}
Hariya and Yor  \cite{HarYor} show the existence, and describe,
the limiting measures, as $t\rightarrow \infty $, of the laws of
$\{ B_s+\mu s; 0 \leq s \leq t\}$ perturbed by the Radon-Nikodym
density consisting of either the normalized functionals $\exp
(-\alpha A_t ^{(\mu)})$, or $1/\Big(A_t ^{(\mu)}\Big)^m$, where
$\displaystyle A_t ^{(\mu)}=\int _0^t ds \exp\{2(B_s+\mu s)\}$.
The results exhibit different regimes according to whether $\mu
\geq 0$ or $\mu<0$ in the first case, and a partition of the $(\mu
, m)$-plane in the second case.
\end{rem}

\begin{prooff}\ {\bf of lemma \ref{lem3Bes}}.
Let $(R_t; t \geq 0)$ and $(X_t; t \geq 0)$ be defined as
solutions of :
\begin{equation}\label{8Bes}
 X_t=x+B_t+\frac{2\mu+1}{2}\int _0^t\frac{ds}{X_s}+\int _0^t
    \frac{(\varphi^{(\mu )}_\lambda)'}{\varphi^{(\mu )}
    _\lambda}(X_s)ds, \quad t \geq 0,
\end{equation}
\begin{equation}\label{9Bes}
    R_t=x+B_t+\frac{2\mu+1}{2}\int _0^t\frac{ds}{R_s},\quad t \geq 0,
\end{equation}
\noi with the same underlying Brownian motion $(B_t; t \geq 0)$.

\noi Since $(\varphi^{(\mu )}_\lambda)'/\varphi^{(\mu )}_\lambda$
is a bounded function, we may apply Girsanov's theorem :
$$E_x\Big[ F(X_s;0\leq s \leq t)\Big ]=E_x^\mu \Big[ Y_t F(R_s;0 \leq s \leq t) \Big
],
$$
\noi where
\begin{equation}\label{10Bes}
    Y_t=\exp \Big\{\int _0^t
\frac{(\varphi^{(\mu )}_\lambda)'}{\varphi^{(\mu
)}_\lambda}(R_s)dB_s-\frac{1}{2}\int _0^t
\Big(\frac{(\varphi^{(\mu )}_\lambda)'}{\varphi^{(\mu
)}_\lambda}\Big)^2(R_s)ds\Big\}.
\end{equation}
%
%
%
\noi Applying by now standard arguments (cf, formula
(\ref{4Aintro})), we obtain :
$$
Y_t=\frac{\varphi^{(\mu )}_\lambda(R_t)}{\varphi^{(\mu
)}_\lambda(x)} \exp \Big\{- \int _0^t
\frac{L^{(\mu)}(\varphi^{(\mu )}_\lambda )}{\varphi^{(\mu
)}_\lambda}(R_s) ds \Big\}.
$$
\noi We conclude from (\ref{ident3Bes}) that (\ref{6Bes}) holds.
\end{prooff}

\begin {lemma} \label {lem4Bes}
Let $x \geq 0$. Let  us denote by  $(X_t; t \geq 0)$ the diffusion
:
\begin{equation}\label{12Bes}
 X_t=x+B_t+\frac{2\mu+1}{2}\int _0^t\frac{ds}{X_s}+\int _0^t
    \frac{(\varphi^{(\mu )}_\lambda)'}{\varphi^{(\mu )}_\lambda}(X_s)ds, \quad t \geq 0,
\end{equation}
\noi and $(R_t^\mu ; t \geq 0)$ (resp.  $(R^{\mu +2n_\mu}_t; t
\geq 0)$ the Bessel processes with index $\nu=\mu$, (resp. $\nu
=\mu +2n_\mu$) solving  :
\begin{equation}\label{13Bes}
    R^\nu_t=x+\tilde{B}_t+\frac{2\nu+1}{2}\int _0^t\frac{ds}{R^\nu_s},\quad t \geq 0,
\end{equation}
\noi where $(\tilde{B}_t ; t \geq0)$ is the Brownian motion :
$\displaystyle \tilde{B}_t = \int_0^t \mbox{sgn}(X_s)dB_s, \ t\geq
0$.

\noi Then a.s. for any $t\geq 0$ :
\begin{equation}\label{ineg3Bes}
(R^\mu_t)^2 \leq X_t^2\leq (R^{\mu +2n_\mu}_t)^2.
\end{equation}

\end {lemma}

\begin{prooff}\ {\bf of lemma \ref{lem4Bes}}.
Applying It\^o's formula to the squares of the  processes $X$ and
$R^\nu$, we obtain :
\begin{equation}\label{15Bes}
X_t^2=x^2+2\int _0^t \sqrt{X_s^2}d\tilde{B}_s+ \int _0^t \Big
(2X_s \frac{(\varphi^{(\mu )}_\lambda)'}{\varphi^{(\mu
)}_\lambda}(X_s)+2(\mu +1)\Big) ds ,
\end{equation}
\begin{equation}\label{16Bes}
    (R^\nu_t)^2=x^2+2\int _0^t \sqrt{(R^\nu_s)^2}d\tilde{B}_s +2(\nu
+1)t.
\end{equation}
\noi We observe that the function $\displaystyle x\mapsto
\frac{2x(\varphi^{(\mu )}_\lambda)'(x)}{\varphi^{(\mu )}_\lambda
(x)}+2(\mu +1)$ may be written as $h(x^2)$, a function of $x^2$,
hence:
\begin{equation}\label{17ABes}
    X_t^2=x^2+2\int _0^t \sqrt{X_s^2}d\tilde{B}_s+ \int _0^t
h(X_s^2)ds.
\end{equation}

 \noi Inequalities (\ref{3Bes}) and
(\ref{5Bes}) imply :
%
%
\begin{equation}\label{17BBes}
    2(\mu +1) \leq h(x) \leq2(2n_\mu +\mu+1).
\end{equation}

\noi The inequalities (\ref{ineg3Bes}) are a direct consequence of
comparison results for solutions of one-dimensional stochastic
differential equations \cite[Prop 2.18, Chapter 5]{KaratShrev}.
\end{prooff}

\begin{rem} \label{Rem4bBes} In (\ref{13Bes}) we have chosen the
Brownian motion $(\tilde{B})$ instead of $(B)$ in order to obtain
 the inequalities (\ref{ineg3Bes})). Replacing $(\tilde{B_t})$ by $(B_t)$
 in (\ref{13Bes}) does not change the law of $(R^\nu_t)$.
\end{rem}
\begin {lemma} \label {lem5Bes}
Let $a>0$ be fixed. Let  $(Y_t^a; t \geq 0)$ and $(Z_t^a; t \geq
0)$ denote the processes :
\begin{equation}\label{17Bes}
 Y_t^a =\frac{1}{a}(R^\nu_{at})^2;  \quad Z_t^a=\frac{1}{a}X^2_{at}, \quad  t \geq 0,
\end{equation}
 with $\nu=\mu +2n_\mu$, $(X_t; t \geq 0)$ (resp. $(R^\nu_t; t \geq 0)$) solving
 (\ref{12Bes}) (resp. (\ref{13Bes})).

 \noi Then for any $1\leq p <\infty$, we have :
 \begin{equation}\label{18Bes}
  \lim_{a\rightarrow+\infty}\Big(  E[(Y^a_1 - Z^a_1
  )^p]\Big)=0.
\end{equation}
\end {lemma}

\begin{prooff}\ {\bf of Lemma \ref{lem5Bes}}.
From (\ref{15Bes}) and (\ref{16Bes}), we deduce :
$$
0 \leq E[Y^a_t - Z^a_t] \leq E \Big[ \int _0^{t} \Big (4n_\mu-
2X_{as} \frac{(\varphi^{(\mu )}_\lambda)'}{\varphi^{(\mu
)}_\lambda}(X_{as})\Big) ds \Big].
$$
\noi Using successively (\ref{5Bes}) and (\ref{ineg3Bes}), we
obtain :
$$E[Y^a_t - Z^a_t] \leq 2\theta_\mu E \Big[ \int _0^{t}\frac{1}{1+(X_{as})^2}ds \Big]
\leq 2\theta_\mu \int _0^{t}E \Big[\frac{1}{1+(R^\nu_{as})^2}
\Big]ds ,
$$
\noi The scaling property of Bessel processes yields to :
$$
E \Big[\frac{1}{1+(R^\nu_{as})^2}
\Big]=E^\nu_x\Big[\frac{1}{1+(R_{as})^2}
\Big]=E^\nu_{x/\sqrt{a}}\Big[\frac{1}{1+a(R_{s})^2} \Big].
$$
\noi Obviously the right hand-side of the previous inequality
tends to $0$, as $a\rightarrow +\infty$. The dominated convergence
theorem implies that :
$$
\lim_{a\rightarrow +\infty}\Big(E[Y^a_t - Z^a_t]\Big)=0.
$$
Let $1\leq p<\infty$ and $t=1$. Since $0\leq Y^a_1-Z^a_1$, then
$(Y^a_1-Z^a_1)^p$ goes to $0$ in probability, as
$a\rightarrow\infty$.

\noi We claim that for any $\alpha \in ]1,\infty[$ :
\begin{equation}\label{180ABes}
    \sup_{a \geq 1}E[(Y^a_1-Z^a_1)^{p\alpha}]<\infty.
\end{equation}

 \noi This will prove  (\ref{18Bes}),
because (\ref{180ABes}) implies that $((Y^a_1-Z^a_1)^p; a \geq 1)$
is uniformly integrable.

\noi To prove (\ref{180ABes}), we use the definition of $Y^a$ and
(\ref{ineg3Bes}). Denoting  $\beta=p\alpha$, we have :
$$E[(Y^a_1-Z^a_1)^{\beta}]\leq E[(Y^a_1)^{\beta}]
\leq E\Big[\Big(\frac{R^\nu_a}{\sqrt{a}}\Big)^{\beta}\Big].
$$
\noi Recall that with our notations :
$$E\Big[\Big(\frac{R^\nu_a}{\sqrt{a}}\Big)^{\beta}\Big]=
E^\nu_x\Big[\Big(\frac{R_a}{\sqrt{a}}\Big)^{\beta}\Big]=
E^\nu_{x/\sqrt{a}}[(R_1)^{\beta}].$$
The second equality follows from  the scaling property of Bessel
processes. Comparison theorem tells us :
$$E^\nu_{x/\sqrt{a}}\Big[(R_1)^{\beta}\Big]\leq
E^\nu_{x}\Big[(R_1)^{\beta}\Big]<\infty ,$$
for any $a\geq 1$. This proves (\ref{180ABes}).

%

  %
%
%
%
%
\end{prooff}

\vskip 10 pt
\begin{prooff} \ {\bf of Theorem \ref{theo2Bes}}
1) Using the scaling property of Bessel processes, we get :
\begin{equation}\label{180Bes}
 E^{\mu}_x\Big [ \exp
\Big \{ -\frac{\lambda}{2} \int _0 ^t\frac{ds}{\theta+R_s^2} \Big
\} \Big ] =E^{\mu}_{x/\sqrt{\theta}}\Big [ \exp \Big \{
-\frac{\lambda}{2} \int _0 ^{t/\theta}\frac{ds}{1+R_s^2} \Big \}
\Big ].
\end{equation}
\noi Thus, it suffices to prove (\ref{ident4Bes}) when $\theta
=1$.

2) Let $\theta =1$. From (\ref{7Bes}), it remains to prove :
\begin{equation}\label{19Bes}
    \lim _{t\rightarrow\infty} \Big ( t^{n_\mu}
    E_x\Big[\frac{1}{\varphi^{(\mu)}_\lambda (X_t)}\Big] \Big)
    =\frac{1}{2^{n_\mu}}\frac{\Gamma
(\mu +n_\mu +1)} {\Gamma (\mu +2n_\mu +1)}.
\end{equation}
\noi By (\ref{17Bes}) and (\ref{5ABes}), we get :
$$E_x\Big[\frac{1}{\varphi^{(\mu)}_\lambda (X_t)}\Big] =\frac{1}{t^{n_\mu}}
E_x\ \Big[ \frac{1}{D^{(\mu)}_\lambda(t,\sqrt{Z_1^t})} \Big] ,
$$
\noi 3) Let us prove :
\begin{equation}\label{190Bes}
\lim _{t\rightarrow\infty}E_x
\Big[\frac{1}{D^{(\mu)}_\lambda(t,\sqrt{Z_1^t})}-
\frac{1}{D^{(\mu)}_\lambda(t,\sqrt{Y^t _1})}\Big ]=0.
\end{equation}
\noi Using (\ref{5CBes}), we obtain :
$$E_x\Big[\frac{1}{D^{(\mu)}_\lambda(t,\sqrt{Z_1^t})}-
\frac{1}{D^{(\mu)}_\lambda(t,\sqrt{Y^t _1})}\Big ] \leq C
E_x\Big[\frac{\sqrt{Y_1^t}-\sqrt{Z_1^t}}{(Z_1^t)^{n_\mu +1/2}}
\Big]\leq C E_x\Big[\frac{\sqrt{Y_1^t -Z_1^t}}{(Z_1^t)^{n_\mu
+1/2}} \Big].
$$
\noi Since $\lambda <8\mu^2+6\mu+1$ and
$n_\mu=\frac{-\mu+\sqrt{\mu^2+\lambda}}{2}$, then
$n_\mu+1/2-\mu<1$. Hence we may find $\varepsilon>0$ such that :
\begin{equation}\label{191Bes}
n_\mu+1/2+\varepsilon-\mu<1 .
\end{equation}
Let $\displaystyle
q=q(\varepsilon)=\frac{n_\mu+1/2+\varepsilon}{n_\mu+1/2}>1$ and
$p$ be the conjugate exponent of $q$.

\noi Applying H\"{o}lder's inequality leads to :
$$E_x\Big[\frac{1}{D^{(\mu)}_\lambda(t,\sqrt{Z_1^t})}-
\frac{1}{D^{(\mu)}_\lambda(t,\sqrt{Y^t _1})}\Big ] \leq C
\Big\{E_x\big[(Y_1^t-Z_1^t)^{p/2}\big]\Big\}^{1/p}
\Big\{E_x\big[\frac{1}{(Z_1^t)^{n_\mu +1/2+\varepsilon}}\big]
\Big\}^{1/q}.
$$
Property (\ref{190Bes}) will be a direct consequence of
(\ref{18Bes}), once we have proved that : $\displaystyle
t\rightarrow E_x\big[\frac{1}{(Z_1^t)^{n_\mu
+1/2+\varepsilon}}\big]$ is bounded.

\noi Using the definition of $Z_1^t$,  (\ref{ineg3Bes}) and the
scaling property of Bessel processes, we obtain :
$$E_x\big[\frac{1}{(Z_1^t)^{n_\mu +1/2+\varepsilon}}\big]
=E_x^\mu \big[\Big(\frac{t}{R_t^2}\Big)^{n_\mu
+1/2+\varepsilon}\big]=E_{x/\sqrt{t}}^\mu
\big[\frac{1}{(R_1^2)^{n_\mu +1/2+\varepsilon}}\big].
$$
Comparison theorem implies that :
$$E_{x/\sqrt{t}}^\mu
\big[\frac{1}{(R_1^2)^{n_\mu +1/2+\varepsilon}}\big] \leq E_0^\mu
\big[\frac{1}{(R_1^2)^{n_\mu +1/2+\varepsilon}}\big].$$
\noi Under $P_0^\mu$, the distribution of $R_1^2/2$ is gamma$ (\mu
+1)$. Hence $\displaystyle E_0^\mu \big[\frac{1}{(R_1^2)^{n_\mu
+1/2+\varepsilon}}\big] <\infty$ as soon as :
$$
\int_0^1\frac{y^\mu}{y^{n_\mu +1/2+\varepsilon}}dy <\infty .$$
This integral is finite since condition (\ref{191Bes}) holds.

%
%
%

\noi 4) Due to the scaling property of Bessel processes,
$$
E_x\Big [\frac{1}{D^{(\mu)}_\lambda(t,\sqrt{Y^t _1})}\Big ]= E^
\nu _{x/\sqrt{t}}\Big [\frac{1}{D^{(\mu)}_\lambda (t,R_1)}\Big ],
$$
\noi where $\nu=\mu+2n_\mu$.

 \noi Applying (\ref{5BBes}), we
obtain :
$$
\lim _{t\rightarrow\infty} E^ \nu _{x/\sqrt{t}}\Big
[\frac{1}{D^{(\mu)}_\lambda(t, R _1)}\Big ]= E^ \nu _0 \Big
[\frac{1}{R^{2n_\mu} _1}\Big ].
$$
\noi Relation (\ref{19Bes}) now follows from (\ref{0Bes}).
\end{prooff}

\begin{rem} Note that condition (\ref{002Bes}) has only  been  used
in the last part of the proof of Theorem \ref{theo2Bes}. It may
not be  necessary but we have not been able to justify this.
\end{rem}

\begin{theorem}\label{theo3Bes}
Assume that $\lambda,\ \theta >0$ obey (\ref{002Bes}).  Let $Q
_{x,t}$ be the probability defined on ${\cal F} _t$ via :
\begin{equation} \label{def4Qbil}
Q _{x,t}(\Lambda  _t) = \frac{E _x \Big [1 _{\Lambda  _t}\ \exp
\Big \{ -\frac{\lambda}{2} \int _0 ^t\frac{ds}{\theta+R_s^2}\Big\}
\Big ]}
 {E _x \Big [ \exp
\Big \{ -\frac{\lambda}{2} \int _0 ^t\frac{ds}{\theta+R_s^2}
\Big\}\Big ]}
  ,\  \Lambda  _t \in {\cal F} _t.
\end{equation}
\noi Then, for any $\Lambda  _s$ in ${\cal F} _s$, $Q
_{x,t}(\Lambda  _s)$ converges to $ P^{\varphi^{(\mu)}_\lambda}
_{x}(\Lambda _s)$
 as $t \rightarrow \infty$, where $P^{\varphi^{(\mu)}_\lambda}_ x$ is the probability
 defined on ${\cal F} _\infty$ by :
\begin{equation} \label{defQbilA0}
P^{\varphi^{(\mu)}_\lambda} _x(\Lambda  _s
)=\frac{1}{\varphi^{(\mu)}_\lambda(x/\sqrt{\theta})} E _x\Big [1
_{\Lambda _s} \varphi^{(\mu)}_\lambda(R_ s/\sqrt{\theta}) \exp
\Big \{-{\frac{1}{2}} \int _0 ^s\frac{dv}{\theta+R_v^2} \Big
\}\Big ],
\end{equation}
\noi for any $s>0$ and $\Lambda  _s \in {\cal F} _s$.
 \item
 Let  $(X_t ^x;t\geq 0)$ be the solution to :
\begin{equation} \label{def3XBes}
X_t=x+B_t+\frac{2\mu+1}{2}\int
_0^t\frac{ds}{X_s}+\frac{1}{\sqrt{\theta}}\int _0^t
    \frac{(\varphi^{(\mu)}_\lambda)'}{\varphi^{(\mu)}_\lambda}(X_s/\sqrt{\theta})ds, \quad t \geq
    0.
\end{equation}
 Then the law of $(X_t ^x;t\geq 0)$ is $P^{\varphi^{(\mu)}_\lambda} _ x$.

\end{theorem}

 \noi The proof of Theorem \ref{theo3Bes} is based on Theorem
\ref{theo2Bes} and the estimate (\ref{5BBes}), the details are
left to the reader.

\vskip 10 pt \noi Let us mention two consequences of Theorem
\ref{theo3Bes}.

\begin{coro}Let $\mu \geq -1/2, \theta >0$ satisfying (\ref{002Bes}). Then
\begin{equation}
    \lim _{\varepsilon\rightarrow 0} \Big ( \frac{1}{\varepsilon^{n_\mu}}
    E^{\mu}_{x\sqrt{\varepsilon}} \Big [ \exp
\Big \{ -\frac{\lambda}{2} \int _0
^1\frac{ds}{\theta\varepsilon+R_s^2} \Big \} \Big ] \Big
)=\theta^{n_\mu}\varphi^{(\mu)}_\lambda (x/\sqrt{\theta})
\frac{1}{2^{n_\mu}}\frac{\Gamma (\mu +n_\mu +1)} {\Gamma (\mu
+2n_\mu +1)}.
\end{equation}
\end{coro}

\begin{coro} Assume $\mu \geq -1/2, \theta >0$ and
 (\ref{002Bes}) holds. Let $\psi_\lambda$ be the unique solution
of :
\begin{equation}\label{20Bes}
\left \{
    \begin{array}{ll}
\frac{\partial\psi}{\partial t}(t,x)=\frac{1}{2}\frac{\partial^2
\psi}{\partial x^2}(t,x)+\frac{2\mu+1}{2x} \frac{\partial
\psi}{\partial x}(t,x)-\frac{\lambda}{2}\frac{\psi
(t,x)}{\theta+x^2},\quad t>0,x \geq 0,
\\
\psi (0,x)=1.
\end{array}
\right.
\end{equation}
\noi Then
\begin{equation}\label{21Bes}
\lim _{t\rightarrow \infty} \Big (t^{n_\nu}\psi_\lambda (t,x)\Big
) =\theta^{n_\mu}\varphi^{(\mu)}_\lambda (x/\sqrt{\theta})
\frac{1}{2^{n_\mu}}\frac{\Gamma (\mu +n_\mu +1)} {\Gamma (\mu
+2n_\mu +1)}.
\end{equation}
\end{coro}


\section {On the use of large deviations}\label {LD}
\setcounter {equation}{0}

In this section we will be concerned with $\lambda V $, where
$\lambda  >0$ and :
\begin{equation}\label{1GD}
    V(x)=\frac{1}{1+|x|^\alpha} ;\qquad  x \in \mathbb{R},\
    0<\alpha<2.
\end{equation}

\noi We investigate the asymptotic behaviour of :
\begin{equation}\label{2GD}
    Z^{\lambda V}_t (x)=E_x\Big [ \exp \Big \{
    -\frac{\lambda}{2}\int_0^t\frac{ds}{1+|B_s|^\alpha} \Big \}
\Big ] ,
\end{equation}
\noi when $t\rightarrow\infty$.

 \noi Notice that if $\alpha>2$ then $\displaystyle \lambda \int
_\mathbb{R}V(x)|x|dx<\infty$, hence we may apply the results of
section 3. The critical case $\alpha =2$ has been treated in the
previous  section.


\begin{theorem}\label{theo1GD}
Let $0<\alpha<2,\ \lambda >0$. Then
\begin{equation}\label{3GD}
    \lim_{t\rightarrow\infty} \Big( t^{\frac{\alpha-2}{\alpha+2}}
\ln \big (Z^{\lambda V}_t (x)\big )
\Big)=-\frac{1}{2}I_0(\lambda),
\end{equation}
\noi where
\begin{equation}\label{4GD}
I_0(\lambda)=\inf _{\psi \in {\cal C}_0}\Big\{ \int
_0^1\dot{\psi}^2(s)ds +\lambda \int _0^1\frac{ds}{|\psi
(s)|^\alpha} \Big \},
\end{equation}
\noi belongs to $]0,+\infty[$, and ${\cal C}_0$ is the set of
continuous functions $f :[0,1] \rightarrow \mathbb{R}$ vanishing
at $0$.
\end{theorem}


\begin{rem} We observe that the limit in (\ref{3GD}) does not depend on
$x$, which may be   due to the fact this result  only gives a
logarithmic equivalent to $Z^V_t(x)$. Indeed, consider  the
equivalent of $Z^V_t(x)$ given by Theorem \ref{thbase1} :
$\displaystyle \lim
_{t\rightarrow\infty}\Big(\sqrt{t}Z^V_t(x)\Big)=\varphi_V(x)$.
Then $\displaystyle \frac{\ln(Z^V_t(x))}{\ln t}$ converges to the
constant $-1/2$, as $t\rightarrow\infty$.
\end{rem}

\begin{lemma}\label{lem1GD}
Let $\eta >0$. Let us denote :
\begin{equation}\label{5GD}
    I_\eta(\lambda)=\inf _{\psi \in {\cal C}_0}\Big\{ \int
_0^1\dot{\psi}^2(s)ds +\lambda \int _0^1\frac{ds}{\eta +|\psi
(s)|^\alpha} \Big \}.
\end{equation}
\noi Then $I_\eta(\lambda)$ is a positive real number,
$\eta\mapsto I_\eta(\lambda)$ is decreasing and
\begin{equation}\label{6GD}
\lim_{\eta \rightarrow0} I_\eta(\lambda) = I_0(\lambda).
\end{equation}
\end{lemma}

\begin{prooff} \ {\bf of Lemma \ref{lem1GD}}.
Let $\psi _\eta$ be a function in ${\cal C}_0$ such that :
\begin{equation}\label{60GD}
    I_\eta(\lambda)=\int _0^1\dot{\psi}_\eta^2(s)ds +\lambda \int
_0^1\frac{ds}{\eta +|\psi_\eta (s)|^\alpha}.
\end{equation}
\noi Then $\psi_\eta\geq 0$ and the  Euler equation associated
with $\psi _\eta$ is :
\begin{equation}\label{7GD}
    \left \{
    \begin{array}{l}
    2\ddot{\psi_\eta}+\frac{\alpha\lambda \psi _\eta ^{\alpha-1}}{(\eta+\psi _\eta
    ^\alpha)^2}=0
\\
\psi _\eta (0)=0, \qquad \dot{\psi _\eta}(1)=0 .
\end{array}
\right.
\end{equation}
\noi Consequently  $\psi _\eta$ is a positive, increasing and
convex function. Multiplying the first line of (\ref{7GD}) by
$\dot{\psi _\eta }$ and integrating, we obtain :
\begin{equation}\label{8GD}
    \dot{\psi _\eta}^2(t)=\lambda \Big (
    \frac{1}{\eta+\psi _\eta
    (t)^\alpha}-\frac{1}{\eta+\psi _\eta
    (1)^\alpha}\Big).
\end{equation}
 \noi Let $H_\eta$ be the function :
$$
H_\eta (C,x)=\frac{1}{\sqrt{\lambda}}\int _0
^x\frac{dy}{\sqrt{\frac{1}{\eta+y^\alpha}-\frac{1}{\eta+
   C ^\alpha}}}, \quad x\in [0,C], \ C>0.
$$
\noi Since the derivative of $\psi _\eta$ is  positive, the
relation (\ref{8GD}) is equivalent to :
$$
\frac{d}{dt}H_\eta (C,\psi _\eta (t)))=1; \quad 0\leq t\leq 1,
$$
\noi  or
$$
H_\eta (C,\psi_\eta (t))=t, \quad 0\leq t\leq 1,
$$
\noi with $C=\psi_\eta (1)$.

 \noi This implies that $\psi_\eta$ is the inverse of $t \ (\geq
0) \mapsto H_\eta (C,t)$. As for $C$, we observe that  it remains
to take into account the condition : $C= \psi _\eta (1)$. Let
$C_\eta$ be the unique solution in $]0,+\infty[$ of :
\begin{equation}\label{9GD}
\frac{1}{\sqrt{\lambda}}\int _0
^{C_\eta}\frac{dy}{\sqrt{\frac{1}{\eta+y^\alpha}-\frac{1}{\eta+
   C _\eta^\alpha}}}=1.
\end{equation}
\noi Taking  $C=C_\eta$, we have $C= \psi _\eta (1)$.

\end{prooff}

\begin{lemma}\label{lem2GD}
Let $\alpha \in $]0,2[ and $\psi_0$ be defined by the relation (
\ref{60GD}), with $\eta =0$. Then
\begin{equation}\label{15GD}
    \liminf_{\varepsilon\rightarrow 0}\Big(
    \varepsilon \ln \big(P\big\{
 \sqrt{\varepsilon}|B_s|^\alpha+\varepsilon^{\frac{2\alpha}{2-\alpha}}
    \geq \psi_0 (s)^\alpha  ; \ \forall s \in [0,1]\big\}\big)\Big)
    \geq -\frac{1}{2}\int_0^1\dot{\psi}_0^2(s)ds.
\end{equation}
\end{lemma}

\begin{prooff} \ {\bf of Lemma \ref{lem2GD}}.
We suppose for simplicity $\alpha=1$, the general case being only
slightly  more complicated.

\noi Let us introduce   the set :
$$
\Gamma_\varepsilon =\big\{
 \sqrt{\varepsilon}|B_s|
    \geq \psi_0 (s) -\varepsilon^2 ; \ \forall s \in [0,1]\big\}.
$$

\noi Since $\psi_0$ is an increasing and positive function, there
exists $\delta(\varepsilon)>0$ such that $\psi_0(s) \geq
\varepsilon^2$ if and only if $s\geq \delta(\varepsilon)$ and
$\delta(\varepsilon)$ goes to $0$ as $\varepsilon\rightarrow 0$.
Consequently :
$$
\Gamma_\varepsilon =\big\{
 \sqrt{\varepsilon}B_s
    \geq \psi_0 (s) -\varepsilon^2 ; \ \forall s \in [\delta(\varepsilon),1]\big\}
    \cup
    \big\{
 \sqrt{\varepsilon}B_s
    \leq -\psi_0 (s) +\varepsilon^2 ; \ \forall s \in
    [\delta(\varepsilon),1]\big\}.
$$
Then, computing the probability of $\Gamma_\varepsilon$, we obtain
:
$$
P(\Gamma_\varepsilon)=P\big\{B_s
    \geq \frac{\psi_0 (s)}{\sqrt{\varepsilon}} -\varepsilon^{3/2} ; \ \forall s \in
    [\delta(\varepsilon),1]\big\} ,
$$
$$
\qquad \geq P\big\{B_s -\frac{\psi_0 (s)}{\sqrt{\varepsilon}}
    \geq  -\varepsilon^{3/2} ; \ \forall s \in
    [0,1]\big\}.
$$
\noi Let us denote by $\displaystyle \Lambda_\varepsilon$ the set
: $\Lambda_\varepsilon= \{\inf _{s\in [0,1]} B_s\geq
-\varepsilon^{3/2}\}$. Using  Girsanov's theorem, we have :
$$
P(\Gamma_\varepsilon)\geq E \Big[1_{\Lambda_\varepsilon}\
\exp\Big\{ -\frac{1}{\sqrt{\varepsilon}}\int _0 ^1
\dot{\psi_0}(s)dB_s -\frac{1}{2\varepsilon}\int _0 ^1
\dot{\psi_0}(s)^2ds \Big\}\Big],
$$
$$
\qquad \qquad \geq \exp\Big\{-\frac{1}{2\varepsilon}\int _0 ^1
\dot{\psi_0}(s)^2ds \Big\} \ E \Big[1_{\Lambda_\varepsilon}\
\exp\Big\{ -\frac{1}{\sqrt{\varepsilon}}\int _0 ^1
\dot{\psi_0}(s)dB_s  \Big\}\Big].
$$
\noi Jensen's inequality applied to $x\mapsto e^{-x}$ leads to :
$$\frac{1}{P(\Lambda_\varepsilon)}E \Big[1_{\Lambda_\varepsilon}\
\exp\Big\{ -\frac{1}{\sqrt{\varepsilon}}\int _0 ^1
\dot{\psi_0}(s)dB_s  \Big\}\Big] \geq \exp\Big\{
-\frac{1}{P(\Lambda_\varepsilon)}E
\Big[1_{\Lambda_\varepsilon}\frac{1}{\sqrt{\varepsilon}}\int _0 ^1
\dot{\psi_0}(s)dB_s\Big]\Big\}.
$$
\noi Holder's inequality yields to :
$$E
\Big[1_{\Lambda_\varepsilon}\Big|\int _0 ^1
\dot{\psi_0}(s)dB_s\Big| \ \Big] \leq  C(p)
P(\Lambda_\varepsilon)^{1/p},
$$
\noi where $p>1$ and
$$C(p)=\Big(E \Big[\big|\int _0 ^1
\dot{\psi_0}(s)dB_s\big|^q\Big] \Big)^{1/q}=c_q\Big( \int _0 ^1
(\dot{\psi_0}(s))^2ds \Big)^{1/2} , \quad
\frac{1}{p}+\frac{1}{q}=1,
$$
\noi where $c_q=\Big(E[|B_1|^q]\Big)^{1/q}$ is a universal
constant.

\noi Recall that $\displaystyle{ (-\inf _{s\in [0,1]} B_s) }$ is
distributed as $|B_1|$, hence, there exist two positive constants
$c_0$ and $c_1$ such that :
$$c_0\varepsilon^{3/2}\leq P(\Lambda_\varepsilon)= P(|B_1|\leq
\varepsilon^{3/2})\leq c_1 \varepsilon^{3/2}, \quad \varepsilon
\in ]0,1].
$$
Let $0<\delta <1/2$, we choose $p>1$ such that $\displaystyle
1-{\frac{1}{p}=\frac{2}{3}\delta}$, then :
$$\frac{1}{P(\Lambda_\varepsilon)}E
\Big[1_{\Lambda_\varepsilon}\frac{1}{\sqrt{\varepsilon}}\int _0 ^1
\dot{\psi_0}(s)dB_s\Big] \leq \frac{c'_q\Big( \int _0 ^1
(\dot{\psi_0}(s))^2ds \Big)^{1/2}}{\varepsilon^{\delta+1/2}},
$$
$$
P(\Gamma_\varepsilon)\geq c_0\varepsilon^{3/2}
\exp\Big\{-\frac{1}{2\varepsilon}\int _0 ^1 \dot{\psi_0}(s)^2ds
\Big\} \  \exp\big\{ -\frac{c'_q\Big( \int _0 ^1
(\dot{\psi_0}(s))^2ds \Big)^{1/2}}{\varepsilon^{\delta+1/2}}
\big\},
$$
\noi where $\displaystyle c'_q=\frac{c_q}{c_0^{2\delta/3}}$. Then
(\ref{15GD}) follows immediately.

\end{prooff}

\begin{prooff} \ {\bf of Theorem \ref{theo1GD}}
\noi Suppose that $x=0$. Setting $\displaystyle
{\varepsilon=t^{\frac{\alpha-2}{2+\alpha}}}$ and using the scaling
property of Brownian motion and definition (\ref{2GD}), we have :
\begin{equation}\label{16GD}
    Z^{\lambda V}_t(0)=E_0\Big [ \exp \Big \{
    -\frac{\lambda t}{2}\int_0^1\frac{ds}{1+t^{\alpha/2}|B_s|^\alpha} \Big \}
\Big ]=E_0\Big [ \exp \Big \{-\frac{\lambda }{2\varepsilon}
\int_0^1\frac{ds}{\varepsilon^{\frac{2\alpha}{2-\alpha}}+|\sqrt{\varepsilon}B_s|^\alpha}
\Big \}\Big ] .
\end{equation}
%

1)We first prove :
\begin{equation}\label{15AGD}
    \limsup _{t\rightarrow +\infty}\Big (
t^{\frac{\alpha-2}{\alpha+2}}\ln \big(Z^{\lambda
V}_t(0)\big)\Big)\leq -\frac{1}{2}I_ 0(\lambda),
\end{equation}
\noi where $I_ 0(\lambda)$ is defined by (\ref{4GD}).

\noi Let $\eta>0$ be a fixed real number, and $\varepsilon>0$ such
that $\displaystyle{\varepsilon^{\frac{2\alpha}{2-\alpha}}<\eta}$.
Hence :
\begin{equation}\label{17GD}
    Z^{\lambda V}_t(0) \leq \exp \big
    \{-\frac{1}{2\varepsilon}\Phi_\eta
    (\lambda,\sqrt{\varepsilon}B_\cdot)\big\},
\end{equation}
\noi where
\begin{equation}\label{18GD}
\Phi_\eta (\lambda,f)=\lambda
\int_0^1\frac{ds}{\eta+|f(s)|^\alpha}.
\end{equation}
\noi  Varadhan's theorem \cite{DeuStro} yields to:
$$
\lim _{\varepsilon\rightarrow 0}\Big(\varepsilon \ln \Big(E_0\Big[
\exp \big
    \{-\frac{1}{2\varepsilon}\Phi_\eta
    (\lambda,\sqrt{\varepsilon}B_\cdot)\big\}\Big]\Big)\Big)=-\frac{1}{2}I_
    \eta (\lambda),
$$
\noi where $I_ \eta (\lambda)$ is defined by (\ref{5GD}).

\noi Consequently,
$$
\limsup _{t\rightarrow +\infty}\Big (
t^{\frac{\alpha-2}{\alpha+2}}\ln \big(Z^{\lambda
V}_t(0)\big)\Big)\leq -\frac{1}{2}I_ \eta (\lambda),
$$
\noi for any $\eta >0$.

\noi Lemma \ref{lem1GD}  implies (\ref{15AGD}).

2) We claim that :
\begin{equation}\label{20GD}
    \liminf _{t\rightarrow +\infty}\Big (
t^{\frac{\alpha-2}{\alpha+2}}\ln \big(Z^{\lambda V}_t(0)\big)\Big)
\geq -\frac{1}{2}I_ 0(\lambda),
\end{equation}
\noi Starting from (\ref{16GD}), we have :
$$
Z^{\lambda V}_t(x)\geq E_0\Big [1_{\Gamma_\eta} \exp \Big
\{-\frac{\lambda }{2\varepsilon}
\int_0^1\frac{ds}{\varepsilon^{\frac{2\alpha}{2-\alpha}}+|\sqrt{\varepsilon}B_s|^\alpha}
\Big \}\Big ],
$$
\noi with $\eta>0$ and
$$\Gamma_\eta=
\big\{
 \sqrt{\varepsilon}|B_s|^\alpha+\varepsilon^{\frac{2\alpha}{2-\alpha}}
    \geq \psi_0 (s)^\alpha  ; \ \forall s \in
    [0,1]\big\}\big)\Big).
    $$
    \noi Hence,
$$
Z^{\lambda V}_t(x)\geq \exp \Big \{-\frac{\lambda }{2\varepsilon}
\int_0^1\frac{ds}{\psi_0 (s)^\alpha}\Big\} P(\Gamma_\eta).
$$
\noi Relation (\ref{20GD}) is a direct consequence of Lemma
\ref{lem2GD} and (\ref{60GD}).
\end{prooff}




\section { The bilateral case}\label{bilateral}

\setcounter{equation}{0}

In this section it is required  that $V(x)$ goes to $+\infty$ as
$|x|\rightarrow \infty$. The asymptotic behaviour of
$\displaystyle Z^V_t(x):=E_x\Big [ \exp \big \{-{\frac{1}{2}} \int
_ 0^t V(B_s)ds \big \} \Big ]$, when $t \rightarrow \infty$ has
been initiated by Kac \cite{Kac} . Let us briefly recall (cf
\cite{Titch}) the main result.  Let us consider the second order
differential equation :
\begin{equation}\label{Kac1}
\frac{1}{2}\psi " -V\psi =-\lambda \psi .
\end{equation}
Then there exist a sequence $(\lambda _n)_{n\geq 1}$ of positive
numbers, $0<\lambda_1 <\lambda_2<\cdots$ and  an orthonormal basis
of functions $(\psi _n)_{n\geq 1}$ in $L^2(\mathbb R)$ such that
for any $n$,  $\psi _n$ is a solution to (\ref{Kac1}) with
$\lambda =\lambda _n$, and the ground state $\psi _1>0$. This
spectral gap property (i.e. $\lambda_1 >0 $, cf \cite{Titch})
plays a central role. With some additional assumptions, Kac proved
:
\begin{equation}\label{Kac2}
E_x\Big [ \exp \big \{-{\frac{1}{2}} \int _ 0^t V(B_s)ds \big \}
\Big ] \sim \rho e^{-\lambda_1 t} \quad t\rightarrow \infty,
\end{equation}
\noi where
$$\displaystyle \rho = \sum _{\lambda_j =\lambda_1}\psi _j
(0)\int_\mathbb R \psi _j(x)dx .$$
\noi  R. Carmona \cite{Carm1}, \cite{Carm2} generalized this
result to the case where $V$ may be written as the sum  $V_1+V_2$,
where $ V_2 \in L^p(\mathbb R),\ V_1$ being larger than a constant
and fulfilling for any $\beta>0$ :
$$\displaystyle
\lim _{|x|\rightarrow \infty} \Big ( \int _{x-\alpha}^{x +\alpha}
e^{-\beta V_1 (y)} dy =0 \Big), \mbox{ for some } \alpha >0.$$
\noi The proof is based on the compactness of the family of
operators $(T_t)_{t\geq 0}$, and the discrete spectrum of the
generator $L$ of the  semi group $(T_t)_{t\geq 0}$ :
\begin{equation}\label{Carmona1}
    L(f)=\frac{1}{2}f " -Vf .
\end{equation}
\noi Here, we investigate the case where
\begin{equation}\label{hypothesubil1}
 V : \mathbb{R} \rightarrow \mathbb{R} \mbox { is an even
 function, non-decreasing on } [0,+\infty[,
\end{equation}
\begin{equation}\label{hypothesubil2}
    V(x) \mbox { converges to a real number } \bar{V} , \mbox { as }
    |x|\rightarrow\infty .
\end{equation}
\noi Notice that, to our knowledge,  this setting   was neither
considered by Kac nor Carmona.

\noi Our approach is direct. We prove that there exists a unique
$\gamma_0$ such that a solution $\varphi_{V-\gamma_0}$ to the
Sturm-Liouville equation $ \varphi''=(V-\gamma_0)\varphi$ with
suitable boundary conditions,
 satisfies the condition of Section (\ref{Csmall}). This allows us
 to prove  the exponential decay of $Z^V_t(x)$, as
 $t\rightarrow\infty$.

 \noi  In order to present our main result in Theorem \ref{thbasebil} below,we need
 to define  the parameter $\gamma_0$.

\noi We start with the following definitions :
\begin{equation}\label{borninfubil}
   \underline{V}=\inf_{x\in \mathbb{R}}V(x)=V(0) \ ,\
 \bar{V}= \sup_{x\in \mathbb{R}}V(x) <\infty .
\end{equation}
\noi Notice that we do not require that $V(x)$ is non-negative.

\noi If $V$ is constant, the result is obvious.  Therefore we
suppose in the sequel that : $\underline{V} <\bar{V}$.

\noi The function $V$ being monotone on  $[0,\infty[$ (resp.
$]-\infty,0]$),   coincides with its right continuous modification
 $V_0$, except on an at most countable set. Then,   a.s. :
$$\int_0^t V(B_s)ds=\int_0^t V_0(B_s)ds, \quad \mbox{ for any }
t\geq 0.$$
Then $V$ may be assumed to be right continuous.

\noi Let $V^{-1}$ be the right inverse of the restriction of $V$
 to $[0,\infty[$ :
$$V^{-1}(\gamma) =\inf \{t\geq 0 ; V(t)>\gamma \}, \quad
 \gamma \in ]\underline{V} , \bar{V} [.$$
\noi Then:
\begin{equation}\label{propu-1bil}
    V(V^{-1}(\gamma))\geq \gamma, \quad \underline{V} < \gamma
    < \bar{V}.
\end{equation}
 \noi For any $\gamma \in ]\underline{V} , \bar{V} [$, let
$F_\gamma$ (resp. $G_\gamma$) be the unique solution to
\begin{equation}\label{Sturmgamma1bil}
    Y" =(V-\gamma)Y,
\end{equation}
\noi on $[0,V^{-1}(\gamma)]$ (resp. $[V^{-1}(\gamma) ,\infty [$)
with the  boundary conditions
\begin{equation}\label{Fgamma2bil}
F_\gamma (V^{-1}(\gamma))=1,  \qquad F'_\gamma (0)=0.
\end{equation}
    \begin{equation}\label{Ggamma3bil}
    \displaystyle
G_\gamma (V^{-1}(\gamma))=1,  \qquad G_\gamma (+\infty):= \lim_{x
\rightarrow \infty}G_\gamma (x)=0.
\end{equation}
\noi We set:
\begin{equation}\label{Cgamma4bil}
\varphi_{V-\gamma }(x) = \left \{
\begin{array}{lll}
F_\gamma (x) & \mbox { if } & x\in [0,V^{-1}(\gamma)]\\
G_\gamma (x) & \mbox { if } & x \geq V^{-1}(\gamma).
\end{array}
\right.
\end{equation}

\noi We extend $\varphi_{V-\gamma } $ to the whole line, setting :
$\varphi_{V-\gamma } (-x)=\varphi_{V-\gamma }(x)$. Then
$\varphi_{V-\gamma } $ is a continuous and even function defined
on $\mathbb{R}$. Notice that $\varphi_{V-\gamma } (x)$ is
differentiable for any $x\not = \pm V^{-1}(\gamma)$.


\vskip 10 pt \noi
\begin{theorem} \label{thbasebil}
Let $V$ be a function  fulfilling (\ref{hypothesubil1}) and
(\ref{hypothesubil2}).
\begin{enumerate}
 \item
 There exists a unique   $\gamma_0 \in ]\underline{V},\bar{V}[$ such that the
 function $\varphi_{V-\gamma_0 } $
 defined by (\ref{Cgamma4bil}) is
 differentiable on $\mathbb{R}$.

\item
 The quantity:
$$e^{ \gamma_0 t/2} \ E _x\Big [ \exp \Big \{-{\frac{1}{2}}
\int _0^t V(B_s)ds \Big \} \Big ],
$$
 converges as $t\rightarrow \infty$, to $\displaystyle \big (
 \int_\mathbb{R}\varphi_{V-\gamma_0 }(y)dy \big ) \varphi_{V-\gamma_0 }(x)$.
\item
 Let us define the probability $Q _{x,t}^V$ on ${\cal F} _t$ via :
\begin{equation} \label{def0Qbil}
Q _{x,t}^V(\Lambda  _t) = \frac{E _x \Big [1 _{\Lambda  _t}\ \exp
\Big \{-\frac{1}{2}
 \int _0^t V(B_h)dh \Big \} \Big ]}
 {E _x \Big [ \exp \Big \{-\frac{1}{2} \int _0^t V(B_h)dh\Big \} \Big ]}
  ,\  \Lambda  _t \in {\cal F} _t.
\end{equation}
\noi Then, for any $\Lambda  _s$ in ${\cal F} _s$, $Q
_{x,t}^V(\Lambda  _s)$ converges to $ P _{x}^{\varphi_{V-\gamma_0
}}(\Lambda _s)$
 as $t \rightarrow \infty$, where $P _{x}^{\varphi_{V-\gamma_0
}}$ is the probability
 defined on ${\cal F} _\infty$ by :
\begin{equation} \label{defQbil}
P _{x}^{\varphi_{V-\gamma_0 }}(\Lambda  _s )=\frac{e^{\gamma_0
s/2}}{\varphi_{V-\gamma_0 }(x)}  E _x\Big [1 _{\Lambda _s}
\varphi_{V-\gamma_0 }(B_ s) \exp \Big \{-{\frac{1}{2}} \int _0^s
V(B_h)dh \Big \} \Big ],
\end{equation}
\noi for any $s>0$ and $\Lambda  _s \in {\cal F} _s$.
 \item
 Let  $(X_t ^x;t\geq 0)$ be the solution to :
\begin{equation} \label{defXbil}
X _t=x+B _ t+\int  _0 ^t \frac{\varphi'_{V-\gamma_0 }}{
\varphi_{V-\gamma_0 }}(X _s)ds, \ t \geq 0.
\end{equation}
 Then the law of $(X_t ^x;t\geq 0)$ is $P _{x}^{\varphi_{V-\gamma_0
}}$.

\item
 The process $(X_t ^x;t\geq 0)$ is recurrent with invariant
 finite measure   $  (\varphi_{V-\gamma_0
})^2(x)dx$.

 \noi

\end{enumerate}
\end{theorem}

 \vskip 20 pt \noi We begin by proving two
preliminary results in the form of the next  Lemmas \ref{lem1bil}
and \ref{lem2bil}.


\begin {lemma} \label{lem1bil}
Let $\theta_1 , \theta_2$ be two  functions defined on $[a,b), \
\theta_1 \geq \theta_2 \geq 0, \ \varphi_i$ a solution of
$\varphi_i''=\theta_i\varphi_i,i=1,2$ on $[a,b)$, such that
$\varphi_1(a)=\varphi_2(a)$ and $\varphi_1(b)=\varphi_2(b)\geq 0$
(If $b=+\infty$, $\varphi_i(b)$  has to be understood as
$\displaystyle \lim _{ x\rightarrow \infty} \varphi_i(x)$). Then
$\varphi_2 \geq \varphi_1$.
\end {lemma}


\begin {prooff} \ {\bf of Lemma \ref{lem1bil}}.
Suppose there exists $x_0 $ in $(a,b)$ such that $\varphi_2 (x_0)
< \varphi_1 (x_0)$. We can find a non-empty interval $[\alpha
,\beta[$ included in $[a,b)$ such that
$\varphi_1(\alpha)=\varphi_2(\alpha), \ \varphi_1(\beta)
=\varphi_2(\beta)$ and $\varphi_1 >\varphi_2$ on $]\alpha
,\beta[$.

\noi Let $h=\varphi_1 -\varphi_2$. Then
$h''=\theta_1\varphi_1-\theta_2\varphi_2$.

\noi But on $[\alpha ,\beta [$:
$$ \theta_1\geq \theta_2 \geq 0, \quad \varphi_1\geq \varphi_2 \geq 0, \quad
\Rightarrow h''\geq 0.$$
\noi This generates a contradiction because then $h$ is a
non-constant and non-negative convex function on $[\alpha ,\beta[
$ such that  $h(\alpha)=h(\beta)=0$.
\end {prooff}

\begin {lemma} \label{lem2bil}
Let $\gamma \in ]\underline{V} ,\bar{V}[$.
\begin{enumerate}
    \item
    There exists two positive constants $k_1,k_2$ such that
    \begin{equation}\label{ineg1bil}
    \varphi_{V-\gamma}(x)\leq k_1e^{-k_2|x|}.
\end{equation}
\item The function : $\gamma :\in ]\underline{V},\bar{V}[ \
\mapsto \varphi_{V-\gamma} '(V^{-1}(\gamma)_-)$ is continuous,
increasing and
\begin{equation}\label{ineg2bil}
\lim_{\gamma\rightarrow \underline{V}}
\varphi_{V-\gamma}'(V^{-1}(\gamma)_-) =0, \qquad
\liminf_{\gamma\rightarrow \bar{V}} \big ( -\varphi_{V-\gamma}
'(V^{-1}(\gamma)_-)\big ) >0
\end{equation}
\noi (  $\varphi_{V-\gamma}' (V^{-1}(\gamma) _-)$ denotes the left
derivative of $\varphi_{V-\gamma}'$ at point  $V^{-1}(\gamma)$).

\item We have :
\begin{equation}\label{ineg3bil}
\lim_{\gamma\rightarrow \bar{V}} \varphi_{V-\gamma}
'(V^{-1}(\gamma)_+) =0, \quad \lim_{\gamma\rightarrow
\underline{V}} \varphi_{V-\gamma} '(V^{-1}(\gamma)_+) <0 .
\end{equation}
\end{enumerate}
\end {lemma}


\begin {prooff} \ {\bf of Lemma \ref{lem2bil}}.
 1) Let $\underline{u}<\gamma
<\gamma'<\bar{u}$ and $\theta$ be the solution to
$\theta''=(\gamma'-\gamma)\theta$, on $[V^{-1}(\gamma
'),+\infty[$, with the boundary conditions : $\theta
(V^{-1}(\gamma') )=\varphi_{V-\gamma} (V^{-1}(\gamma')),\  \theta
(+\infty)= \varphi_{V-\gamma} (+\infty)=0$.

\noi Because $0< \gamma' -\gamma \leq V(x) -\gamma$, for any
$x\geq V^{-1}(\gamma')$, Lemma \ref{lem1bil} implies
$\varphi_{V-\gamma} \leq \theta $.

\noi But
$$\theta (x) = \varphi_{V-\gamma} (V^{-1}(\gamma'))
e^{-k_1 (x - V^{-1}(\gamma'))},$$
\noi with $k_1=\sqrt{\gamma'-\gamma}$.

\noi This proves (\ref{ineg1bil}).

\noi 2) Since $\varphi_{V-\gamma}$ coincides with $F_\gamma$ on
$[0,V^{-1}(\gamma)]$,
$$\varphi_{V-\gamma}'(V^{-1}(\gamma) _-)=\int_0^{V^{-1}(\gamma) )}
(V(x)-\gamma ) \varphi_{V-\gamma} (x) dx.$$
\noi The function $\varphi_{V-\gamma} $ is concave on
$[0,V^{-1}(\gamma)]$, consequently, if $x \in [0,V^{-1}(\gamma)]$,
then
\begin{equation}\label{ineg4bil}
\varphi_{V-\gamma}  (x) \geq 1 ,
\end{equation}
 \begin{equation}\label{ident3bil}
|\varphi_{V-\gamma}' (V^{-1}(\gamma)_- )|=-\varphi_{V-\gamma}'
(V^{-1}(\gamma) _-) ,
\end{equation}
\begin{equation}\label{ident4bil}
|\varphi_{V-\gamma}' (V^{-1}(\gamma)_- )| =\int_0^{V^{-1}(\gamma)
)}(\gamma -V(x))\varphi_{V-\gamma} (x) dx ,
\end{equation}
\begin{equation}\label{ident5bil}
|\varphi_{V-\gamma}'(V^{-1}(\gamma)_- )|\geq
\int_0^{V^{-1}(\gamma) )}(\gamma -V(x))dx,
\end{equation}
$$\varphi_{V-\gamma} (x) \leq 1 +\varphi_{V-\gamma}'(V^{-1}(\gamma) _-)
(x -V^{-1}(\gamma)), \qquad \qquad\qquad\qquad\qquad\qquad
\qquad$$
$$\varphi_{V-\gamma} (x) \leq 1 +|\varphi_{V-\gamma}'(V^{-1}(\gamma)_- )|
 (V^{-1}(\gamma)-x).
\qquad\qquad\qquad\qquad\qquad\qquad$$
\noi Minoring $\varphi_{V-\gamma} (x)$ in (\ref{ident4bil}), we
obtain:
\begin{equation}\label{ineg5bil}
|\varphi_{V-\gamma} (V^{-1}(\gamma) _-)|\Big (
1-\int_0^{V^{-1}(\gamma) )}(\gamma -V(x))(V^{-1}(\gamma) -x)dx
\Big )\leq \int_0^{V^{-1}(\gamma) )}(\gamma -V(x))dx .
\end{equation}
\noi From now on, we suppose for simplicity that the restriction
of $V$
 to $[0,+\infty[$ is strictly increasing. In particular $V^{-1}$
is a continuous function. Relation (\ref{ineg5bil}) implies that
$$\lim_{\gamma\rightarrow \underline{V}} \varphi_{V-\gamma}
'(V^{-1}(\gamma)_-) =0.$$
\noi From (\ref {ident5bil}), it easily follows that
$$\liminf_{\gamma\rightarrow \bar{V}}
\varphi_{V-\gamma} '(V^{-1}(\gamma)_-) >0 .$$

3) Let $h$ be  the solution to $h''=(\bar{V}-\gamma)h$, on
$[V^{-1}(\gamma),+\infty[$,  with the boundary conditions : $h
(V^{-1}(\gamma) )=\varphi_{V-\gamma} (V^{-1}(\gamma)),\  h
(+\infty)= \varphi_{V-\gamma} (+\infty)=0$.

\noi Since $0< V -\gamma \leq \bar{V} -\gamma$, for any $x\geq
V^{-1}(\gamma)$,  Lemma \ref{lem1bil} shows that
$\varphi_{V-\gamma}\geq h $.

\noi But $h (V^{-1}(\gamma) )=\varphi_{V-\gamma}
(V^{-1}(\gamma))$, hence
$$|\varphi_{V-\gamma}' (V^{-1}(\gamma)_+)|=-\varphi_{V-\gamma}' (V^{-1}(\gamma)_+) \leq -
h'(V^{-1}(\gamma))=\sqrt{\bar{V}-\gamma},$$
$$\lim_{\gamma\rightarrow \bar{V}} \varphi_{V-\gamma} '(V^{-1}(\gamma)_+)
=0.$$
\noi It is easy to check that :
$$\lim_{\gamma\rightarrow \underline{V}} \varphi_{V-\gamma}
'(V^{-1}(\gamma)_+) =H'(0),$$
\noi where $H$ is the solution to $H"(x)=(V(x)-\gamma )H(x)\ x\geq
0$, with the boundary conditions : $H(0)=1$ and $H(\infty )=0$.

\noi Because $H'(\infty)=0$,
$$-H'(0)=\int _0^\infty(V(y)-\gamma)H(y)dy >0.$$
This gives (\ref{ineg3bil}).
\end {prooff}


\vskip 10 pt \noi  \begin {prooff} \ {\bf of   Theorem
\ref{thbasebil}}.

\noi The existence of $\gamma _0$ such that $\varphi_{V-\gamma_0}$
is of class $C^1$ can be derived using  the continuity of the
functions $\gamma \mapsto \varphi_{V-\gamma}' (V^{-1}(\gamma)_+)$
and $\gamma \mapsto \varphi_{V-\gamma}' (V^{-1}(\gamma)_-)$, and
properties (\ref{ineg2bil}), (\ref{ineg3bil}).

\noi $\varphi_{V-\gamma_0}$ is an even function, non-decreasing on
$(0,\infty[$. Inequality (\ref{ineg1bil}) leads to $\displaystyle
\int _\mathbb R (\varphi_{V-\gamma_0}(x))^p dx<\infty$, for any
$p>0$. Then we may apply Theorem \ref {theoC} :
\begin{equation}\label{ident8bil}
\lim _{t\rightarrow\infty}E _x\Big [ \exp \Big \{-{\frac{1}{2}}
\int _0^t ( V(B_s)-\gamma_0)ds \Big \} \Big ]=
\frac{\int_\mathbb{R}\varphi_{V-\gamma_0}(y)dy}
{\int_\mathbb{R}(\varphi_{V-\gamma_0})^2(y)}
  \varphi_{V-\gamma_0}(x).
\end{equation}
\noi Consequently point 2 of Theorem \ref {thbasebil} holds.

\noi Obviously the probability measure $Q _{x,t}^V$ defined by
relation (\ref{def0Qbil}) is also given by the following :
$$Q _{x,t}^V(\Lambda  _t) = \frac{E _x \Big [1 _{\Lambda  _t}\
\exp \Big \{-\frac{1}{2}
 \int _0^t (V(B_h)-\gamma_0)dh \Big \} \Big ]}
 {E _x \Big [ \exp \Big \{-\frac{1}{2} \int _0^t (V(B_h)-\gamma_0)dh\Big \} \Big ]}
  ,\  \Lambda  _t \in {\cal F} _t.
$$
From Theorem \ref {theoC} we may  conclude that $Q
_{x,t}^V(\Lambda  _s)$ converges to $Q _x^V(\Lambda  _s )$, for
any positive $s$ and $\Lambda  _s$ in ${\cal F}_s$, and $P
_x^{\varphi_{V-\gamma_0}}$ is defined by (\ref{defQbil}).

\noi Parts 3. and 4. of Theorem \ref{thbasebil} are direct
consequences of
 Theorem \ref {theoC}.

\end{prooff}



\begin {exam} \label {ex1bil}
We will denote by    $V$  the function : $V(x)=1_{\{|x|>a\}}$,
where $a>0$. Let $\gamma _0$ be the unique solution in $]0,1
\wedge (\frac{\pi^2}{4})[$ to
$$ \tan (a \sqrt{\gamma})=\sqrt{\frac{1-\gamma}{\gamma}}.$$
\noi Then  :
\begin{equation}\label{defCbil}
    \varphi_{V-\gamma _0}(x)= \left \{
    \begin {array}{ll}
 \displaystyle e^{-\sqrt{1-\gamma_0}|x-a|} & \mbox {if } |x|>a \\
\displaystyle \frac{\cos (\sqrt{\gamma_0} x) }{\cos
(\sqrt{\gamma_0} a)} & \mbox {if } |x|\leq a .
\end {array}
\right.
\end{equation}
\end {exam}

%
\def\refname{References}
\bibliographystyle{plain}
\bibliography{refexpoRVYoct03}
\end{document}